\documentclass[notitlepage,11pt, a4paper]{article}
\usepackage[utf8]{inputenc}

\usepackage{amsfonts, amssymb, amsmath}
\usepackage{mathtools}
\usepackage{fullwidth}
\usepackage{graphicx}

\usepackage{caption}
\usepackage{subcaption}

\usepackage{algorithm2e}

\usepackage[colorlinks=true, pdfstartview=FitV, linkcolor=blue, citecolor=blue, urlcolor=blue,pagebackref=false]{hyperref}
\usepackage{varioref}
\usepackage[capitalize]{cleveref}

\usepackage[title]{appendix}

\usepackage[numbers]{natbib}
\usepackage{amsthm}
\usepackage{thmtools, thm-restate}

\declaretheorem[numberwithin=section]{theorem}
\declaretheorem[sibling=theorem]{lemma}
\declaretheorem[sibling=theorem]{corollary}
\declaretheorem[sibling=theorem]{proposition}
\declaretheorem[sibling=theorem]{remark}
\declaretheorem[sibling=theorem]{definition}

\usepackage{xcolor}
\definecolor{darkred}{rgb}{0.9,0.1,0.1}

\def \D {{\mathbb D}}
\def \E {{\mathbb E}}
\def \N {{\mathbb N}}
\def \P {{\mathbb P}}

\def \R {{\mathbb R}}
\def \Z {{\mathbb Z}}

\def \cC {{\mathcal C}}

\def \cE {{\mathcal E}}

\def \cG {{\mathcal G}}

\def \cI {{\mathcal I}}

\def \cL {{\mathcal L}}
\def \cM {{\mathcal M}}

\def \cT {{\mathcal T}}

\def \vD {{\mathbf D}}
\def \vX {{\mathbf X}}
\def \vZ {{\mathbf Z}}
\def \va {{\mathbf a}}
\def \vc {{\mathbf c}}
\def \vd {{\mathbf d}}

\def \vk {{\mathbf k}}
\def \vs {{\mathbf s}}
\def \vu {{\mathbf u}}
\def \vx {{\mathbf x}}
\def \vy {{\mathbf y}}

\def \subset {\subseteq}

\DeclarePairedDelimiter{\abs}{\lvert}{\rvert}
\DeclarePairedDelimiter{\norm}{\lVert}{\rVert}
\DeclarePairedDelimiter{\floor}{\lfloor}{\rfloor}

\DeclareMathOperator{\var}{Var}
\DeclareMathOperator{\cov}{Cov}
\DeclareMathOperator{\corr}{Corr}
\DeclareMathOperator{\iso}{Iso}
\DeclareMathOperator{\supp}{supp}

\def\one{\rlap{\mbox{\small\rm 1}}\kern.15em 1} % indicator function

\def \todist {\xrightarrow[]{\mathrm{(d)}}}

\usepackage{stmaryrd} % needed for double square brackets
\DeclarePairedDelimiter{\pathbtw}{\llbracket}{\rrbracket}
\def \littleo {o}
\def \littleomega {\omega}
\def \bigo {O}
\newcommand \diff[2] {\frac{\mathrm{d}#1}{\mathrm{d}#2}}
\newcommand*\dif{\mathop{}\!\mathrm{d}}
\def \lattice {\Lambda}
\newcommand \vectwo[2] {\begin{pmatrix} #1 \\ #2 \end{pmatrix}}
\def \maxeval {{\lambda_{\text{max}}}}
\def \mineval {{\lambda_{\text{min}}}}
\def \centeredX {\hat{\vX}}
\def \distmdm {d_{\vec{G}}}

\def \zeroloop {\mathfrak L}

\def \lltEvent {\mathcal A}
\def \omegaEvent {\mathcal B}

\def \walkDeviationEvent {\mathcal E}

\def \toprob {\xrightarrow[]{\mathrm{(p)}}}

\title{Universality for the directed configuration model: metric space convergence of the strongly connected components at criticality}
\author{ \and Zheneng Xie}

\author{Serte Donderwinkel\thanks{\href{mailto:serte.donderwinkel@stats.ox.ac.uk}{serte.donderwinkel@stats.ox.ac.uk}, University of Oxford, Oxford, United Kingdom. ORCID iD: 0000-0001-8148-8631.} \and Zheneng Xie \thanks{\href{zheneng.xie@stats.ox.ac.uk}{zheneng.xie@stats.ox.ac.uk}, University of Oxford, Oxford, United Kingdom. ORCID iD: 0000-0002-7252-3631. This research has been supported by the EPSRC Centre for Doctoral Training in Mathematics of Random Systems: Analysis, Modelling and Simulation (EP/S023925/1).}}
\date{\today}

\begin{document}

\maketitle
\begin{abstract}
  We consider the strongly connected components (SCCs) of a uniform directed graph on $n$ vertices with i.i.d. in- and out-degree pairs distributed as $(D^-,D^+)$, with $\mathbb E[D^+]=\mathbb E[D^-]=\mu$. We condition on equal total in- and out-degree. A phase transition for the emergence of a giant SCC is known to occur at the critical value $\mathbb E[D^-D^+] = \mu$. We study the model at this critical value and, additionally, require that $\mathbb E[(D^-)^i(D^+)^j]<\infty$ for all $i+j\leq 3$, and for $(i,j)=(1,3)$ and $(i,j)=(3,1)$. We show that, under these conditions, the SCCs ranked by decreasing number of edges with distances rescaled by $n^{-1/3}$ converge in distribution to a sequence of finite strongly connected directed multigraphs with edge lengths, and that these are either $3$-regular or loops. The limit objects lie in a $3$-parameter family, which contains the scaling limit of the SCCs in the directed Erd\H{o}s-R\'enyi model at criticality as found by Goldschmidt and Stephenson (2019). This is the first universality result for the scaling limit of a critical directed graph model and the first quantitative result on the directed configuration model at criticality. As a trivial consequence, the largest SCCs at criticality contain $\Theta(n^{1/3})$ vertices and edges in probability, and the diameter of the directed graph at criticality is $\Omega(n^{1/3})$ in probability. We use a metric on the space of weighted multigraphs in which two multigraphs are close if there are compatible isomorphisms between their vertex and edge sets which roughly preserve the edge lengths. We use the product topology on the sequence of multigraphs. Our method of proof involves a depth-first exploration of the directed graph, resulting in a spanning forest with additional identifications, of which we study the limit under rescaling.
\end{abstract}
% \tableofcontents

\section{Introduction}

\subsection{Overview}

Edges in real-world networks are often directed, such as links on the world wide web, ``follows'' on Twitter, financial transactions or disease transmission in a social network. When analysing networks, the first quantity that is often considered is the distribution of the degrees of nodes in the network.  In this paper we will consider sampling an i.i.d.\ sequence of in- and out-degrees, conditional on the total in-degree being equal to the total out-degree. We will then sample a uniform directed graph (digraph) with the given degree sequence. Results on such graphs are a useful benchmark, exposing additional underlying structure of a real-world network compared to a uniformly random graph with its degree sequence.

\begin{figure}
    \centering
    \includegraphics[scale=0.6]{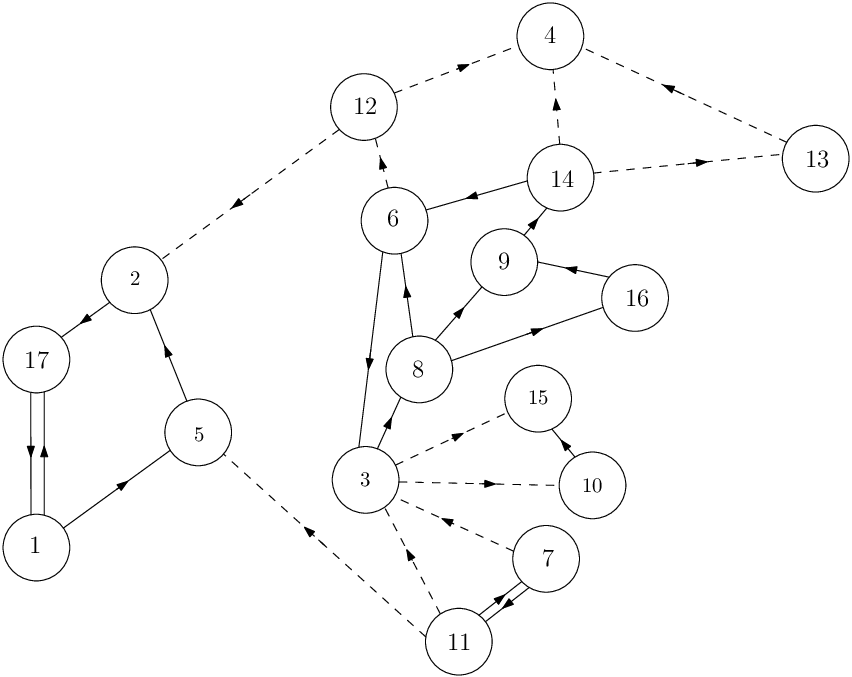}
    \caption{A directed graph on [17]. The strongly connected components have vertex sets $\{1,2,5,17\}$, $\{3,6,8,9,14,16\}$, $\{7,11\}$, $\{4\}$, $\{10\}$, $\{12\}$, $\{13\}$, and $\{15\}$. Edges that are not part of an SCC are depicted as dashed arrows. Taken from \cite{goldschmidtScalingLimitCritical2019} with permission of the authors.}
    \label{fig.SCCs}
\end{figure}

When considering such models, previous work by \citet{cooperSizeLargestStrongly2004} (which we will discuss in more detail in Section \ref{sec.previouswork}) shows that there exists a phase transition in the strong directed connectivity of the graph. Two vertices are part of the same \emph{strongly connected component} (SCC) if and only if there exists a directed cycle that contains both of them. Above some threshold, there will exist a unique giant SCC that occupies a positive proportion of the vertices whereas, below the threshold no SCC will occupy a positive proportion of the vertices. In Figure \ref{fig.SCCs}, a directed graph and its strongly connected components are depicted. In this paper we will prove the first detailed results about the critical case - specifically, that there exists a sequence of random weighted directed multigraphs that can be understood as the scaling limit of the SCCs when viewed in decreasing order of size.

\subsection{Directed graphs}

\begin{figure}
    \centering
    \includegraphics{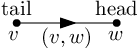}
    \caption{An edge $(v,w)$ will be depicted as an arrow from $v$ to $w$.}\label{fig.tailhead}
\end{figure}

There are two notions of connectivity when working with a directed graph: weak and strong connectivity. We will be working with the strong notation. We say a vertex $v$ leads to a vertex $w$, written $v \rightarrow w$, if there exists a directed path from $v$ to $w$ in the graph. We say $v$ is \emph{strongly connected to $w$}, written $v \leftrightarrow w$, if $v$ leads to $w$ and $w$ leads to $v$. By convention, $v$ leads to itself. A graph is \emph{strongly connected} if all pairs of vertices in the graph are strongly connected. The relation $v \leftrightarrow w$ is an equivalence relation; the digraphs induced by the equivalence classes of $\leftrightarrow$ are referred to as the \emph{strongly connected components} (SCCs). For each vertex $v$ in a directed graph $\vec{G}$, we will use the notation $d^-(v)$ for the in-degree of $v$ and $d^+(v)$ for the out-degree of $v$. Moreover, a directed edge $(v,w)$ has \emph{tail} $v$ and \emph{head} $w$ (see Figure \ref{fig.tailhead}).

\subsection{Description of the model}

\label{subsec:model-description}

First consider a deterministic degree sequence $\vd_1, \ldots, \vd_n$ where $\vd_i = (d_i^-, d_i^+)\in \N \times \N$ for $i = 1, \ldots, n$. We say a directed graph with vertex set $[n]$, where $[n] = \{1,\dots,n\}$, has degree sequence $\vd_1, \ldots, \vd_n$ if $(d^-(i), d^+(i)) = (d_i^-, d_i^+)$ for $i = 1, \ldots, n$.

In order to sample a uniformly random graph with a given degree sequence, we first consider the \emph{directed configuration model} introduced by \citet{cooperSizeLargestStrongly2004}. Take $n$ vertices $v_1, \ldots, v_n$ such that $v_i$ has $d^-_i$ in-half-edges and $d^+_i$ out-half-edges. Then construct a multigraph by choosing a uniformly random pairing of the in-half-edges with the out-half-edges. \citet[Sec.\ 2.1]{cooperSizeLargestStrongly2004} proved that if we condition on the resulting multigraph being simple, we obtain a uniformly chosen random digraph with the given degree sequence.

In this paper we will consider the case where the degree sequence  consists of $n$ i.i.d.\ random variables conditioned on the total in-degree being equal to the total out-degree. Let $\nu$ be a distribution on $\N \times \N$, and let $\vD_1, \ldots, \vD_n$ be a sequence of i.i.d.\ random variables with distribution $\nu$. We condition on the event
\begin{equation*}
    \left\{ \textstyle \sum_{i=1}^n D_i^- = \sum_{i=1}^n D_i^+ \right\},
\end{equation*}
observing that this is an asymptotically singular event as $n\to\infty$. Let $\vec{G}_n(\nu)$ be a digraph chosen uniformly at random from all digraphs with degree sequence $\vD_1, \ldots, \vD_n$. We are interested in the limit under rescaling of the SCCs of $\vec{G}_n(\nu)$ as $n\to \infty$.

Suppose $(D^-, D^+)$ has law $\nu$. We will require the following assumptions to hold:
\begin{enumerate}
    % \item $\E[(D^-) + D^+)^3] < \infty$,
    \item $\E[(D^-)^i(D^+)^j]< \infty$ for $1 \leq i+j\leq 3$, $(i, j) = (1, 3)$ and $(i, j) = (3, 1)$.
    \item $\E[D^-] = \E[D^+]$.
    \item $D^- - D^+$ is strongly aperiodic. This means that for all $p > 1$, there does not exist $k \in \Z$ such that 
    \begin{equation*}
        \P(D^- - D^+ \in k + p\Z) = 1.
    \end{equation*}
    \item $\E[D^-D^+] = \E[D^{\pm}]$.
\end{enumerate}

The first condition is required to ensure that the steps of a random walk used in the proof have finite variance, so that the random walk will convergence under rescaling to a Brownian motion. It also ensures similar regularity of other random variables that we use to encode the directed graph. (We discuss relaxing the moment conditions in Subsection \ref{subsec.openproblems}.)

The second and third conditions make sure the event $\{\sum_{i=1}^n D^-_i = \sum_{i=1}^n D^+_i\}$ is well-behaved. The second condition ensures that it is not a large deviation event. Using a result from \citet[Page 42, P1]{spitzerPrinciplesRandomWalk1964}, the third condition ensures that the event has positive probability for all sufficiently large $n \geq 1$. This condition can be relaxed to assuming that $D^- - D^+$ is non-constant by taking limits for $n \in p \N$ rather than $n \in \N$ where $p$ is the periodicity of $D^- - D^+$. However, for simplicity of presentation, we will keep it as an assumption.

The fourth assumption is the criticality condition. To understand how this arises, consider the directed configuration model and let $(V_n, W_n)$ be a uniformly chosen edge. For now, ignore the conditioning on the total in- and out-degrees being equal. We consider the distribution of the in- and out-degree of $W_n$. Because the degree sequence is an i.i.d.\ sequence, $W_n$ is equally likely to be any vertex $i$. Thus for any $\vk = (k^-, k^+)$,
\begin{align*}
    \P(d^-(W_n) = k^-, d^+(W_n) = k^+)
    &= n \P(W_n = 1, \vD_1 = \vk) \\
    &= n \E[\P(W_n = 1 \mid \vD_1 = \vk, \vD_2, \ldots, \vD_n)] \P(\vD_1 = \vk)
\end{align*}
Conditionally on the degree sequence, we have that $W_n = i$ with probability proportional to $D^-_i$ since we used an uniform pairing of the in- and out-half-edges. Therefore
\begin{align*}
    \P(W_n = 1 \mid \vD_1 = \vk, \vD_2, \ldots, \vD_n)
    &= \frac{k^-}{k^- + \sum_{i=2}^n D_i^-}.
\end{align*}
Thus
\begin{equation*}
    \P(d^-(W_n) = k^-, d^+(W_n) = k^+) = \E\left[ 
        \frac{k^-}{\frac{1}{n}\left( k^- + \sum_{i=2}^n D_i^- \right)}
    \right]
    \P \left[ D^- = k^-, D^+ = k^+ \right].
\end{equation*}
Using the law of large numbers, the above will converge to
\begin{equation*}
    \frac{k^-}{\E[D^-]} \P\left[ D^- = k^-, D^+ = k^+ \right].
\end{equation*}
Let $(Z^-, Z^+)$ be such that $P(Z^- = k^-, Z^+ = k^+)$ is given by the above expression. We say $(Z^-, Z^+)$ has the law of the \emph{degree distribution size-biased by in-degree}. For large $n$, any other fixed out-edge of $W_n$ is then also distributed approximately like a uniformly chosen edge (here we are ignoring the fact that we have already sampled an edge) since we chose the in- and out-edge pairing uniformly at random. Therefore the out-degree of the head will have approximately the same distribution as $Z^+$. Thus if we were to look at the graph of all vertices leading from $W_n$, it would look approximately like a Bienaymé tree\footnote{For $\mu$ a probability distribution on $\N$, a Bienaymé tree with offspring distribution $\mu$ is the family tree of a branching process with offspring distribution $\mu$. Bienaymé trees are often referred to as Galton-Watson trees, but we decide to follow the name change suggested by \citet{addario-berryUniversalHeightWidth2021}.} with offspring distribution $Z^+$. It is well known that such trees exhibit critical behaviour in whether or not the tree is finite at $\E[Z^+] = 1$. This is equivalent to assuming $\E[D^-D^+] = E[D^-]$.

\citet{cooperSizeLargestStrongly2004} studied this phase transition for a deterministic degree sequence $\vd_1, \ldots, \vd_n$. They defined the parameter
\begin{equation*}
    d = \frac{\sum_{i=1}^n d_i^+ d_i^-}{\sum_{i=1}^n d_i^-}
\end{equation*}
which is a counterpart of $\E[Z^-]$ for deterministic degree sequences. They then showed that, under additional assumptions, there exists a phase transition for the existence of a giant SCC depending on whether $d$ is strictly greater than or less than 1. Our work in this paper shows our corresponding condition, $\E[Z^-] = 1$, is also the correct criticality condition to take for i.i.d.\ random degree sequences.

We define the following parameters that will determine the behaviour of the SCCs in the limit.
\begin{enumerate}
    \item $\mu:=\E[D^-]=\E[D^+]=\E[D^-D^+]$
    \item $\nu_-:= \E[Z^-] - 1 = \frac{\E[(D^-)^2]-\mu}{\mu}$ 
    \item $\sigma_-^2 := \var(Z^-) = \frac{\mu\E[(D^-)^3]-\E[(D^-)^2]^2}{\mu^2}$ 
    \item $\sigma_+^2 := \var(Z^+) = \frac{\E[D^-(D^+)^2]-\mu}{\mu}$ 
    \item $\sigma_{-+} := \cov(Z^-, Z^+) = \frac{\E[(D^-)^2D^+]-\E[(D^-)^2]}{\mu}$ 
\end{enumerate}
% \begin{remark}
% Conditions \ref{cond.beta} and \ref{cond.gamma} ensure that the Central Limit Theorem applies to the fluctuations of the first explored in-degrees around their mean. Condition \ref{cond.critical} ensures that the branching process corresponding to the depth-first exploration (i.e. the exploration of the out-components) is critical. Condition \ref{cond.rho} ensures that this branching process has Brownian scaling. Condition \ref{cond.tau} ensures that the covariance of the in- and out-degrees that are discovered first is finite. Condition \ref{cond.iota} ensures that the strongly connected components are $3$-regular. 
% \end{remark}
\subsection{Metric directed multigraphs and kernels}\label{subsec.mdmkernels}

\begin{figure}[htbp]
    \centering

    \includegraphics[width=\textwidth]{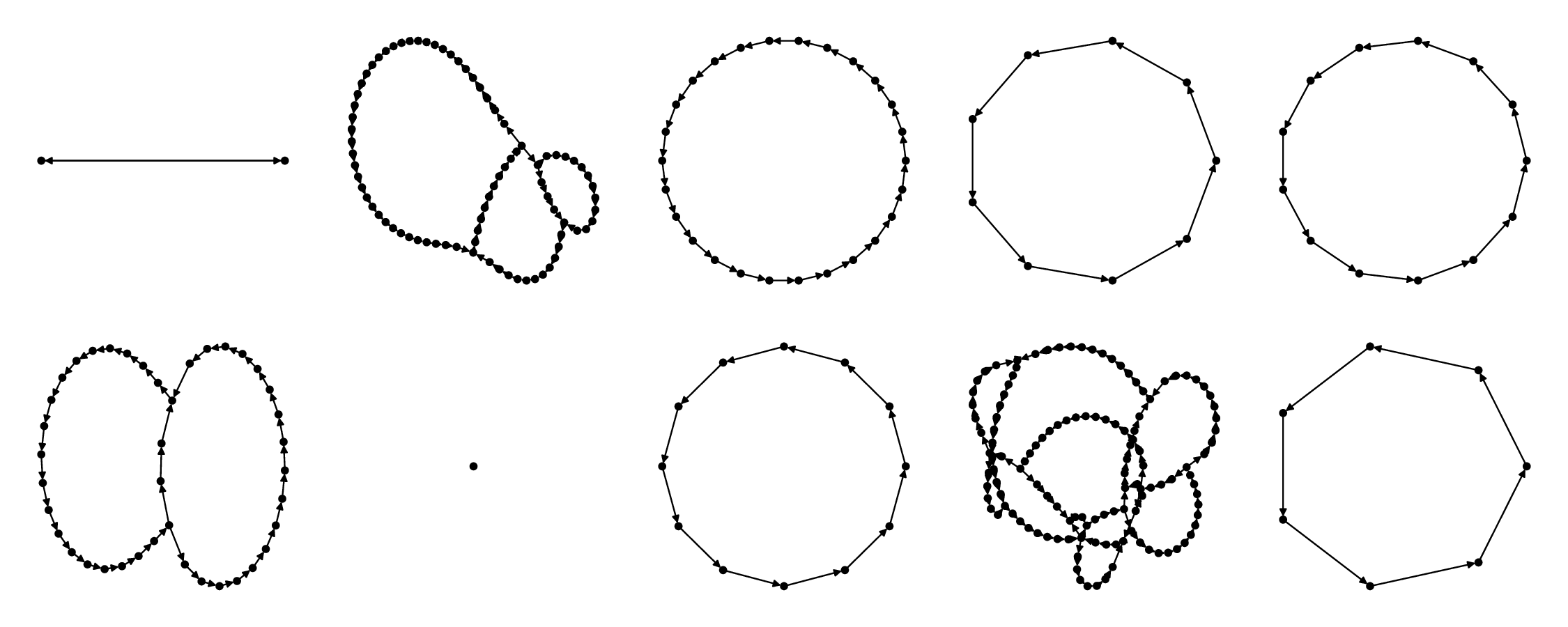}
    
    \caption{The largest SCC from samples of a directed configuration model with independent $\text{Poisson(1)}$ in- and out-degrees}
    \label{fig:largest-sccs}
\end{figure}

Figure \ref{fig:largest-sccs} shows the largest SCC from samples of a directed configuration model. As can be seen, while the lengths of paths in the SCC are long, the actual structure of the SCC is often quite simple. Previous work by \citet{goldschmidtScalingLimitCritical2019} shows that this is true for the \emph{directed Erdős--Rényi model} in the \emph{critical window}. The \emph{directed Erd\H{o}s--Rényi model} on $n$ vertices with parameter $p$, denoted by $\vec{G}(n,p)$, is a random digraph with vertex set $[n]$ in which each of the $n(n-1)$ possible directed edges is included with probability $p$ independently. The cases $p=(1+\lambda n^{-1/3})/n$ for $\lambda\in\R$ are referred to as \emph{the critical window}, and the case $p=1/n$ is called \emph{criticality}. In \cite{goldschmidtScalingLimitCritical2019}, it was shown that, for the directed Erdős--Rényi model in the critical window, while the lengths of paths in the SCCs scale like $n^{1/3}$, the combinatorial structure of the SCCs remains finite. The same turns out to be true in our setting.

This idea was formalised in \cite{goldschmidtScalingLimitCritical2019}, and we will use the same formalism in this work. We will first introduce \emph{metric directed multigraphs} (MDMs). These are simply weighted directed multigraphs, but in our context it is more appropriate to think of the weights as lengths, which motivates the change in naming. Formally, a \emph{directed multigraph} is a tuple $(V, E, r)$ where
\begin{enumerate}
    \item $V$ is a set of \emph{vertices},
    \item $E$ is a set of \emph{edges}, and
    \item $r: E \to V \times V$ is a function mapping each edge to its \emph{head} and \emph{tail}; associated with $r$ are two functions $r_1: E \to V$ and $r_2: E \to V$ such that
\begin{equation*}
    r(e) = (r_1(e), r_2(e))
\end{equation*}
for all $e \in E$. $r_1(e)$ is the tail of the edge $e$ and $r_2(e)$ is the head of the edge $e$.
\end{enumerate}
 Then a \emph{metric directed multigraph (MDM)} is a tuple $M = (V, E, r, l)$ where $(V, E, r)$ is a directed multigraph and $l:E \to [0, \infty)$. Let $\zeroloop$ denote the MDM consisting of a single vertex with a self-loop of length 0.

An \emph{isomorphism} between two MDMs $M = (V, E, r, l)$ and $M' = (V', E', r', l')$ is a pair of functions $(i_V, i_E)$ where $i_V: V \to V'$ and $i_E: E \to E'$ are bijections satisfying the relation
\begin{equation*}
    r'(i_E(e)) = (i_V(r_1(e)), i_V(r_2(e)))
\end{equation*}
for all $e \in E$. We say two MDMs are \emph{isomorphic} if there exists an isomorphism between them. In other words, isomorphic MDMs have the same graph structures for their underlying directed multigraphs up to a relabelling of the edges and vertices. Write $\iso(M, M')$ for the set of all isomorphisms between $M$ and $M'$.

\begin{figure}[htbp]
    \centering
    \begin{subfigure}[htbp]{0.45\textwidth}
        \centering
        \includegraphics[width=0.95\textwidth]{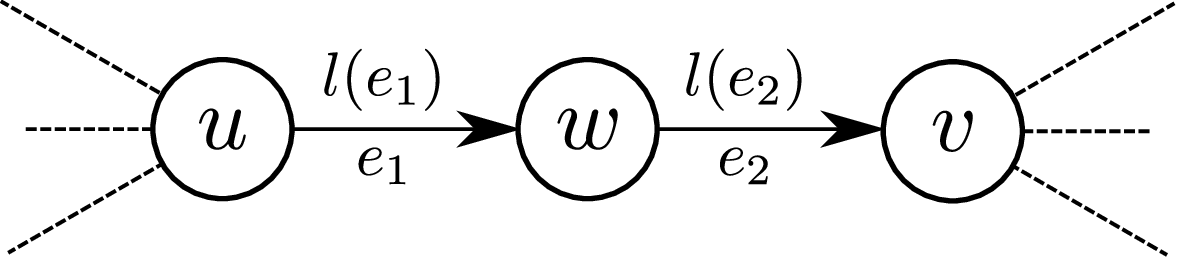}
        \caption{The graph before smoothing $w$}
    \end{subfigure}
    \hfill
    \begin{subfigure}[htbp]{0.45\textwidth}
        \centering
        \includegraphics[width=0.95\textwidth]{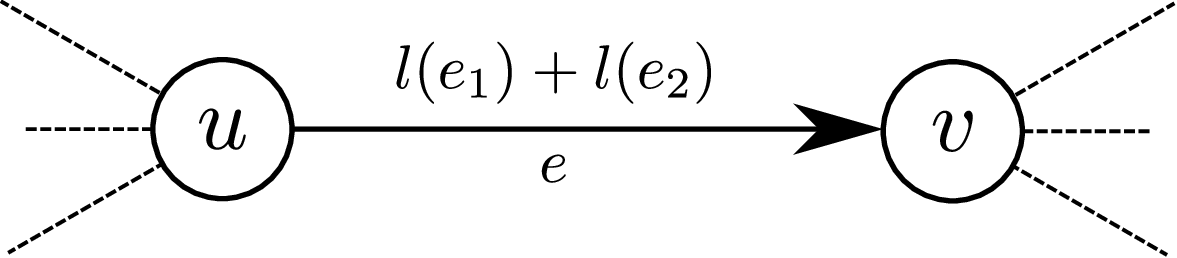}
        \caption{The graph after smoothing $w$}
    \end{subfigure}
    \caption{Smoothing a vertex $w$}
    \label{fig:smoothing}
\end{figure}

We now define a distance $\distmdm$ between two MDMs $M$ and $M'$.  Any isomorphism between $M$ and $M'$ gives a correspondence between the edges of $M$ and the edges of $M'$. We can then take an $\ell_{\infty}$ distance between the lengths of the edges and finally take the isomorphism which minimizes this distance. If $M$ and $M'$ are not isomorphic, we set the distance to be infinite. Formally,
\begin{equation*}
    \distmdm(M, M') = \begin{cases}
        \inf_{(i_V, i_E) \in \iso(M, M')} \sup_{e \in E} \abs{l(e) - l'(i_E(e))} & \text{if $M$ and $M'$ are isomorphic,} \\
        \infty & \text{otherwise.}
    \end{cases}
\end{equation*}
Consider an MDM $M$ and a vertex $w \in M$ with in-degree 1 and out-degree 1 which is not a self-loop. Let $u$ and $v$ be the unique in-neighbour and out-neighbour of $w$ respectively. The MDM obtained by \emph{smoothing} $w$ is obtained by deleting the edges $e_1$ and $e_2$ such that $r(e_1) = (u, w)$ and $r(e_2) = (w, v)$, then adding an edge $e$ such that $r(e) = (u, v)$ and assigning it length $l(e) = l(e_1) + l(e_2)$. This is illustrated in Figure \ref{fig:smoothing}. 

\begin{figure}[htbp]
    \centering
    \begin{subfigure}[htbp]{0.45\textwidth}
        \centering
        \includegraphics[width=0.90\textwidth]{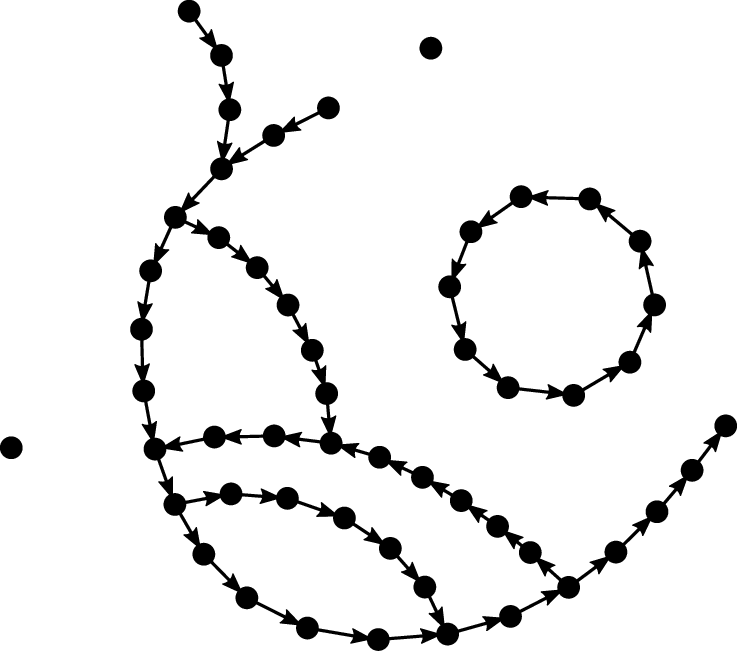}
        \caption{$\vec{G}$}
    \end{subfigure}
    \hfill
    \begin{subfigure}[htbp]{0.45\textwidth}
        \centering
        \includegraphics[width=0.90\textwidth]{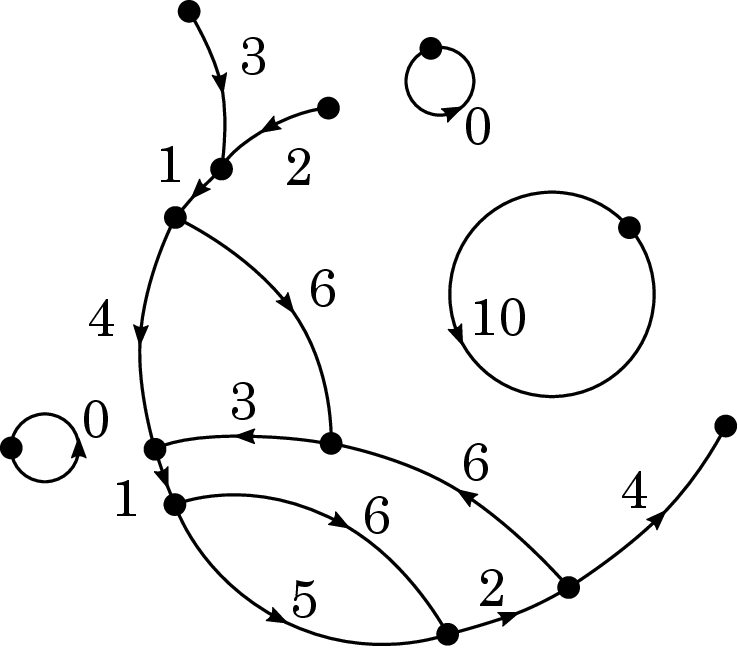}
        \caption{Kernel of $\vec{G}$}
    \end{subfigure}
    \caption{An example of a digraph $\vec{G}$ and its kernel. The numbers indicate the edge lengths.}
    \label{fig:kernel}
\end{figure}

Then the kernel of a digraph $\vec{G}$ is obtained by doing the following:
\begin{enumerate}
    \item Assign length $1$ to each edge.
    \item Iteratively smooth vertices with in-degree 1 and out-degree 1 that are not self-loops until there are none remaining.
    \item Replace all singletons by $\zeroloop$.
\end{enumerate}
An example is shown in Figure \ref{fig:kernel}.

\subsection{Our results}

For $M$ an MDM and $c\in (0,\infty)$, let $cM$ be equal to $M$ with all lengths multiplied by $c$. Let $C_i(n)$ for $i\geq 1$ be the kernels of the SCCs of $\vec{G}_n(\nu)$, listed in decreasing order of number of edges, breaking ties arbitrarily. Complete the list with an infinite repeat of $\zeroloop$. Then, our main theorem is as follows.
\begin{theorem}\label{thm.main}
There exists a sequence $\cC=(\cC_i,i\in \N)$ of random strongly connected MDMs such that 
$$\left(n^{-1/3}C_i(n),i\in \N\right)\todist\left(\cC_i,i\in \N\right)$$
as $n\to \infty$, with respect to the product $d_{\vec{\cG}}$-topology. The law of $\cC=(\cC_i,i\in \N)$ depends only on the parameters $\mu$, $\sigma_+$, and $(\sigma_{-+}+\nu_-)/\mu$. Further, for each $i\geq 1$, $\cC_i$ is either $3$-regular or a loop.
\end{theorem}
We will describe the limit object and some of its further properties in Subsection \ref{subsec.limitobject}.

The law of the limit object places some particular cases of our model in the universality class of the directed Erd\H{o}s--Rényi model as studied by \citet{goldschmidtScalingLimitCritical2019}. This is the content of the following corollary. Note however that their result holds in a stronger topology: they use an $\ell_1$-like topology on the space of sequences of MDMs, whereas we show our result in the product topology. Due to this, it is important in their paper to consider singletons as loops of length zero. For any fixed $k$, the $k$th largest SCC will not be a singleton with high probability as $n \to \infty$. Therefore, no component of our limiting object will be a singleton. Thus they need to pad their SCCs by $\zeroloop$ and consider the kernel of singletons to be $\zeroloop$, to prevent the $\ell_1$-distance, as defined by $d_{\vec{G}}$, between $\left(n^{-1/3}C_i(n),i\in \N\right)$ and $\left(\cC_i,i\in \N\right)$ being infinite. We follow the same convention.

\begin{theorem}\label{cor.erdosrenyi}
Consider $\vec{G}_n(\nu)$, with $\nu$ such that $$\mu=\sigma_+=\sigma_{-+}+\nu_-=1.$$ 
Let $(C^\nu_i(n), i\geq 1)$ be the kernels of the SCCs of $\vec{G}_n(\nu)$. Furthermore, let $(C^{ER}_i(n), i\geq 1)$ be the kernels of the SCCs of $\vec{G}(n,1/n)$. Then, $\left(n^{-1/3}C^\nu_i(n),i\in \N\right)$ and 
$\left(n^{-1/3}C^{ER}_i(n),i\in \N\right)$ have the same limit in distribution in the product-$d_{\vec{\cG}}$-topology as $n\to \infty$. 
\end{theorem}
Note that the condition in Corollary \ref{cor.erdosrenyi} is satisfied by $\nu(k^-,k^+)=\nu_1(k^-)\nu_2(k^+)$, with $\nu_1$ and $\nu_2$ the law of a $\operatorname{Poisson}(1)$ random variable.

Moreover, Theorem \ref{thm.main} has the following trivial corollaries, which were previously unknown. 
\begin{corollary}\label{cor.componentsizes}
Let $E^i_n$ and $V^i_n$ be the number of edges and vertices in $C_i(n)$ respectively, both appended with infinite repeats of $0$. Then there exists a random sequence $(E_i,i\in \N)\in \R_+^\infty$, such that
$$\left(n^{-1/3}E^n_i,n^{-1/3}V^n_i, i\in \N\right)\todist\left(E_i,E_i,i\in \N\right)$$
as $n\to \infty$ in the product topology on $(\R^2)^\infty$. 
\end{corollary}
In particular, note that, in the above corollary, the number of vertices and number of edges have the exact same scaling limit.
\begin{corollary}\label{cor.diameter}
For $v,w\in \vec{G}_n(\nu)$ such that $v\to w$, let $d(v,w)$ denote the length of the shortest directed path from $v$ to $w$, and let $$\operatorname{Diam}\left(\vec{G}_n(\nu)\right)=\max_{v,w\in V}\{d(v,w):v\to w\}$$ be the \emph{diameter} of $\vec{G}_n(\nu)$. Then, for any $\epsilon>0$, there is a $\delta>0$ such that $$\P\left( n^{-1/3}\operatorname{Diam}\left(\vec{G}_n(\nu)\right)>\delta\right)>1-\epsilon$$ for all $n$ large enough. Equivalently, $\operatorname{Diam}\left( \vec{G}_n(\nu) \right) = \Omega_p(n^{1/3})$.
\end{corollary}

\subsection{Previous work}\label{sec.previouswork}
The configuration model was introduced by \citet{bollobasProbabilisticProofAsymptotic1980} to sample a uniformly random undirected graph with a given degree sequence. (For a discussion of the configuration model and proofs of standard results, we refer the reader to \cite[Chapter 7]{hofstadRandomGraphsComplex2017}.)

Most results on the configuration model are proved for models with a deterministic degree sequence. The phase transition for the undirected setting was shown in \cite{molloyCriticalPointRandom1995, molloySizeGiantComponent1998, jansonNewApproachGiant2009}. The law of component sizes at criticality and in the critical window were obtained by \citet{riordanPhaseTransitionConfiguration2012} under the assumption that the degrees are bounded. Dhara, van der Hofstad, van Leeuwaarden and Sen showed convergence of the size and surplus edges in the critical window with a finite third moment \cite{dharaCriticalWindowConfiguration2017} and in the heavy-tailed regime \cite{dharaHeavytailedConfigurationModels2020}.  Bhamidi, Dhara, van der Hofstad and Sen obtained metric space convergence in the critical window in \cite{bhamidiUniversalityCriticalHeavytailed2020}, a result that the authors later improved to a stronger topology in \cite{bhamidiGlobalLowerMassbound2020}.

Configuration models with a random degree sequence are considered in \cite{josephComponentSizesCritical2014}, \cite{conchon--kerjanStableGraphMetric2021}, and \cite{Donderwinkel2021heightprocess}. \citet{josephComponentSizesCritical2014} showed convergence of the component sizes and surpluses of the large components under rescaling at criticality, both for degree distributions with finite third moments and for the heavy-tailed regime. \citet{conchon--kerjanStableGraphMetric2021} show Gromov-Hausdorff-Prokhorov convergence of the rescaled components ordered by decreasing size at criticality in these two regimes. The results in \cite{conchon--kerjanStableGraphMetric2021} in the heavy-tailed regime are extended to the critical window by the first author in \cite{Donderwinkel2021heightprocess}. Our techniques are closely related to the techniques introduced in \cite{conchon--kerjanStableGraphMetric2021}. 

Some results have been obtained for other directed graph models. \citet{caoConnectivityGeneralClass2019} consider a class of inhomogeneous directed random graphs. Their results include a phase transition for the existence of a giant SCC. This is a generalisation of work by Bloznelis, Götze and Jaworski in \cite{bloznelisBirthStronglyConnected2012}, in which a smaller class of inhomogeneous directed graphs is considered. Samorodnitsky, Resnick, Towsley, Davis, Willis and Wan \cite{samorodnitskyNonstandardRegularVariation2016} studied the tails of the degree distribution in the directed preferential attachment model. \citet{goldschmidtScalingLimitCritical2019} studied the directed Erd\H{o}s-R\'enyi model, and were the first to obtain metric space convergence of the SCCs of a directed graph. Our methods build on their techniques.

The directed configuration model was first considered by \citet{cooperSizeLargestStrongly2004}. They consider a deterministic degree sequence under a number of conditions. As discussed previously in \cref{subsec:model-description}, a phase transition for the SCCs occurs when a parameter $d=1$. They show that for $d<1$, with high probability, all SCCs contain $O(\Delta\log(n))$ vertices, for $\Delta$ the maximal degree. On the other hand, for $d>1$, there is a unique SCC that contains a positive proportion of the vertices and edges. Their conditions are restrictive, and include finite second moments for both the in- and out-degree of a uniformly chosen vertex, and a bound of size $n^{1/12}/\log(n)$ on the largest degree. Their proofs are based on an algorithm to explore the directed graph. The condition on the largest degree was later relaxed to $O(n^{1/4})$ by \citet{grafStronglyConnectedComponents2016}. These results are in contrast with the critical case, with Corollary \ref{cor.componentsizes}, which says that in our set-up the number of vertices and edges in the largest strongly connected components are $\Theta(n^{1/3})$ in probability.

Recently, Cai and Perarnau have obtained a number of results on the directed configuration model with deterministic degrees. In \cite{caiDiameterDirectedConfiguration2020}, they show, under first and second moment conditions of the degree of a uniformly picked vertex, for $d\neq 1$ (i.e. not at criticality), that the diameter of the model on $n$ vertices, rescaled by $\log(n)$ converges to a constant that they identify. This is in contrast with Corollary \ref{cor.diameter}, which says that in our set-up the diameter is $\Omega(n^{1/3})$ in probability at criticality. Then, in \cite{caiGiantComponentDirected2021}, they show a law of large numbers for the number of vertices and edges in the largest SCC, under slightly stronger moment conditions, and again away from the critical point. In \cite{caiMinimumStationaryValues2021}, they study the behaviour of a random walk on a directed configuration model.
 
A necessary and sufficient condition for the existence of a giant weakly connected component for the directed configuration model with a deterministic degree sequence is discussed in the physics literature by \citet{kryvenEmergenceGiantWeak2016}. He also studies the distribution of the in- and out-components in \cite{kryvenFiniteConnectedComponents2017}.

The directed configuration model with random in- and out-degrees is also considered by \citet{chenDirectedRandomGraphs2013} although, importantly, they do not allow for the in- and out-degree of a vertex to be dependent. The authors consider a model in which the in- and out-degrees are two independent sequences of i.i.d.\ random variables drawn from different probability distributions. They propose an algorithm to sample degree sequences that correspond to a simple graph and show the limiting distribution of the degrees generated by this algorithm.

\subsection{Proof outline}\label{sec:proofoutline}

\def \exploredvertices {\mathcal V}
\def \explorededges {\mathcal E}
\def \forest {F}
\def \edgestack {\mathcal Q}

The techniques we will use to investigate the graph model are a combination of the techniques introduced by Conchon-Kerjan and Goldschmidt in \cite{conchon--kerjanStableGraphMetric2021} and the strategy of Goldschmidt and Stephenson in \cite{goldschmidtScalingLimitCritical2019}, adapted to our set-up. The former work finds the scaling limit of an undirected uniform graph with i.i.d.\ degrees at criticality, and the latter finds the scaling limit of the SCCs of a directed Erd\H{o}s-Rényi graph at criticality. 

Our techniques use height processes and \L ukasiewicz paths, which are standard objects used to encode trees and forests (see for instance \cite[Chapter 0]{AST_2002__281__R1_0}). We will introduce these here. Let $T=(V,E,\rho)$ be an ordered rooted finite tree with vertex set $V$, edge set $E$ and root vertex $\rho$; say $|V|=n$. Let $v_0,\dots,v_{n-1}$ denote the vertices of the tree visited in depth-first order, so that $v_0=\rho$. We can view $T$ as a metric space by regarding all edges as line segments of length $1$ that are connected via the vertices. The distance $d_T$ between points $a_1$ and $a_2$ on line segments $l_1$ and $l_2$ respectively is then defined as the length of the unique non-self-intersecting path between $a_1$ and $a_2$ that traverses the line segments of the tree. Denote $(T,d_T)$ by $\mathrm{T}$.\\
We will define the height process and \L ukasiewicz path of $\mathrm{T}$. Both of these functions uniquely characterize $\mathrm{T}$. The height process of $\mathrm{T}$, referred to as $h$, is defined as $$h(i)=d_T(v_i,v_0),$$ i.e.  for all $i$, $h(i)$ equals the distance from $v_i$ to the root.
Moreover, for all $i=1,\dots,n$, let $y_i$ be the number of children of $v_{i-1}$, and set $y_0=1$. Then, the \L ukasiewicz path of $\mathrm{T}$ is defined by $$s(i)=\sum\limits_{j\leq i} (y_j-1)$$ for $i=0,\dots,n$. Then, $s(i)$ is the total number of younger siblings of $v_i$ and its ancestors.
For a sequence of ordered rooted finite trees, we define its height process by concatenating the height processes of the trees in the sequence. The \L ukasiewicz path is defined similarly. 

\begin{figure}
    \begin{subfigure}[htbp]{\textwidth}
        \centering
        \includegraphics[scale=1]{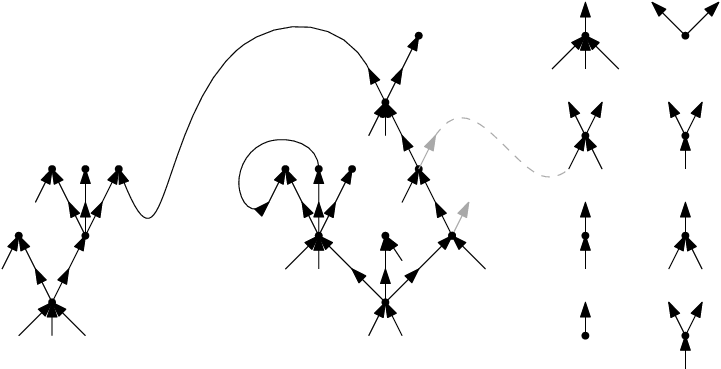}
        \caption{The gray arrows represent unpaired out-half-edges of vertices that have been discovered. One by one, in depth first order, these are paired to a uniform unpaired in-half-edge.}
        \label{fig.configuration model}
    \end{subfigure} \\

    \vspace{2em}

    \begin{subfigure}[htbp]{\textwidth}
        \centering
        \includegraphics[scale=1]{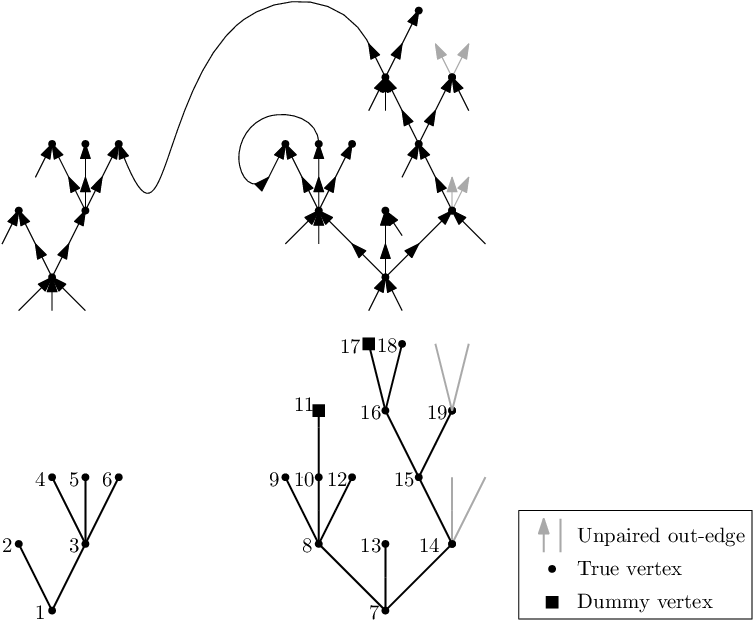}
        \caption{The out-forest is defined based on the exploration of the digraph. For each surplus edge, we add a dummy leaf. The labels of the vertices correspond to the time step in the exploration at which the vertex is added. The gray edges lead to vertices of which we do not know whether it is a dummy vertex, and if not, what its degree is. }
        \label{fig.configuration modeloutforest}
    \end{subfigure}

    \caption{Partial constructions of the configuration model and out-forest}
\end{figure}

We will study the law of the SCCs of a uniform directed graph with degree sequence $(\mathbf{D}_1,\dots,\mathbf{D}_n)$, conditional on $\sum_{i=1}^n D^-_i=\sum_{i=1}^n D^+_i$ by exploring the configuration model in a depth-first manner. This sampling naturally gives rise to a directed subforest of the resulting multigraph, which we call the \emph{out-forest}. The sampling procedure is  described in \cref{alg:edfs}, and is also illustrated in Figure \ref{fig.configuration model}. The definition of the out-forest is illustrated in Figure \ref{fig.configuration modeloutforest}.  

% This procedure differs from the exploration used by \citet{goldschmidtScalingLimitCritical2019}.  The eDFS takes as input a directed multigraph and uses a stack (ordered list) of edges. When the stack is empty, we are at the start of a new out-component and pick a new vertex $w$ with probability proportional to its in-degree. Otherwise we remove the last edge $(v, w)$ from the stack. In both cases, if $w$ is not yet in the list of discovered vertices, we add the out-edges from this vertex to the \emph{end} of the stack of edges (this choice is what makes the exploration depth-first) and add $w$ to the list of discovered vertices $\exploredvertices$.  The order in which vertices are added to $\exploredvertices$ is referred to as their \emph{order of discovery}. Note that we will not discover any vertex with in-degree 0. From the perspective of finding the SCCs, this is not a problem since such vertices will form singleton SCC.

The sampling procedure uses a queue of unpaired out-edges (represented by the label of their corresponding vertex). When the queue is empty, we are at the start of a new out-component and pick a new vertex $w$ with probability proportional to its in-degree if there are vertices with positive in-degree remaining. Else, we pick a new vertex uniformly at random. If the queue is not empty, we pair the first out-edge in the queue to a uniform unpaired in-edge and call the corresponding vertex $w$. In both cases, if $w$ is not yet in the list of discovered vertices, we add the out-edges from this vertex to the \emph{front} of the queue of edges (this choice is what makes the exploration depth-first) and add $w$ to the list of discovered vertices. The order in which vertices are added to the list of discovered vertices is referred to as their \emph{order of discovery}. 

This procedure will discover vertices with in-degree 0 last. This is fine since such vertices form singleton SCCs, so we have discovered the non-trivial SCCs first.

\begin{algorithm}[htbp]
    \SetAlgoLined
    \KwData{A set of vertices $V=\{v_1,\dots, v_n\}$ with degree pairs $(d^-(v_1),d^+(v_1)),\dots, (d^-(v_n),d^+(v_n))$ satisfying $\sum d^-(v_i)=\sum d^+(v_i)$}
    $\exploredvertices \leftarrow$ an empty ordered list of vertices \tcp*[f]{the list of discovered vertices}\;
    $\edgestack \leftarrow$ an empty ordered list of vertices \tcp*[f]{the queue}\;
    $(d^-_{\mathrm{unpaired}}(v_1),\dots,d^-_{\mathrm{unpaired}}(v_n))\leftarrow (d^-(v_1),\dots,d^-(v_n))$  \tcp*[f]{the number of unpaired in-edges per vertex}\;
    $k \leftarrow 0$ \tcp*[f]{the index of the current step} \;
    $\hat{s}^- \leftarrow 0$ \tcp*[f]{the number of unpaired in-edges of discovered vertices}\;
    $\hat{s}^+ \leftarrow 1$ \tcp*[f]{the queue size minus the number of explored out-components} \;
    $\forest \leftarrow$ a directed forest with vertices $V$ and no edges
    \tcp*[f]{current out-forest}\;
     $M \leftarrow$ a directed multigraph with vertices $V$ and no edges
    \tcp*[f]{current di-multigraph}\;
    \While(){there exist undiscovered vertices OR $\edgestack$ is non-empty}{
        \eIf(\tcp*[f]{we start a new out-component}){$\edgestack$ is empty}{
            \eIf(){there exist undiscovered vertices with positive in-degree}{
            $w\leftarrow$ a random vertex not in $\exploredvertices$ chosen with prob.\ proportional to $d^-(w)$ \;}
            {$w\leftarrow$ a uniformly random vertex not in $\exploredvertices$}
            $\hat{s}^+\leftarrow \hat{s}^+-1$\tcp*[f]{we have explored a component}\;}{
            $v \leftarrow$ first entry in $\edgestack$ \tcp*[f]{we will pair an unpaired out-edge of $v$}\;
            remove first entry from $\edgestack$ \;
            $\hat{s}^+\leftarrow \hat{s}^+-1$\tcp*[f]{the queue size decreases by $1$}\;
            $w \leftarrow$ a random vertex chosen with prob.\ proportional to $d^-_{\mathrm{unpaired}}(w)$ \;
            add $(v,w)$ to $M$\tcp*[f]{we pair the out-edge of $v$ with a uniform unpaired in-edge}\;
            $d^-_{\mathrm{unpaired}}(w)\leftarrow d^-_{\mathrm{unpaired}}(w)-1$\;
            $\hat{s}^-\leftarrow \hat{s}^--1$\tcp*[f]{we have paired an in-edge}\;
            \eIf(\tcp*[f]{we sampled a surplus edge}){$w\in \exploredvertices$}{
            add a dummy leaf to $F$ and an edge from $v$ to the leaf\;}{
            add $(v,w)$ to $F$\;
            }
        }
        \If(){$w\not\in \exploredvertices$}{
        append $w$ to the end of $\exploredvertices$ \tcp*[f]{vertex $w$ is now discovered}\;
        append $d^+(w)$ repeats of $w$ to the start of $\edgestack$\;
        $\hat{s}^+\leftarrow \hat{s}^++d^+(w)$\tcp*[f]{the queue size has increased}\;
        $\hat{s}^-\leftarrow \hat{s}^-+d^-(w)$\tcp*[f]{the number of unpaired in-edges of discovered vertices has increased}\;
        }
        $k\leftarrow k+1$\;
        $\hat{s}^+_k\leftarrow \hat{s}^+$\;
        $\hat{s}^-_k\leftarrow \hat{s}^-$ \;
        }

    \caption{The edge depth-first configuration model \label{alg:edfs}}
\end{algorithm}

At each step we also track two natural numbers $\hat{s}^-(k)$ and $\hat{s}^+(k)$. The first one, $\hat{s}^-(k)$ keeps track of the number of unpaired in-edges of discovered vertices at time $k$. The second one, $\hat{s}^+(k)$ is akin to a \L{}ukasiewiscz path. At any given step it is equal to the size of the queue after subtracting the number of fully explored out-components.

We also construct a directed forest for which $\hat{s}^+(k)$ will be the true \L{}ukasiewicz path. At each step of the process we will examine a vertex $w$. If $w$ has not been discovered yet then either we are at the start of a new out-component, in which case we make $w$ the root of the next out-component, or we added an edge $(v, w)$ to the multigraph with $v$ already discovered, in which case we add the edge $(v, w)$ to the out-forest as well. If $w$ has already been explored we cannot add $(v, w)$ to the out-forest without creating cycles or connecting two different components. We instead add a \emph{dummy leaf} to the out-forest and an edge from $v$ to the dummy leaf.  We call any vertex that is not a dummy leaf a  \emph{true vertex}. This is illustrated in Figure \ref{fig.configuration modeloutforest}.

Consider an edge $(v,w)$ in the directed multigraph. If $(v, w)$ is not in the out-forest we refer to the edge as \emph{surplus}. Such an edge will instead correspond to an edge $(v, d)$ in the out-forest where $d$ is a dummy leaf. 

An important motivation for studying the out-forest is the fact that the vertex set of any SCC is contained in one of the components of the out-forest. This is a straightforward property which we will prove below as part of \cref{lem:whatispartofscc}. Moreover, we defined the out-forest in such a way that every time step in the exploration corresponds to one vertex in the out-forest.

Our technique relies on dismissing surplus edges that cannot be part of a strongly connected component (for example, surplus edges between two different out-components cannot form a directed cycle and are never part of a strongly connected component). We define a necessary condition for a surplus edge to be part of an SCC (see Definition \ref{def.candidate} and \cref{prop:edgesinSCCs}), and we call dummy leaves that correspond to surplus edges with this property \emph{candidates}. Then, we define a procedure to sample only the out-forest and the edges corresponding to candidates, which allows us to find the SCCs.

A key fact is that the order in which the true vertices are discovered does not depend on the positions of the dummy leaves. Similarly, the positions of the dummy leaves do not depend on the position of the heads of the surplus edges. Finally, whether a dummy leaf is a candidate does not depend on the position of the heads of the surplus edges. This allows us to define the following step-by-step sampling procedure. 
\begin{enumerate}
    \item We sample the order of discovery of the true vertices.
    \item We sample at which time steps we add a dummy leaf instead of a true vertex. 
    \item For each dummy leaf we sample whether it is a candidate.
    \item For each candidate we sample the position of the head of the corresponding surplus edge.
\end{enumerate}
For an exact description of the sampling procedure, see Subsection \ref{subsec.discrete}. The analogous sampling procedure for the limit object is described in Subsection \ref{subsubsec.samplecontinuousobject}. Then, our approach to show convergence is as follows.

\begin{enumerate}
    \item We find the limit under rescaling of the \L ukasiewicz path and height process of the out-forest up to time $m_n=\Theta(n^{2/3})$ conditional on the event $\left\{\sum_{i=1}^n D^-_i=\sum_{i=1}^n D^+_i\right\}$. This is the content of  \cref{prop:convoutforest}. Note that we condition on an asymptotically singular event, which causes significant difficulties.   Our method relies on a measure change between the sequence of degrees in order of discovery under this conditioning and a sequence of i.i.d.\ random variables in $\N\times\N$. In Section \ref{sec:measure-change} , we show the convergence of the measure change under rescaling.
    \item We establish that the positions of the tails of the surplus edges corresponding to the candidates converge. This is the content of Proposition \ref{prop.convergenceancestraledges}, Lemma \ref{prop.extractexcursions}, and Proposition \ref{prop.convergencestartingpointscandidates}.
    \item We show that the positions of the heads of the surplus edges corresponding to the candidates converge, which is the content of Proposition \ref{prop.convergenceheadscandidates}. 
    \item We identify the tails and heads of the surplus edges corresponding to the candidates, and recover the SCCs from the resulting digraph via a cutting procedure. We use a result from \cite{goldschmidtScalingLimitCritical2019} to show that the cutting procedure converges. This summarised in Corollary \ref{prop:sccordereduptotimeT}.
    \item We show that conditioning on the resulting multigraph being simple does not affect the sampling procedure on the time scale $O(n^{2/3})$. This is the content of Proposition \ref{prop.anomalousedges}. 
    \item We prove that for any $\delta>0$, with high probability, all SCCs with more than $\delta n^{1/3}$ edges are contained in the exploration up to time $O(n^{2/3})$. Therefore, we can choose $m_n$ such that, with high probability, we do not miss any large SCCs by not considering the exploration beyond time $m_n$, which finishes the proof of the convergence in the product topology. This is the content of Lemma \ref{lemma.largesccfoundfirst}. 
\end{enumerate}

\section{Sampling the MDM in the discrete and the continuum}
If we forget about the directions of the edges in $\vec{G}_n(\nu)$, the resulting undirected graph is supercritical, and, with high probability, the graph contains a unique giant component with surplus going to infinity as $n\to \infty$ (see e.g. \cite{molloyCriticalPointRandom1995,molloySizeGiantComponent1998,jansonNewApproachGiant2009} for a discussion of the phase transition in the undirected configuration model). This suggests that if we do not dismiss a large amount of edges, we will not be able to study the digraph in enough detail to find a metric space scaling limit of the SCCs. Therefore, we will not try to sample the entire digraph, but focus on the information that we need to find the SCCs. We start by studying the discrete digraph model, with the goal of identifying which edges can be part of an SCC, and how to sample them. In Subsection \ref{subsubsec.defcandidates}, we establish necessary conditions for an edge to be part of an SCC. These conditions imply that we only need to study the out-forest, and the surplus edges corresponding to a small subset of the dummy leaves. We call these dummy leaves \emph{candidates}. In Subsections \ref{subsubsec.samplingoutforest} and \ref{subsubsec.samplecandidates} we study the law of the out-forest and the surplus edges corresponding to the candidates respectively, and we define a procedure to sample them both. This yields a sequence of directed multigraphs with edge lengths in which the SCCs are embedded.  
In Subsection \ref{subsec.limitobject}, we define the continuous counterpart of the sampling procedure. The resulting object will be the limit under rescaling of the sequence of directed multigraphs with edge lengths in which the SCCs are embedded that was constructed in Subsections \ref{subsubsec.samplingoutforest} and \ref{subsubsec.samplecandidates}. 
\subsection{The discrete case}\label{subsec.discrete}
We will discuss the different type of edges that we can encounter in the exploration. Recall from Subsection \ref{sec:proofoutline} that by slight abuse of terminology, we call the dummy leaf that corresponds to a surplus edge its tail.

\subsubsection{Necessary conditions for an edge to be part of an SCC}\label{subsubsec.defcandidates}
\begin{figure}
\centering
\begin{subfigure}{0.8\textwidth}
 \centering
    \includegraphics[scale=1.2]{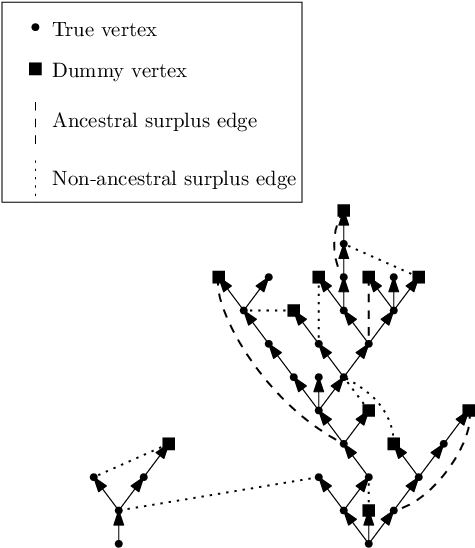}
    \caption{This figure illustrates an example of a depth-first exploration of two out-components with the different type of surplus edges highlighted. The ancestral surplus edges point from a vertex $v$ to one of its ancestors. They are always part of an SCC.}
    \label{subfigure.typesofsurplusedges} 
\end{subfigure}\\
\begin{subfigure}{0.8\textwidth}
  \centering
  \includegraphics[scale=1.2]{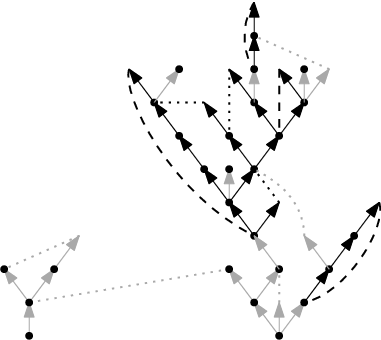}
  \caption{The edges that are part of an SCC are depicted in black. Two vertices are in the same SCC if and only if they are connected by black edges. }
    \label{subfigure.sccinexample}
\end{subfigure}
\caption{We illustrate the different types of surplus edges and how they affect the structure of the SCCs.}
\end{figure}

Amongst the surplus edges, \emph{ancestral surplus edges}, which are surplus edges that point from a vertex to one of its ancestors, play a special role. All other surplus edges are called \emph{non-ancestral}. This is illustrated in Figure \ref{subfigure.typesofsurplusedges}. In Figure \ref{subfigure.sccinexample} we show how surplus edges affect the structure of the SCCs. This is the content of the next lemma.

\begin{lemma}\label{lem:whatispartofscc}
The following facts hold for SCCs. 
\begin{enumerate}
\item \label{item.factsonsccs1}The vertices of an SCC are contained in precisely one of the components of the out-forest. 
\item \label{item.factsonsccs2} Ancestral surplus edges are always part of an SCC.
\item \label{item.factsonsccs4} A non-ancestral surplus edge is part of an SCC only if its head is an ancestor of the tail of a surplus edge that is part of an SCC.
\item \label{item.factsonsccs4andabit} An edge in the out-forest is part of an SCC only if its head is an ancestor of the tail of a surplus edge that is part of an SCC.
\item \label{item.factsonsccs5} For any non-trivial SCC, the first surplus edge of the SCC that is explored is an ancestral surplus edge, and a component of the out-forest contains an SCC if and only if it contains an ancestral surplus edge.
\end{enumerate}
\end{lemma}
\begin{proof}
We start with \ref{item.factsonsccs1}. Let $v$ and $w$ be two vertices in the same SCC. Without loss of generality, $v$ is explored first in depth-first order in the out-direction. Since $v$ and $w$ are part of the same SCC, we know that there is a path from $v$ to $w$ in the out-direction. This implies that $w$ will be part of the out-subtree consisting of the descendants of $v$. This implies that they are part of the same component of the out-forest.

To prove \ref{item.factsonsccs2}, suppose there is an ancestral surplus edge from $v$ to $w$. This implies that $w$ is an ancestor of $v$ in an out-component, which implies that there is a path from $w$ to $v$ as well. It follows that $w$ and $v$ are in the same SCC and that the ancestral surplus edge from $v$ to $w$ is in this SCC as well. 

To prove \ref{item.factsonsccs4} and \ref{item.factsonsccs4andabit}, suppose there is a non-ancestral surplus edge from $v$ to $w$ that is part of an SCC, or that $(v,w)$ is an edge in the out-forest that is part of an SCC. Then, there is some directed path $(x_0,\dots, x_m)$ with $x_0=w$ and $x_m=v$. Let $k$ be minimal such that $x_k$ is not a descendant of $w$ (such a $k$ exists, because by assumption, $v$ is not a descendant of $w$). Then, $(x_{k-1},x_k)$ is a surplus edge that is in the same SCC as $v$ and $w$, and $x_{k-1}$ is a descendant of $w$ by definition of $k$.

Finally, \ref{item.factsonsccs2} and \ref{item.factsonsccs4} imply \ref{item.factsonsccs5}. 
\end{proof}
 \cref{lem:whatispartofscc} motivates the following definition.
\begin{definition}\label{def.candidate}
A dummy vertex is a \emph{candidate} if one of the following statements holds for the surplus edge that it corresponds to. 
\begin{itemize}
    \item It is an ancestral surplus edge, or
    \item Its head is the ancestor of a candidate.
\end{itemize}
\end{definition}
The following proposition is at the core of our strategy to study the SCCs.
\begin{proposition}\label{prop:edgesinSCCs}
Any edge that is part of an SCC is either a surplus edge corresponding to a candidate, or is contained in the subforest of the out-forest that is spanned by the candidates and the roots of the out-components.
\end{proposition}
\begin{proof}
This follows from Definition \ref{def.candidate} and \cref{lem:whatispartofscc}.
\end{proof}

\cref{prop:edgesinSCCs} implies that to sample the SCCs, we do not need to sample the heads corresponding to all dummy leaves. Instead, for every dummy leaf, we only need to know whether it is a candidate, and if so, where its head is. 
\subsubsection{Sampling the out-forest}\label{subsubsec.samplingoutforest}
This subsection discusses how to obtain the out-forest conditional on the order in which the vertices are discovered. We will study the law of the degrees in order of discovery in \cref{sec:measure-change}. The out-forest is obtained in the following way. Let $(\mathbf{\hat{D}}_{n,1},\dots,\mathbf{\hat{D}}_{n,n})$ be the degree pairs in order of discovery (i.e. the order given by $\exploredvertices$ in Algorithm \ref{alg:edfs}). Up to time-step $k$, suppose we have discovered the first $m\leq k<n$ elements of  $(\mathbf{\hat{D}}_{n,1},\dots,\mathbf{\hat{D}}_{n,n})$. Then, at time $k+1$,
\begin{enumerate}
    \item If we have finished a component of the out-forest, let the next out-component have a root with out-degree $\hat{D}_{n,m+1}^+$. 
    \item Otherwise,
    \begin{enumerate}\item With probability proportional to the total in-degree of the undiscovered vertices, i.e. $\sum_{i={m+1}}^{n} \hat{D}_{n,i}^-$, let the next vertex in depth-first order be a true vertex with out-degree $\hat{D}_{n,m+1}^+$. 
    \item With probability proportional to the number of unpaired in-half-edges of the $m$ discovered vertices, let the next vertex in depth-first order be a dummy leaf, and reduce the total number of unpaired in-edges of the $m$ discovered vertices by $1$.
\end{enumerate}
\end{enumerate}
We make this rigorous in the following proposition.
\begin{proposition}\label{prop:sampleoutforest}
Suppose that the sequence of degrees in order of discovery $(\mathbf{\hat{D}}_{n,1},\dots,\mathbf{\hat{D}}_{n, n})$ is given. Suppose that after time-step $k$, there are still unpaired out-half-edges. Suppose that for $1\leq l\leq k$, that up to time $l$, $\hat{P}_n(l)$ surplus edges have been sampled. Then, $$\left(\hat{S}^+_n(l),1\leq l\leq k \right):=\left(\sum_{i=1}^{l-\hat{P}_n(l)}\hat{D}^+_{n,i}-l,1\leq l\leq k\right)$$ is the \L ukasiewicz path of the out-forest up to time $k$. Moreover, for $$\left(\hat{I}^+_n(l),1\leq l\leq k\right):=\left(\min\left\{\hat{S}^+_n(m):1\leq m \leq l\right\},1\leq l \leq k \right),$$
define 
$$\left(\hat{S}^-_n(l),1\leq l \leq k\right):=\left(\sum_{i=1}^{l-\hat{P}_n(l)}\hat{D}^-_{n,i}-l-\hat{I}^+_n(l)+1,1\leq l\leq k\right),$$
so that $\hat{S}^-_n(k)$ is equal to the number of unpaired in-half-edges of discovered vertices at time $k$. Then, the probability that we sample a surplus edge at the $(k+1)$th time-step is given by
$$\frac{\hat{S}^-_n(k+1)}{\sum_{i=1}^n D^-_i-k-\hat{I}^+_n(k)+1}\one_{\left\{\hat{I}^+_n(k)=\hat{I}^+_n(k-1)\right\}}.$$
We do not need to know the position of the heads of the surplus edges in order to sample the out-forest.
\end{proposition}
\begin{proof}
Note that if up to time $k$, $\hat{P}_n(k)$ surplus edges have been sampled, this implies that $k-\hat{P}_n(k)$ true vertices have been discovered. Thus, up to time $k$, the out-forest contains $\hat{P}_n(k)$ dummy leaves, and true vertices with degrees $(\hat{D}^+_{n,1},\dots,\hat{D}^+_{n,k-\hat{P}_n(k)})$, so by definition of the \L ukasiewicz path, its value is indeed equal to $\hat{S}^+_n(k)$ at time $k$. Moreover, up to time $k$, the total in-degree of the discovered true vertices is equal to $\sum_{i=1}^{k-\hat{P}_n(k)}\hat{D}^-_{n,i}$. At every time-step, we pair one in-half-edge of a discovered vertex, unless we start a new component. The value $-\hat{I}^+_n(k)$ corresponds to the number of out-components that are fully explored up to time $k$, so the total number of unpaired in-half-edges of discovered vertices at time $k$ is equal to $\hat{S}^-_n(k)$. By the same reasoning, the total number of unpaired in-half-edges is equal to $\sum_{i=1}^n D^-_i-k-\hat{I}^+_n(k)+1$. The probability of sampling a surplus edge at step $(k+1)$ follows. We note that this probability does not depend on the positions of the heads of the surplus edges, but only on their number, which implies that we can sample the out-forest without sampling the positions of the heads.
\end{proof}
\subsubsection{Sampling the candidates}\label{subsubsec.samplecandidates}
We will now study the law of the candidates and their heads conditional on the out-forest. We will first identify the candidates amongst the dummy leaves, and then we will sample the positions of their heads.

If the vertex discovered at time $k$ is a dummy leaf, the head of the corresponding surplus edge is a uniform pick from the $\hat{S}^-(k)$ unpaired in-half-edges of discovered vertices at time $k$. Therefore, the probability that a dummy leaf added at time $k$ corresponds to an ancestral surplus edge is given by the number of unpaired in-edges on its path to the root divided by $\hat{S}^-(k)$. This implies that to understand the law of the position of ancestral surplus edges, we need to understand where the unpaired in-edges are.

We will study this by modifying the edge lengths in the out-forest. We extend our definitions in \cref{sec:proofoutline} to trees with edge lengths as follows. Suppose $T=(V,E,\rho)$ is an ordered rooted finite tree, and suppose we have a function $\ell:E\to [0,\infty)$. Then, we can view $T$ as a metric space by regarding an $e$ as a line segment with length $\ell(e)$. The distance $d^\ell_T$ between points $a_1$ and $a_2$ on line segments $l_1$ and $l_2$ respectively is then defined as the length of the unique non-self-intersecting path between $a_1$ and $a_2$ that traverses the line segments of the tree, and we denote the resulting metric space $(T,d^\ell_T)$ by $\mathrm{T}^\ell$, and call it a \emph{ordered rooted finite tree with edge lengths}. This gives rise to an alternative height process, referred to as $h^\ell$, which is defined $$h^\ell(i)=d^\ell_T(v_i,v_0),$$ i.e.  for all $i$, $h^\ell(i)$ equals the distance from $v_i$ to the root in $\mathrm{T}^\ell$. We set the \L ukasiewicz path of $\mathrm{T}^\ell$ equal to the \L ukasiewicz path of $\mathrm{T}$.

We will now study the positions of the unpaired in-edges by modifying the edge lengths as follows: for a vertex $v$ with in-degree $m$, the edges connecting it to its children will all have length $m-1$ (unless $v$ is the root of an out-component, in which case the edges connecting to its children will be assigned length $m$). The height of vertex $w$ in this forest with modified edge lengths corresponds to the number of in-half-edges that can be used to form an ancestral surplus edge with tail $w$. We assign lengths to all edges in the out-forest and call the resulting forest with edge lengths \emph{the out-forest with edge lengths}. Denote the height process of the out-forest with edge lengths by $(\hat{H}_n^\ell(k),k\geq 1)$. 
Recall from Lemma \ref{lem:whatispartofscc} that the surplus edge corresponding to the first candidate in any component of the out-forest is ancestral. The following proposition illustrates the importance of $\hat{H}^\ell$ in finding the first ancestral surplus edges in the out-components.

\begin{proposition}\label{prop:probancestral}
Consider the exploration of the out-forest at time $k$. If no ancestral surplus edge has been sampled in the current component, then the probability that the $k$th vertex in depth-first order is a candidate is given by 
$$\frac{\hat{H}^\ell(k)}{\hat{S}_n^-(k)}\one_{\{\hat{P}_n(k)-\hat{P}_n(k-1)=1\}}.$$
This event is conditionally independent of the positions of the heads of the surplus edges that were found before time $k$, given that none of them were ancestral in the current component.
% If $k$ is the tail of an ancestral surplus edge, then the position of the end point $Y$ has the following law. Let $U$ be uniform on $[0,\hat{H}^\ell(k)]$. Then, let $Y_k$ be the height of youngest ancestor $l$ of $k$ such that $H^{\ell}(l)<U$. 
\end{proposition}
\begin{proof}
We claim that if no ancestral surplus edge has been sampled in the current component, none of the ancestors of $k$ are the head of a surplus edge. Indeed, for $x$ an ancestor of $k$, all vertices that are discovered since the discovery of $x$ up to time $k$ are descendants of $x$, because the out-forest is explored in a depth-first manner. Therefore, any surplus edge with head $x$ sampled up to time $k$ is ancestral. This implies that for $d^-$ the in-degree of $x$, the number of unpaired in-half-edges of $x$ at time $k$ is equal to $d^--1$ (unless $x$ is the root of the out-component, in which case it has $d^-$ unpaired in-half-edges).

Therefore, the number of unpaired in-half-edges corresponding to ancestors of $k$ is equal to $H^\ell(k)$. Moreover, note that, by definition of the dummy leaves, $k$ is the tail of a surplus edge if and only if $k$ is a dummy leaf, i.e. if and only if $\hat{P}_n(k)-\hat{P}_n(k-1)=1$. In that case, the probability that it connects to given unpaired in-half-edge of a discovered vertex is equal to $1/\hat{S}_n^-(k)$. The stated probability follows. The independence of the positions of the heads of earlier surplus edges is immediate.
\end{proof}

We now illustrate how to find the other candidates in a component of the out-forest. 

Let $T^n_{g_n}$ be a component of the out-forest with root $g_n+1$ and component size $\sigma_n$. Suppose the first ancestral surplus edge with vertices in $T^n_{g_n}$ corresponds to a dummy leaf $V^n_1\in [g_n+2,g_n+\sigma_n]$. Let $V^n_1<k\leq g_n+\sigma_n$, and suppose the candidates found up to time $k$ are given by $V^n_1,\dots,V^n_m$. Let $T^{n,\mathrm{mk}}_k$ be the subtree of $T^n_{g_n}$ spanned by $\{g_n+1,V^n_1,\dots,V^n_m,k\}$, and let $\ell(T^{n,\mathrm{mk}}_k)$ be its total length with edge lengths as encoded by $(\hat{H}^\ell(i),i\in [g_n+1,g_n+\sigma_n])$.
\begin{proposition}
\label{prop:samplecandidates}
 The probability that $k$ is a candidate is given by 
$$\frac{\ell\left(T^{n,\mathrm{mk}}_k\right)-m}{\hat{S}^-(k)}\one_{\{\hat{P}_n(k)-\hat{P}_n(k-1)=1\}}.$$
\end{proposition}
\begin{proof}
Note that if $k$ is a dummy leaf, it gets paired to a uniform pick from the $\hat{S}^-(k)$ as-yet unpaired in-half-edges of discovered vertices. By Definition \ref{def.candidate}, in that case, $k$ is a candidate if and only if the head of its corresponding surplus edge is in $T^{n,\mathrm{mk}}_k$. Observe that $\ell\left(T^{n,\mathrm{mk}}_k\right)$ is equal to the number of in-half-edges of $T_{k}$ that can be used to form surplus edges. By the definition of a candidate, exactly $m$ of those have been paired: one for each element in $\{V^n_1,\dots,V^n_m\}$. This implies that $\ell\left(T^{n,\mathrm{mk}}_k\right)-m$ of the $\hat{S}^-(k)$ options will cause $k$ to be a candidate.
\end{proof}

Note that the probability that a dummy leaf is a candidate only depends on the out-forest and the number of candidates that have been found in the component so far. The position of the heads of the surplus edges corresponding to candidates can be found as follows.

Let $T^n_{g_n}$ be a component of the out-forest with root $g_n+1$ and component size $\sigma_n$. Suppose its candidates are given by $\{V^n_1,\dots,V^n_{N_n}\}$. Then, for $1\leq i\leq {N_n}$, suppose the heads of the surplus edges corresponding to $V^n_1,\dots,V^n_{i-1}$ are given by $W_1^n,\dots,W^n_{i-1}$ respectively. 
\begin{proposition}
\label{prop.sampleheadcandidates}
The in-half-edge that $V^n_{i}$ gets paired to is a uniform pick from the $$\ell\left(T^{n,\mathrm{mk}}_{V^n_i}\right)-(i-1)$$ unpaired in-half-edges of $T^{n,\mathrm{mk}}_{V^n_i}$ that remain.
\end{proposition}
\begin{proof}
Given that $V^n_{i}$ is a candidate, its head will be in $T^{n,\mathrm{mk}}_{V^n_i}$. Then, the distribution follows.
\end{proof}

\begin{figure}
\centering
\begin{subfigure}{.7\textwidth}
 \centering
    \includegraphics[width=0.8\linewidth]{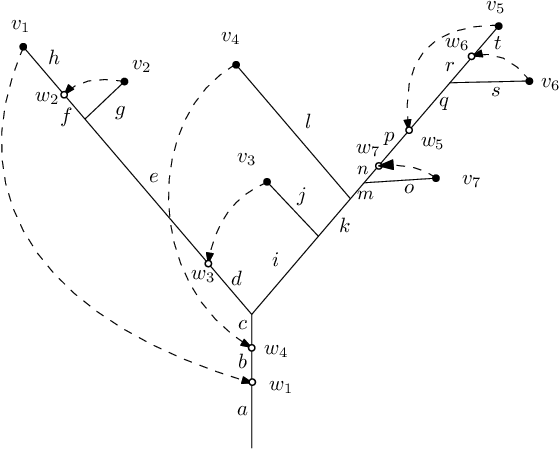}
    \caption{This is a subtree of an out-component spanned by the root of the out-component and the candidates $(v_1,\dots,v_7)$. Call the marked tree $T^{\mathrm{mk}}$. The heads of the surplus edges corresponding to candidates are denoted by $(w_1,\dots,w_7)$. }
\label{figure.extractSCCs1}
\end{subfigure}\\
\vspace{1.5em}
\begin{subfigure}{.8\textwidth}
  \centering
  \includegraphics[width=0.9\linewidth]{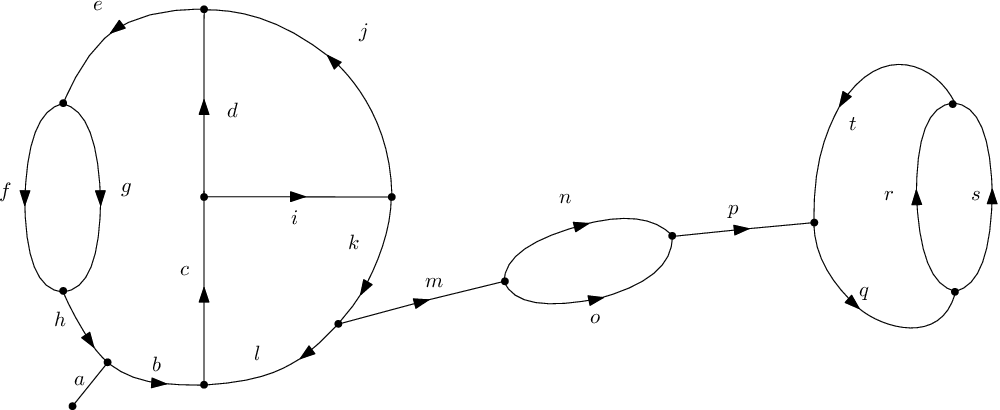}
  \caption{Identifying $v_i$ with $w_i$ for $i\in [7]$ gives $M$.}
  \label{figure.extractSCCs2}
\end{subfigure}\\
\vspace{1.5em}
\begin{subfigure}{.8\textwidth}
  \centering
  \includegraphics[width=0.9\linewidth]{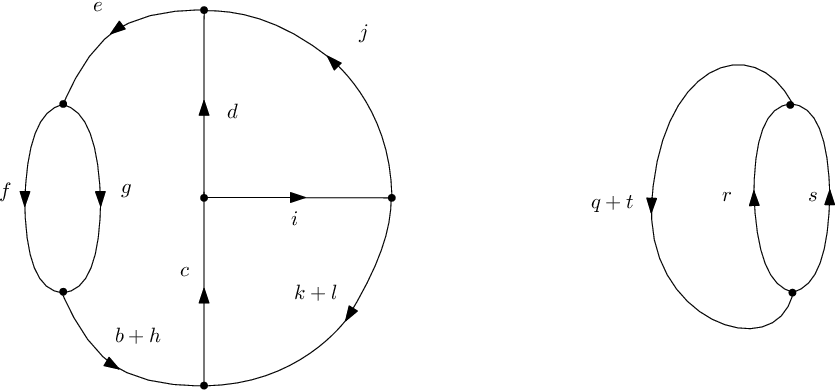}
  \caption{We find the SCCs that are contained in $M$.}
  \label{figure.extractSCCs3}
\end{subfigure}

\caption{We illustrate the procedure to find the SCCs in a component of the out-forest after finding the candidates. Taken from \cite{goldschmidtScalingLimitCritical2019} with permission of the authors.}
\end{figure}

Propositions \ref{prop:sampleoutforest}, \ref{prop:probancestral}, \ref{prop:samplecandidates}, and \ref{prop.sampleheadcandidates} justify the following sampling procedure.
\begin{enumerate}
    \item Sample the out-forest, and suppose it has $N$ vertices.
    \item Define a counting process $(A_n(k),k\leq N)$, with the probability of an increment at time $k$ given by $$\frac{\hat{H}_n^\ell(k)}{\hat{S}_n^-(k)}\one_{\{\hat{P}_n(k)-\hat{P}_n(k-1)=1\}}.$$
    \item For $i\geq 1$, let $X_i^n=\min\{k:A_n(k)=i\}$ be the time that the $i$th ancestral surplus edge is sampled. For $i\geq 1$, let $G_i^n$ be the left endpoint of the excursion of $\hat{S}^{+}_n$ above its running infimum that encodes the out-component that contains the $i$th ancestral surplus edge, and let $\Sigma_i^n$ be the length of this excursion, i.e. 
\begin{align*}G_i^n&=\min\left\{k\geq 1:\hat{S}^{+}_n(k)=\min\{\hat{S}^{+}_n(l):l\leq X_i^n\}\right\}\\
\Sigma_i^n&=\min\left\{k \geq 1: \min\left\{\hat{S}^{+}_n(l):l\leq G_i^n+k\right\} < \min\left\{\hat{S}^{+}_n(l):l\leq X_i^n\right\}\right\},
\end{align*}
so that for each $i\geq 1$, the excursion $\left(\hat{S}^+(k),k\in [G_i^n+1,G_i^n+\Sigma_i^n]\right)$ encodes the out-tree containing $X_i^n$. For each $(g_n,\sigma_n)\in \{(G_i^n,\Sigma_i^n)\}$, let $T^n_{g_n}$ be the tree in out-forest with root $g_n+1$, and do the following.
    \begin{enumerate}
    \item \label{item.procedure3} Set $V_1^n=\min\{m\geq 1:A_n(m)=A_n(g_n)+1\}$, and find the other candidates $\{V_2^n,\dots ,V_{N_n}^n\}$ according to the procedure described in the statement of \cref{prop:samplecandidates}.
    \item \label{item.procedure4} For the tails $V_1^n,\dots, V_{N_n}^n$, sample their corresponding heads $W_1^n,\dots ,W_{N_n}^n$ respectively according to the procedure described in the statement of \cref{prop:samplecandidates}.
    \item Let $T^{n,\mathrm{mk}}(g_n)$ be the subtree of $T^n_{g_n}$ spanned by $\{g_n+1,V_1^n,\dots ,V_{N_n}^n\}$. Then, quotient it by the equivalence relation $\sim$ which identifies $V_i^n$ and $W_i^n$ for each $1\leq i\leq N_n$ to obtain a rooted metric space with surplus $N_n$ $$M^n_{g_n}=T^{n,\mathrm{mk}}(g_n)/\sim.$$ 
\end{enumerate}
\end{enumerate}
Then, all SCCs of $\vec{G}_n(\nu)$ are sub-digraphs of $\left\{M^n_{G_i^n}, i\geq 1 \right\}$. Call the kernels of these SCCs, ordered by decreasing size, $(C_i(n),i\geq 1)$, completed with an infinite repeat of $\mathfrak{L}$. Observe that we may view $M^n_{G_i^n}$ as a finite rooted directed multigraph $M^n_{G_i^n}$ whose edges are endowed with lengths. To be precise, in  $M^n_{G_i^n}$, let the vertex set consist of $G_i^n+1$, $W_i^n$ for $i\leq N_n$, and the branch points $V_i^n\wedge V_j^n$ for $i\neq j\leq N_n$. Then, we obtain $(C_i(n),i\geq 1)$ by ordering the kernels of the non-trivial SCCs in $\left\{M^n_{G_i^n}, i\geq 1 \right\}$ by decreasing size, and completing the list with an infinite repeat of $\mathfrak{L}$. See Figures \ref{figure.extractSCCs1}, \ref{figure.extractSCCs2} and \ref{figure.extractSCCs3} for an illustration of this procedure.

\subsection{The continuum case}\label{subsec.limitobject}

We will define now define the continuous counterpart of the sampling procedure of the out-forest and the candidates. This is a modification of the procedure defined in Subsection 3.2.2 of \cite{goldschmidtScalingLimitCritical2019}.

\subsubsection{\texorpdfstring{$\R$}{R}-trees and their encoding}
The continuum analogue of discrete trees are given by $\R$-trees. We give the basic definitions here and refer the reader to the survey paper \cite{legallRandomTreesApplications2005} for more details. An $\R$-tree is a compact metric space $(\cT, d)$ such that for every $a, b \in \cT$ the following two properties hold:
\begin{enumerate}
    \item There exists a unique isometry $$i_{a, b} : [0, d(a, b)] \to \cT$$ such that $i_{a, b}(0) = a$ and $i_{a, b}(d(a, b)) = b$.
    \item If $q: [0, 1] \to \cT$ is any continuous map such that $q(0) = a$ and $q(1) = b$ then the image of $q$ is the same as the image of $i_{a, b}$.
\end{enumerate}
Let $\pathbtw{a, b}$ denote the image of $i_{a, b}$. This is the unique path between $a$ and $b$.

$\R$-trees are often encoded by continuous excursions which can be seen as a continuous analogue of the height function of a tree. Let $f: [0, \sigma] \to [0, \infty)$ be a continuous excursion, meaning $f$ is continuous, $f(0) = f(\sigma) = 0$ and $f(x) > 0$ for all $x \in (0, \sigma)$. Using $f$ we can define a pseudo-metric
\begin{equation*}
    d_f(x, y) = f(x) + f(y) - 2 \min_{s \in [x \wedge y, x \vee y]} f(s).
\end{equation*}
This allows us to define the quotient space
\begin{equation*}
    \cT_f = [0, \sigma] / \{d_f = 0\}.
\end{equation*}
The space $\cT_f$ equipped with the metric $d_f$ is the \emph{$\R$-tree encoded by the excursion $f$}. Let $p_f: [0, \sigma] \to \cT_f$ be the natural projection function. Then $\cT_f$ inherits a distinguished root vertex $\rho = p(0) = p(\sigma)$.

A sequence of $\R$-trees is referred to as an $\R$-forest.

\subsubsection{The limit object}\label{subsubsec.samplecontinuousobject}
Let $(B_t,t\geq 0)$ be a standard Brownian motion, and set $$\left(\hat{B}_t,t\geq 0\right)=\left(B_t-\frac{\sigma_{-+}+\nu_-}{2\sigma_+\mu}t^2,t\geq 0\right).$$ 
\begin{remark}
We note that the coefficient of the parabolic drift of $\hat{B}$ is negative. Indeed, by definition of $\sigma_{+-}$ and $\nu_-$, the sign of the parabolic drift is the same as the sign of $\mu-\E[(D^-)^2D^+]$, and we note that
$$\frac{\E[(D^-)^2D^+]}{\E[D^+]}-\left(\frac{\E[D^+D^-]}{\E[D^+]}\right)^2=\frac{\E[(D^-)^2D^+]}{\mu}-1$$
is the variance of $D^-$ under the law of $\mathbf{D}$ size-biased by $D^+$, which is positive. Hence $\E[D^+(D^-)^2]/\mu\geq 1$, and the claimed negativity follows. 

\end{remark}
Define the reflected process
$$(\hat{R}_t,t\geq 0)= \left(\hat{B}_t-\inf\left\{\hat{B}_s: s\leq t\right\},t\geq 0\right).$$
Then, it follows from the argument in Section \ref{sec.convoutforest} that $\left(\frac{2}{\sigma_+}\hat{R}_t,t\geq 0\right)$ is the height process corresponding to an $\R$-forest with \L ukasiewicz path $\left(\sigma_+\hat{B}_t,t\geq 0\right)$. Call this forest the out-$\R$-forest. \\
Conditionally on $\hat{R}$, let $(A_t,t\geq 0)$ be a Cox process of intensity $$\frac{2(\sigma_{-+}+\nu_-)}{\sigma_+\mu^2} \hat{R}_t$$ at time $t$. Then, for $i$ in $\left\{1,2, \dots \right\}$, set $X_i=\min\{t:A_t=i\}$. For $i$ in $\left\{1,2,\dots \right\}$, define
\begin{align*}
G_i&=\inf\left\{t\geq 0:\hat{B}_t=\inf\{\hat{B}_s:s\leq X_i\}\right\}\text{ and}\\
\Sigma_i&=\inf\left\{ t\geq 0: \inf\{\hat{B}_s:s\leq G_i+t\} < \inf\{\hat{B}_s:s\leq X_i\}\right\},
\end{align*}
so that for each $i$ in $\left\{1,2,\dots \right\}$, $\left(\frac{2}{\sigma_+}\hat{R}_t,t\in [G_i,G_i+\Sigma_i]\right)$ encodes the $\R$-tree in the out-$\R$-forest that contains $X_i$. For each element of $\{(G_i,\Sigma_i):i=1,2,\dots\}$ we will sample the candidates in the $\R$-tree. Fix $i$, and set $(g,\sigma)=(G_i,\Sigma_i)$. Let $V_1=\inf\{s>0:A(s)=A(g)+1\}$, so that $g\leq V_1\leq g+\sigma$ by definition of $(g,\sigma)$. Let $\cT_g$ be the $\R$-tree encoded by $\left(\frac{2}{\sigma_+}\hat{R}_t,t\in [g,g+\sigma]\right)$ and let $p_g:[g,g+\sigma]\to \cT_l$ be the projection onto $\cT_g$ given by the encoding. Set $$||\cT_g||=\sup\left\{\frac{2}{\sigma_+}\hat{R}_t,t\in [g,g+\sigma]\right\},$$
the \emph{height} of $\cT_g$. \\
Suppose we have found candidates $\{V_1,\dots,V_m\}$. For $V_m\leq s\leq g+\sigma$, let $T^{\mathrm{mk}}_s$ be the subtree of $\cT_g$ spanned by $p_g\left(\{g,V_1,\dots,V_m,s\}\right)$, and let $|T^{\mathrm{mk}}_s|$ be its total length. Then, let $V_{m+1}$ be the first arrival time of a Poisson process on $[V_m,g+\sigma]$ of intensity $$\frac{\sigma_{-+}+\nu_-}{\mu^2}|T^{\mathrm{mk}}_s|ds.$$ If the process does not contain a point, let $\{V_1,\dots,V_m\}$ be the candidates of $\cT_g$, and set $N=m$. Otherwise, we repeat the inductive step for $\{V_1,\dots,V_{m+1}\}.$ If the induction does not terminate, we set $N=\infty$.\\
We show that $\P(N=\infty)=0$, by adapting the argument in Subsection 3.2.2 of \cite{goldschmidtScalingLimitCritical2019} to our set-up. Indeed, note that $V_m\leq s\leq V_{m+1}$ implies that  $|T^{\mathrm{mk}}_s|<(m+1)||\cT_g||$. Therefore, 
$$\P\left(\left.N\geq g+1,V_{m+1}-V_m<t \right|N\geq g\right)\leq \P(E_{m+1}<t),$$
for $(E_{k},k\geq 1)$ a sequence of exponential random variables with respective rates $$\frac{\sigma_{-+}+\nu_-}{\mu^2}k||\cT_g||.$$ 
Then,
$$\P\left(N=\infty \right)=\P\left(N=\infty\text{ and }\sup\{V_i:i\in \N\}<g+\sigma\right)\leq \P\left(\sum_{i=2}^\infty E_k\leq g+\sigma-V_1\right).$$
However, $\sum_{i=2}^\infty E_k=\infty$ a.s., because the harmonic series diverges, so, indeed, $\P\left(N<\infty \right)=1$. \\
Finally, for $1\leq i \leq N$, let the head corresponding to $V_i$, which we call $W_i$, be a uniform pick from the length measure on $T^{\mathrm{mk}}_{V_i}$. \\
Let $T^{\mathrm{mk}}(g)$ be the subtree of $\cT^{g}$ spanned by $\{g,V_1,\dots,V_N\}$. Then quotient $T^{\mathrm{mk}}(g)$ by the equivalence relation $\sim$ which identifies  $V_i$ and $W_i$ for each $1\leq i\leq N$ to obtain a rooted metric space $$\cM_{g}:=T^{\mathrm{mk}}(g)/\sim.$$ 
View $\cM_{g}$ as an element of $\vec{\cG}$ in the natural way. To be precise,  let the vertex set of $\cM_l$ consist of $g$, $W_i$ for $i\leq N$, and the branch points $V_i\wedge V_j$ for $i\neq j\leq N$. The directions are inherited from $\cT^l$, by considering all edges directed away from the root. Remove all edges that do not lie in an SCC of $\cM_{g}$ and delete any isolated vertices that are thus created. Then, apply the smoothing operation as defined in Subsection \ref{subsec.mdmkernels}. This creates a collection $\cC_g$ of strongly connected MDMs. Doing this for each $(g,\sigma)\in \{[G_i,\Sigma_i]\}$ yields the collection of strongly connected MDMs $\cC$ that has the law of the limit in Theorem \ref{thm.main}.

\subsubsection{Properties of the limit object}
We note that the limit object is encoded by $3$ parameters: the out-$\R$-forest is encoded by a Brownian motion with variance $\sigma_+^2$ and parabolic drift with coefficient $-(\sigma_{-+}+\nu_-)/(2\mu)$, and the identifications are a Cox process with intensity $(\sigma_{-+}+\nu_-)/\mu^2$ on the length measure of the subtree spanned by the previously found candidates and the currently explored vertex as described in Subsection \ref{subsubsec.samplecontinuousobject}. The limit object that is studied in \cite{goldschmidtScalingLimitCritical2019} corresponding to $\lambda=0$ (i.e. at criticality) is equal to our limit object in the case $\sigma_+^2=1$, $-(\sigma_{-+}+\nu_-)/(2\mu)=-1/2$, and $(\sigma_{-+}+\nu_-)/(\mu^2)=1$. Note that these three conditions are satisfied if we let $D^-$ and $D^-$ be independent $\operatorname{Poisson}(1)$ random variables. In \cite{goldschmidtScalingLimitCritical2019}, some properties of the limit object corresponding to these specific parameters are shown. A quick check shows that the proofs do not depend on the values of the parameters, so we deduce that the same properties also hold for our limit object. Let $\cM:=\bigcup_{G_i}\cM_{G_i}$.

\begin{proposition}
\begin{enumerate}
    \item The number of complex connected components of $\cM$ has finite expectation.
    \item The number of loops of $\cM$ is a.s. infinite.
\end{enumerate}
\end{proposition}

\begin{proposition}
\label{prop:allengthsaredifferent}
The SCCs of $\cM$ all have different lengths almost surely.
\end{proposition}
Write $\cC$ for the SCCs of $\cM$ and $\mathbf{C}_l$ for those of $\cM_l$, in decreasing order of length, with $\cM_l$ as defined in Subsection  \ref{subsubsec.samplecontinuousobject}. Write $\cC_{\text{cplx}}$ for the list of complex components of $\cC$ in decreasing order of length. For sequences $(K_1,\dots, K_j)$ and $(J_1,\dots,J_k)$ of directed multigraphs, write $(J_1,\dots,J_k)\equiv(K_1,\dots, K_j)$ if $j=k$ and $J_i$ is isomorphic to $K_i$ for each $i\leq j$. Extend this notation naturally to the case where one or both of the sequences has edge lengths by ignoring the edge lengths. 
\begin{theorem}
Let $K_1,\dots, K_j$ be a finite sequence consisting of $3$-regular strongly connected directed multigraphs or loops. We have 
$$\P\left[\cC_l\equiv(K_1,\dots, K_j)\right]>0.$$
Assuming that $K_1,\dots, K_j$ are all complex, we also have that 
$$\P\left[ \cC_{\text{cplx}}\equiv(K_1,\dots, K_j)\right]>0.$$
Let $(e_i,1\leq i \leq M)$ be an arbitrary ordering of the edges of $K_1,\dots, K_j$. Then, conditionally on $\cC_l\equiv(K_1,\dots, K_j)$, (resp. $\cC_{\text{cplx}}\equiv(K_1,\dots, K_j)$), $\cC_l$ (resp. $\cC_{\text{cplx}}$) gives lengths $(\ell(e_i),1\leq i \leq M)$ to these edges, and their joint distribution has full support in 
$$\left\{\mathbf{x}=(x_1,\dots, x_M)\in \R_+^M:\forall 1\leq i\leq j-1, \sum_{k:e_k\in E(K_i)}x_k \geq \sum_{k:e_k\in E(K_{i+1})}x_k\right\}.$$
\end{theorem}

\section{Analysis of the measure change}
\label{sec:measure-change}

Recall $(\hat{\vD}_{n, 1}, \hat{\vD}_{n, 2}, \ldots, \hat{\vD}_{n, n})$ are are the degree pairs of the vertices in order of discovery, and that $R_n$ is the number of vertices with positive in-degree. The behaviour of the $\hat{\vD}_{m, n}$ for $m \leq R_n$ and $m > R_n$ is rather different. Before $R_n$, new vertices are discovered with probability proportional to their in-degree. After $R_n$, all vertices with positive in-degree have already been discovered and we choose to explore the remaining vertices in some uniform order.

Later in \cref{sec.convSCCs}, we show that we only need to consider timescales of the order of $m = \Theta(n^{2/3})$. Let $p = \P(D^- > 0)$ such that $R_n$ is distributed as a $\text{Binomial}(n, p)$ random variable. We show now that the probability that $m \leq R_n$ will converge exponentially to 1.
\begin{lemma}
    If $m = \Theta(n^{2/3})$ then there exists $c > 0$ such that $\P(R_n < m) < e^{-c n}$.
\end{lemma}
\begin{proof}
    If $m = \Theta(n^{2/3})$ then $\E[R_n] - m = pn - m = \Theta(n)$. Thus by Hoeffding's inequality
    \begin{equation*}
        \P(R_n < m)
        = \P\left(\E[R_n] - R_n > \E[R_n] - m\right)
        \leq e^{-\frac{2}{n} (\E[R_n] - m)^2} < e^{-cn}
    \end{equation*}
    for some $c > 0$.
\end{proof}
Hence it is sensible to prove results on the event that $m \leq R_n$.

When discussing the criticality condition, we gave heuristics showing that the limiting distribution of $\hat{D}_{n, 1}$ is given by $\vZ$ where
\begin{equation*}
    \P(Z^- = k^-, Z^+ = k^+) = \frac{k^-}{\mu}\P(D^- = k^-, D^+ = k^+).
\end{equation*}
Similarly, $\hat{D}_{n, 2}$ is also approximately distributed like $\vZ$ for large $n$, and so on. In this section we in fact prove a precise relation between $\hat{\vD}_{n, 1}, \ldots, \hat{\vD}_{n, m}$ and the sequence $\vZ_1, \vZ_2, \ldots$ of i.i.d.\ copies of $\vZ$.

The results proved in this section do not actually require the criticality condition, so let us define notation for the mean of the $Z_i^{\pm}$ and the two corresponding centered random walks. Let
\begin{equation*}
    \lambda_{\pm} = \E[Z^{\pm}_1]
    \quad \text{and} \quad
    V^{\pm}(n) = \textstyle \sum_{i=1}^n (Z^{\pm}_i - \lambda_{\pm}).
\end{equation*}
The criticality condition is then equivalent to assuming $\lambda_+ = 1$. We also define the notation
\begin{equation*}
    \Xi^{\pm}_{n-m} = \sum_{i=m+1}^n D_i^{\pm}
    \quad \text{and} \quad
    \Delta_n = \Xi^-_n - \Xi^+_n
\end{equation*}
such that $\{\Delta_n = 0\}$ is the event that the total out-degree is equal to the total in-degree.

The following lemma asserts the existence of the measure change $\phi^n_m$, and its joint scaling limit with the random walks $V^-$ and $V^+$ when $m = \floor{n^{2/3} T}$ for some $T > 0$.
\begin{proposition}
    \label{prop:measure-change-no-crit}
    For all positive integers $n$ and $m$ such that $m \leq n$, there exists a function $\phi^n_m: (\N \times \N)^n \to [0, \infty)$ such that 
    \begin{equation*}
        \E\left[\left.u(\hat{\vD}_{n, 1}, \ldots, \hat{\vD}_{n, m}) \one\{R_n \geq m\} \right| \Delta_n = 0\right]
        = \E[u(\vZ_1, \ldots, \vZ_m) \phi^n_m(\vZ_1, \ldots, \vZ_m)]
    \end{equation*}
    for all bounded test functions $u: (\N \times \N)^m \to \R$. Define
    \begin{equation*}
        \Phi(n, m) = \phi^n_m(\vZ_1, \ldots, \vZ_m).
    \end{equation*}
    Further, let $(W^-, W^+)$ be a pair of correlated standard Brownian motions with correlation $\corr(Z_1^-, Z_1^+)$ and, for $T > 0$, define
    \begin{equation*}
        \Phi(T) = \exp\left( 
            - \frac{\sigma_-}{\mu} \int_0^T s \dif W_s^- 
            - \frac{\sigma_-^2}{6 \mu^2} T^3
        \right).
    \end{equation*}
    Then for all $T > 0$,
    \begin{align*}
        &\left( 
            \Phi(n, \floor{n^{2/3} T}),
            \left(
                n^{-1/3} V^-\left( \floor{n^{2/3} t} \right),
                n^{-1/3} V^+\left( \floor{n^{2/3} t} \right)
            \right)_{t \in [0, T]}
        \right) \\
        & \hspace{23em} \todist \left( 
            \Phi(T),
            (\sigma_- W^-_t, \sigma_+ W^+_t)_{t \in [0, T]}
        \right)
    \end{align*}
    in $\R \times \D([0, T], \R^2)$ as $n \to \infty$, even in the absence of the criticality condition.
\end{proposition}
The rest of this section is dedicated to proving this proposition.

\subsection{Exact form of the measure change}

To determine the exact form of the measure change, we first need to know the law of the ordering of the first $R_n$ vertices. Let $\cI_n = \{i \in [n] : D_i^- > 0\}$. The first $R_n$ vertices we explore in \cref{alg:edfs} will have positive in-degree, thus there exists a random bijection $\Sigma_n: [R_n] \to \cI_n$ such that $\hat{\vD}_{n, i} = \vD_{\Sigma_n(i)}$ for $i = 1, \ldots, R_n$. 
\begin{lemma}
    $\Sigma_n$ has law
    \begin{equation*}
        \P(\Sigma_n = \sigma \mid \vD_1, \ldots, \vD_n)
        = \prod_{i=1}^{R_n} \frac{D^-_{\sigma(i)}}{\sum_{j=i}^{R_n} D^-_{\sigma(j)}}.
    \end{equation*}
    for all bijections $\sigma: [R_n] \to \cI_n$.
\end{lemma}
\begin{proof}
    In \cref{alg:edfs}, a vertex first becomes explored in two ways. Either it is at the start of an out-component or it is discovered when an out-half-edge is paired to one its in-half-edges.

    Suppose we have explored $m$ vertices and $m < R_n$. If the next vertex is explored by pairing one of its in-half-edges, then we have chosen it with probability proportional to its in-degree since in- and out-half-edges are paired uniformly at random. Otherwise, it is at the start of a new out-component, and since $m < R_n$, there are still vertices of positive in-degree. Thus we still pick a new vertex with probability proportional to its in-degree.

    Therefore in all cases,
    \begin{align*}
        &\P(\Sigma_n(m+1) = \sigma(m+1) \mid \Sigma_n(1) = \sigma(1), \ldots, \Sigma_n(m) = \sigma(m), \vD_1, \ldots, \vD_n) \\
        =&  \frac{D^-_{\sigma(m+1)}}{\sum_{i \in \cI_n} D_i^- - \sum_{j=1}^m D^-_{\sigma(j)}}
        = \frac{D^-_{\sigma(m+1)}}{\sum_{j=m+1}^{R_n} D^-_{\sigma(j)}}.
    \end{align*} 
    From this, repeated applications of the definition of conditional probability yields the desired result.
\end{proof}

Next we establish the form of the measure change when we condition on the exact value of $R_n$ but not $\Delta_n = 0$.

\begin{lemma}
    \label{lem:exact-measure-change-no-conditioning}
    For all integers $0 \leq r \leq n$ and test functions $u: (\N \times \N)^r \times \N \to \R$,
    \begin{equation*}
        \E \left[ u\left( 
            \hat{\vD}_{n, 1}, \ldots, \hat{\vD}_{n, r}, \textstyle \sum_{i \in \cI_n^c} D_i^+
        \right) \,\middle\vert\, R_n = r \right] \\
        =
        \E \left[
            u\left( \vZ_1, \ldots, \vZ_r, \textstyle \sum_{i = 1}^{n-r} E_i^+ \right)
            \psi_r(\vZ_1, \ldots, \vZ_r)
        \right]
    \end{equation*}
    where
    \begin{equation*}
        \psi_r(\vk_1, \ldots, \vk_r) =
        \frac{1}{p^r} \prod_{i=1}^r \frac{(r - i + 1)\mu}{\sum_{j=i}^r k_i^-}.
    \end{equation*}
    and $E_1^+, E_2^+, \ldots$ are i.i.d.\ random variables such that $E_i^+$ has the same distribution as $D^+$ conditioned on $D^- = 0$. We take the sequences $(E_i^+)_{i \geq 1}$ and $(\vZ_i)_{i \geq 1}$ to be independent.
\end{lemma}

\begin{proof}
    For any $\vk_1, \ldots, \vk_m \in \N^+ \times \N$ for all $i$ and $s \in \N$.
    \begin{align*}
        & \P\left(\hat{\vD}_{n, 1} = \vk_1, \ldots, \hat{\vD}_{n, r} = \vk_r, \textstyle \sum_{i \in \cI_n^c} D_i^+ = s, R_n = r \right) \\
        =& \sum_{\substack{I \subset [n] \\ \abs{I} = r}} \sum_{\sigma: [r] \to I}
        \P\left(\vD_{\Sigma_n(1)} = \vk_1, \ldots, \vD_{\Sigma_n(r)} = \vk_r, \textstyle \sum_{i \in \cI_n^c} D_i^+ = s, \cI_n = I, \Sigma_n = \sigma \right)
    \end{align*}
    where the second summation is taken over all bijections $\sigma: [r] \to I$. We examine a single summand.
    \begin{align*}
        &\P\left(\vD_{\Sigma_n(1)} = \vk_1, \ldots, \vD_{\Sigma_n(r)} = \vk_r, \textstyle \sum_{i \in \cI_n^c} D_i^+ = s, \cI_n = I, \Sigma_n = \sigma \right) \\
        =&\P\left(
            \vD_{\sigma(j)} = \vk_j\ \text{for $j = 1, \ldots, r$},
            \textstyle \sum_{i \in I^c} D_i^+ = s,
            D^-_i = 0\ \text{for $i \in I^c$},
            \Sigma_n = \sigma
            \right)  \\
        =& \prod_{i=1}^r \frac{k_i^-}{\sum_{j=i}^r k_j^-}
        \times \prod_{i=1}^r \lambda_{\vk_i}
        \times \P\left( 
            \textstyle \sum_{i \in I^c} D_i^+ = s,
            D^-_i = 0\ \text{for $i \in I^c$}
         \right).
    \end{align*}
    where $\lambda_{\vk} = \P(\vD_1 = \vk)$. We have
    \begin{equation*}
        \P\left( 
            \textstyle \sum_{i \in I^c} D_i^+ = s,
            D^-_i = 0\ \text{for $i \in I^c$}
         \right)
         = (1-p)^{n-r} \P\left( 
             \textstyle \sum_{i=1}^{n-r} E^+_i = s
          \right).
    \end{equation*}
    Also
    \begin{align*}
        \prod_{i=1}^r \frac{k_i^-}{\sum_{j=i}^r k_j^-} \times \prod_{i=1}^r \lambda_{\vk_i}
        &= \prod_{i=1}^r \frac{k_i^-}{\mu} \lambda_{\vk_i} \times \prod_{i=1}^r \frac{\mu}{\sum_{j=i}^r k_j^-} \\
        &= \P(\vZ_1 = \vk_1, \ldots, \vZ_r = \vk_r)
        \times \prod_{i=1}^r \frac{\mu}{\sum_{j=i}^r k_j^-}.
    \end{align*}
    Therefore
    \begin{align*}
        & \P\left(\hat{\vD}_{n, 1} = \vk_1, \ldots, \hat{\vD}_{n, r} = \vk_r, \textstyle \sum_{i \in \cI_n^c} D_i^+ = s, R_n = r \right) \\
        =& \binom{n}{r} \times r!
        \times \prod_{i=1}^r \frac{\mu}{\sum_{j=i}^r k_j^-} \times (1-p)^{n-r}
        \times \P\left(\vZ_1 = \vk_1, \ldots, \vZ_r = \vk_r, \textstyle \sum_{i=1}^{n-r} E_i^+ = s \right) \\
        =& \binom{n}{r} p^r (1-p)^{n-r} \times \frac{1}{p^r} \prod_{i=1}^r \frac{(r-i+1) \mu}{\sum_{j=i}^r k_i^-}
        \times \P\left(\vZ_1 = \vk_1, \ldots, \vZ_r = \vk_r, \textstyle \sum_{i=1}^{n-r} E_i^+ = s \right).
    \end{align*}
    Finally dividing by $\P(R_n = r) = \binom{n}{r} p^r (1-p)^{n-r}$ gives the desired measure change.
\end{proof}

Using the previous lemma we can prove existence and give the exact form of the desired measure change $\phi^n_m$.

\begin{lemma}
    \label{lem:exact-measure-change}
    For all $m \leq n$ and test functions $u: (\N \times \N)^m \to \R$,
    \begin{equation*}
        \E \left[
            u\left( \hat{\vD}_{n, 1}, \ldots, \hat{\vD}_{n, m} \right)
            \one\{R_n \geq m\}
            \,\middle\vert\,
            \Delta_n = 0
        \right] \\
        =
        \E \left[
            u\left( \vZ_1, \ldots, \vZ_m \right)
            \phi^n_m(\vZ_1, \ldots, \vZ_m)
        \right],
    \end{equation*}
    where
    \begin{align*}
        \phi^n_m(\vk_1, \ldots, \vk_m) =
        \frac{1}{\P(\Delta_n = 0)} \E\left[ 
            \one \left\{ \Delta_{n-m} = \sum_{i=1}^m (k_i^+ - k_i^-) \right\}
            \prod_{i=1}^m \frac{(n - i + 1) \mu}{\sum_{j=1}^m k_j^- + \Xi_{n-m}^-}
        \right].
    \end{align*}
\end{lemma}

\begin{proof}
    By \cref{lem:exact-measure-change-no-conditioning}, for all $r \geq m$
    \begin{align*}
        &\E \left[ 
            u\left( 
                \hat{\vD}_{n, 1}, \ldots, \hat{\vD}_{n, m}
             \right)
            \one\{\Delta_n = 0\}
            \, \middle\vert \,
            R_n = r
         \right] \\
        =&\E \left[ 
            u\left( \vZ_1, \ldots, \vZ_m \right)
            \one \left\{ 
                \sum_{i=1}^r (Z_i^- - Z_i^+) - \sum_{i=1}^{n-r} E_i^+ = 0
             \right\}
             \frac{1}{p^r} \prod_{i=1}^r \frac{(r - i + 1) \mu}{\sum_{j=i}^r Z_j^-}
         \right] \\
        =&\E \left[ 
            u\left( \vZ_1, \ldots, \vZ_m \right)
            \E \left[ 
                \one \left\{ 
                    \sum_{i=1}^r (Z_i^- - Z_i^+) - \sum_{i=1}^{n-r} E_i^+ = 0
                \right\}
                \frac{1}{p^r} \prod_{i=1}^r \frac{(r - i + 1) \mu}{\sum_{j=i}^r Z_j^-}
                \, \middle \lvert \,
                \vZ_1, \ldots, \vZ_m
             \right]
         \right] \\
        =&\E \left[ 
            u\left( \vZ_1, \ldots, \vZ_m \right)
            \tilde{\gamma}^{n, m}_r (\vZ_1, \ldots, \vZ_m)
         \right],
    \end{align*}
    where
    \begin{align*}
        \tilde{\gamma}^{n, m}_r (\vk_1, \ldots, \vk_m)
        &= \E \Bigg[ 
            \one \left\{ 
                \sum_{i=m+1}^r (Z_i^- - Z_i^+) - \sum_{i=1}^{n-r} E_i^+ = \sum_{i=1}^m (k_i^+ - k_i^-)
            \right\} \times \\
            &\hspace{5em}
            \frac{1}{p^m} \prod_{i=1}^m \frac{(r - i + 1) \mu}{\sum_{j=i}^m k_j^- + \sum_{j=m+1}^r Z_j^-}
            \frac{1}{p^{m-r}} \prod_{i=m+1}^r \frac{(r - i + 1) \mu}{\sum_{j=i}^r Z_j^-}
        \Bigg] \\
        &= \E \Bigg[ 
            \one \left\{ 
                \sum_{i=1}^{r-m} (Z_i^- - Z_i^+) - \sum_{i=1}^{n-r} E_i^+ = \sum_{i=1}^m (k_i^+ - k_i^-)
            \right\} \times \\
            &\hspace{5em}
            \frac{1}{p^m} \prod_{i=1}^m \frac{(r - i + 1) \mu}{\sum_{j=i}^m k_j^- + \sum_{j=1}^{r-m} Z_j^-}
            \frac{1}{p^{m-r}} \prod_{i=1}^{r-m} \frac{(r - m - i + 1) \mu}{\sum_{j=i}^{r-m} Z_j^-}
        \Bigg],
    \end{align*}
    since $(\vZ_i)_{i=m+1}^r$ has the same law as $(\vZ_i)_{i=1}^{r-m}$. Then applying \cref{lem:exact-measure-change-no-conditioning} again shows that
    \begin{align*}
        \tilde{\gamma}^{n, m}_r (\vk_1, \ldots, \vk_m)
        &= \E \Bigg[ 
            \one \left\{ 
                \sum_{i=1}^{r-m} (\hat{D}_{n-m,i}^- - \hat{D}_{n-m,i}^+) - \sum_{i \in \cI_{n-m}^c} D_i^+ = \sum_{i=1}^m (k_i^+ - k_i^-)
            \right\} \times \\
            &\hspace{8em}
            \frac{1}{p^m} \prod_{i=1}^m \frac{(r - i + 1) \mu}{\sum_{j=i}^m k_j^- + \sum_{j=1}^{r-m} \hat{D}_{n-m,j}^-}
            \Biggm|
            R_{n-m} = r-m
        \Bigg].
    \end{align*}
    Conditional on $R_{n-m} = r-m$, we have
    \begin{equation*}
        \sum_{j=1}^{r-m} (\hat{D}_{n-m,j}^- - \hat{D}_{n-m,j}^+) - \sum_{i \in \cI_{n-m}^c} D_i^+ = \Delta_{n-m}
        \quad \text{and} \quad
        \sum_{j=1}^{r-m} \hat{D}_{n-m, j}^- = \Xi^-_{n-m}.
    \end{equation*}
    Therefore,
    \begin{align*}
        \tilde{\gamma}^{n, m}_r (\vk_1, \ldots, \vk_m)
        &= \E \left[
            \frac{1}{p^m} \prod_{i=1}^m \frac{(r - i + 1) \mu}{\sum_{j=i}^m k_j^- + \Xi_{n-m}^-} \one_{A_n}
            \, \middle \vert \,
            R_{n-m} = r-m
        \right],
    \end{align*}
    where
    \begin{equation*}
        A_n(\vk_1, \ldots, \vk_m) = \left\{ \Delta_{n-m} =  \sum_{i=1}^m (k_i^+ - k_i^-) \right\}.
    \end{equation*}
    Hence,
    \begin{align*}
        \E \left[ 
            u\left( 
                \hat{\vD}_{n, 1}, \ldots, \hat{\vD}_{n, m}
             \right)
            \one\{R_n \geq m, \Delta_n = 0\}
         \right]
        =\E \left[ 
            u\left( \vZ_1, \ldots, \vZ_m \right)
            \tilde{\phi}^n_m(\vZ_1, \ldots, \vZ_m)
         \right],
    \end{align*}
    where
    \begin{align*}
        &\tilde{\phi}^n_m(\vk_1, \ldots, \vk_m) \\
        =& \sum_{r=m}^n \binom{n}{r} p^r (1-p)^{n-r} 
        \E \left[
            \frac{1}{p^m} \prod_{i=1}^m \frac{(r - i + 1) \mu}{\sum_{j=i}^m k_j^- + \Xi_{n-m}^-} \one_{A_n}
            \, \middle \vert \,
            R_{n-m} = r-m
        \right] \\
        =& \sum_{l=1}^{n-m} \binom{n}{l+m} p^{l+m} (1-p)^{n-m-l} 
        \E \left[
            \frac{1}{p^m} \prod_{i=1}^m \frac{(l + m - i + 1) \mu}{\sum_{j=i}^m k_j^- + \Xi_{n-m}^-} \one_{A_n}
            \, \middle \vert \,
            R_{n-m} = l
        \right].
    \end{align*}
    We wish to view the sum as an expectation over $R_{n-m}$. In order to do this, we rewrite the expression so that we are taking a sum over the probabilities of a $\text{Binomial}(n-m, p)$ distribution. We can calculate
    \begin{equation*}
        \frac{\binom{n}{l+m} p^{l+m} (1-p)^{n-m-l}}{\binom{n-m}{l} p^l (1-p)^{n-m-l}}
        = p^m \prod_{i=1}^m \frac{(n-i+1)}{(l+m-i+1)}.
    \end{equation*}
    Therefore,
    \begin{align*}
        \tilde{\phi}^n_m(\vk_1, \ldots, \vk_m)
        &= \sum_{l=1}^{n-m} \binom{n-m}{l} p^l (1-p)^{n-m-l}
        \E \left[
            \prod_{i=1}^m \frac{(n - i + 1) \mu}{\sum_{j=i}^m k_j^- + \Xi_{n-m}^-} \one_{A_n}
            \, \middle \vert \,
            R_{n-m} = l
        \right] \\
        &= \E\left[ 
            \E \left[
                \prod_{i=1}^m \frac{(n - i + 1) \mu}{\sum_{j=i}^m k_j^- + \Xi_{n-m}^-} \one_{A_n}
                \, \middle \vert \,
                R_{n-m}
            \right]
         \right] \\
         &= \E \left[
            \prod_{i=1}^m \frac{(n - i + 1) \mu}{\sum_{j=i}^m k_j^- + \Xi_{n-m}^-} \one_{A_n}
        \right].
    \end{align*}
    Finally, dividing by $\P(\Delta_n = 0)$ yields the desired form of $\phi^n_m$.
\end{proof}

\subsection{Asymptotic lower bound on the measure change}

Recall that our goal in \cref{prop:measure-change-no-crit} is to determine the limiting distribution of
\begin{equation*}
    \Phi(n, m) = \phi^n_m(\vZ_1, \ldots, \vZ_m),
\end{equation*}
as $n \to \infty$, in the regime where $m = \Theta(n^{2/3})$. When dealing with convergence in distribution, it is suffcient and necessary to work on a sequence of events occuring with high probability. In particular, for the proof of \cref{prop:measure-change-no-crit}, we work on the event $\walkDeviationEvent_m$ where
\begin{align*}
    \walkDeviationEvent_m
    &= \Bigg\{ 
    \max_{i=1, \ldots, m } \abs*{
        \textstyle\sum_{j=1}^i (Z^{-}_j - \lambda_{-}) 
    } \leq m^{1/2} \log(m) \\
    &\hspace{13em} \text{and} \quad
    \max_{i=1, \ldots, m } \abs*{
        \textstyle\sum_{j=1}^i (Z^{+}_j - \lambda_{+}) 
    } \leq m^{1/2} \log(m)
    \Bigg\}.
\end{align*}
This says that the centered random walks corresponding to $Z^+_i$ and $Z^-_i$ both do not deviate by more than $m^{1/2} \log(m)$ in the first $m$ steps. The conditions in \cref{subsec:model-description} ensure each $Z^+_i$ and $Z^-_i$ has finite variance, thus this event will occur with high probability. 

The following lemma is an analogue of \citet[Lemma 6.7]{conchon--kerjanStableGraphMetric2021}. In it we prove a deterministic lower bound on $\phi^n_m(\vk_1, \ldots, \vk_m)$, for all $\vk_1, \ldots, \vk_m$ corresponding to the event $\walkDeviationEvent_m$, up to an error which vanishes as $n \to \infty$.
\begin{proposition}
    \label{prop:measure-change-approx}
    Define
    \begin{equation*}
        s^{\pm}(i) = \textstyle{\sum_{j=1}^i (k_i^{\pm} - \lambda_{\pm})}.
    \end{equation*}
    Suppose that $\vk_1, \ldots, \vk_m$ are such that
    \begin{equation}
        \label{eq:s-condition}
        \max_{i=1, \ldots, m} \abs{s^{-}(i)} \leq m^{\frac{1}{2}} \log(m)
        \quad \text{and} \quad
        \max_{i=1, \ldots, m} \abs{s^{+}(i)} \leq m^{\frac{1}{2}} \log(m)
    \end{equation}
    Then in the regime $m = \Theta(n^{2/3})$, as $n \to \infty$,
    \begin{equation*}
        \phi^n_m(\vk_1, \ldots, \vk_m)
        \geq \exp\left( \frac{1}{\mu n} \sum_{i=0}^m (s^-(i) - s^-(m)) - \frac{\sigma_-}{6 \mu^2} \frac{m^3}{n^2} \right) + \littleo(1),
    \end{equation*}
    where the $\littleo(1)$ error term is independent of $\vk_1, \ldots, \vk_m$ satisfying the assumption in \cref{eq:s-condition}.
\end{proposition}

The fact that we only prove a lower bound may seem strange at first. To understand why this is sufficient, first note that all measure changes are non-negative random variables and have expectation 1. Hence if the sequence of lower bounds on the measure changes converge to a limit that also has expectation 1, then we have not have lost a significant amount of probability mass. It follows that the measure changes converge to the same limit as the lower bounds. This is made formal by \citet[Lemma 4.8]{conchon--kerjanStableGraphMetric2021}. In \cref{prop:measure-change-no-crit} we are considering the joint convergence of the measure change with two other random walks, and thus we adapt \cite[Lemma 4.8]{conchon--kerjanStableGraphMetric2021} to allow for an additional coordinate that is converging jointly with the first coordinate.
\begin{lemma}
    \label{lem:sandwiching-lemma}
    Let $(X_n, Y_n, Z_n)_{n \geq 1}$ be a sequence of $[0, \infty) \times [0, \infty) \times S$-valued random variables where $S$ is a metric space. Suppose there exists a $[0, \infty) \times S$-valued random variable $(Y, Z)$ such that the following holds:
    \begin{enumerate}
        \item $(Y_n, Z_n) \todist (Y, Z)$ as $n \to \infty$.
        \item $X_n \geq Y_n$ almost surely for all $n$.
        \item $\E[X_n] = 1$ for all $n$ and $\E[Y] = 1$.
    \end{enumerate}
    Then $(X_n, Z_n) \todist (Y, Z)$ also. Moreover $(X_n)_{n \geq 1}$ is a sequence of uniformly integrable random variables.
\end{lemma}
The proof of this lemma is obtained by simply adding the corresponding $Z_n$ or $Z$ coordinate to quantities in the proof of \cite[Lemma 4.8]{conchon--kerjanStableGraphMetric2021} and so we will not repeat it here.

\subsubsection{Discrete local limit theorem}

To prove \cref{prop:measure-change-approx}, we first need to understand the denominator of $\phi^n_m$, which, as given by \cref{lem:exact-measure-change}, is $\P(\Delta_n = 0)$. The random variable $\Delta_n$ is a sum of independent integer-valued random variables and the asymptotic behaviour of such a sum being equal to some value is described by the discrete local limit theorem. Such a theorem was first proven by \citet{gnedenkoLocalLimitTheorem1948}. Here we borrow the presentation from \citet[Section 3.5]{durrettProbabilityTheoryExamples2019}.

Let $X_1, X_2, \ldots$ be i.i.d.\ integer-valued random variables with mean $\mu$ and finite variance $\sigma^2$. Then, by the central limit theorem,
\begin{equation*}
    \frac{ \sum_{i=1}^n X_i - n \mu }{\sigma \sqrt{n}}
    \todist
    N(0, 1)
\end{equation*}
as $n \to \infty$. Thus we expect $\sum_{i=1}^n X_i$ to be distributed like a $N(n \mu, n \sigma^2)$ random variable for large values of $n$. Therefore the probability mass function of $\sum_{i=1}^n X_i$ should be well approximated by the probability density function of a $N(n \mu, n \sigma^2)$ distribution, i.e.
\begin{equation*}
    \P\big( \textstyle \sum_{i=1}^n X_i = s \big)
    \approx
    \frac{1}{\sqrt{2 \pi n \sigma^2}}  \exp \left(
        \frac{-(s - n \mu)^2}{2 n \sigma^2}
    \right)
\end{equation*}
for all integers $s$. Specifically we hope that
\begin{equation}
    \label{eq:llt-approx}
    \sup_{s \in \Z} \abs*{
        \P\big( \textstyle \sum_{i=1}^n X_i = s \big) 
        -
        \frac{1}{\sqrt{2 \pi n \sigma^2}}  \exp \left(
            \frac{-(s - n \mu)^2}{2 n \sigma^2}
        \right)
    } = \littleo(n^{-1/2}).
\end{equation}

This, however, is not always the case. Suppose, for example, that each $X_i$ is almost surely even such that $\P\left( \sum_{i=1}^n X_i = s \right) = 0$ for all odd $s$. Let $s_n$ be the closest odd integer to $n \mu$. Then
\begin{equation*}
    \sup_{s \in \Z} \abs*{
        \P\big( \textstyle \sum_{i=1}^n X_i = s \big) 
        -
        \frac{1}{\sqrt{2 \pi n \sigma^2}}  \exp \left(
            \frac{-(s - n \mu)^2}{2 n \sigma^2}
        \right)
    }
    \geq
    \frac{1}{\sqrt{2 \pi n \sigma^2}}  \exp \left(
        \frac{-(s_n - n \mu)^2}{2 n \sigma^2}
    \right)
    = 
    \Theta(n^{ -1/2 }).
\end{equation*}

Fortunately this kind of periodic behaviour can be mitigated by normalizing the random variables. A one-dimensional random variable $X$ is \emph{lattice} if it is not almost surely constant, and there exists $h > 0$ and $c \in \R$ such that $X \in c + h\Z$ almost surely. The largest such $h$ is called the \emph{span} of $X$. For example, if $X$ is almost surely even then $X$ has span at least 2. If $X$ is lattice with span $h$ and $c$ is in the support of $X$, then the affine transform $\frac{1}{h}(X - c)$ is an integer-valued random variable with span 1, for which it can be shown that the approximation in \cref{eq:llt-approx} does hold. This gives us the discrete local limit theorem:
\begin{theorem}[Discrete local limit theorem]
    \label{thm:discrete-llt}
    Let $X_1, X_2, \ldots$ be i.i.d.\ $\R$-valued lattice random variables with span $h$ and fix arbitrary $c \in \supp(X_1)$. Then
    \begin{equation*}
        \sup_{s \in n c + h \Z} \abs*{
            \P\big( \textstyle \sum_{i=1}^n X_i = s \big) 
            -
            \frac{h}{\sqrt{2 \pi n \sigma^2}}  \exp \left(
                \frac{-(s - n \mu)^2}{2 n \sigma^2}
            \right)
        } = \littleo(n^{-1/2}).
    \end{equation*}
\end{theorem}

\begin{remark}
    \label{rem:llt-limitations}
    For each sequence of integers $(s_n)_{n \geq 1}$ such that $\abs{s_n - n \mu} = \littleomega(n^{1/2})$, we have that
    \begin{equation*}
        \frac{1}{\sqrt{2 \pi n \sigma^2}}  \exp \left(
            \frac{-(s_n - n \mu)^2}{2 n \sigma^2}
        \right)
        = \littleo(n^{-1/2}).
    \end{equation*}
    Hence the discrete local limit theorem (\cref{thm:discrete-llt}) tells you only that $\P\left(\sum_{i=1}^n X_i = s_n\right) = \littleo(n^{-1/2})$; it gives no precise characterization of the leading order term.
\end{remark}

While this remark will be important later, here $\Delta_n$ is centered and we are interested in the probability $\P(\Delta_n = 0)$. In addition, the strong aperiodicity condition in \cref{subsec:model-description} tells us exactly that the $D^- - D^+$ is lattice with span 1. Thus the following is a direct corollary of the discrete local limit theorem (\cref{thm:discrete-llt}).
\begin{corollary}
    \label{cor:measure-change-denominator-control}
    We have
    \begin{equation*}
        \P(\Delta_n = 0) = \frac{1}{\sqrt{2 \pi \sigma^2 n}} + \littleo(n^{-1/2})
    \end{equation*}
    as $n \to \infty$, where $\sigma$ is the variance of $D^- - D^+$.
\end{corollary}
\begin{remark}
    The exact value of $\sigma^2$ is not important for the asymptotic behaviour of $\phi^n_m$ because we show later that it will cancel with a term in the numerator of $\phi^n_m$.
\end{remark}

\subsubsection{Exponential Tilting}

Next we turn to the numerator of $\phi^n_m$. By \cref{lem:exact-measure-change}, this is given by
\begin{equation}
    \label{eq:measure-change-numerator}
    \E\left[ 
        \one \left\{ \Delta_{n-m} = \sum_{i=1}^m (k_i^+ - k_i^-) \right\}
        \prod_{i=1}^m \frac{(n - i + 1) \mu}{\sum_{j=i}^m k_j^- + \Xi_{n-m}^-}
    \right].
\end{equation}
For convenience, let $\lltEvent_n$ denote the event in the indicator function, i.e.
\begin{equation*}
    \lltEvent_n = \left\{ \Delta_{n-m} = \sum_{i=1}^m (k_i^+ - k_i^-) \right\}.
\end{equation*}
When the $\vk_1, \ldots, \vk_m$ satisfy the condition in \cref{eq:s-condition}, we face two problems in evaluating the expectation in \cref{eq:measure-change-numerator}.

The first problem concerns the event $\lltEvent_n$. To evaluate the expectation we need to understand the asymptotic probability of this event. Unfortunately a naïve application of the discrete local limit theorem will not work in this case, as we now explain. Firstly, note that
\begin{equation*}
    \sum_{i=1}^m (k_i^- - k_i^+) = s^-(m) - s^+(m) + (\lambda_+ - \lambda_-) m.
\end{equation*}
We have that
\begin{equation*}
    \lambda_+ - \lambda_- = \E[Z^- - Z^+] = \tfrac{1}{\mu} \E[D^-D^+ - (D^-)^2]
\end{equation*}
which is, in general, non-zero. Then $m = \Theta(n^{2/3})$ whereas, if $\vk_1, \ldots, \vk_n$ satisfy \cref{eq:s-condition}, $s^{-}(m)$ and $s^+(m)$ are both of order $\bigo(n^{1/3} \log n)$. Therefore
\begin{equation*}
    \sum_{i=1}^m (k_i^- - k_i^+) = \Theta(n^{2/3}).
\end{equation*}
In contrast, $\Delta_{n-m}$ is centered, so $\lltEvent_n$ is looking at the event that $\Delta_{n-m}$ takes a value at distance $\Theta(n^{2/3})$ away from its mean. As stated in \cref{rem:llt-limitations}, the discrete local limit theorem provides no useful information in this regime. 

The second problem is that even in absence of the indicator function, the expectation being evaluated in \cref{eq:measure-change-numerator} is not dictated by the typical fluctuations of the random variables $\Xi^-_{n-m}$. In other words, it is not the case that
\begin{equation}
    \E\left[ 
        \prod_{i=1}^m \frac{(n - i + 1) \mu}{\sum_{j=i}^m k_j^- + \Xi_{n-m}^-}
    \right]
    \not \approx
    \prod_{i=1}^m \frac{(n - i + 1) \mu}{\sum_{j=i}^m k_j^- + \E[\Xi_{n-m}^-]}
\end{equation}

It turns out that both of these issues can be addressed by introducing a sequence of exponentially tilted measures.  The first effect of the exponentially tilted measures will be to shift the mean of $\Delta_{n-m}$ in such a way that, after the tilting, the event $\lltEvent_n$ concerns only a typical deviation of $\Delta_{n-m}$ which can be addressed by a local limit theorem. The second effect is that the expectation being evaluated in \cref{eq:measure-change-numerator} will be dictated by the typical fluctations of $\Xi^-_{n-m}$ under the tilted measure. 

The next result defines this tilt and then gives asymptotic expansions for cumulant generating function of $D^-$, the mean of $D^-$ and the mean of $D^+$ under this tilting. 
\begin{lemma}
    \label{lem:asym-expansions}
    Define an measure $\P_{\theta}$, for $\theta \geq 0$, by its Radon--Nikodym derivative
    \begin{equation*}
        \diff{\P_{\theta}}{\P} = \exp\left( - \theta D^- - \alpha(\theta) \right)
        \quad \text{where} \quad
        \alpha(\theta) = \log \E \left[ e^{-\theta D^-} \right].
    \end{equation*}
    Then as $\theta \downarrow 0$ we have
    \begin{align*}
        \alpha(\theta) &= -\mu \theta + \tfrac{1}{2}\var(D^-) \theta^2 - \tfrac{1}{6} \E \left[ (D^- - \mu)^3 \right] \theta^3 + \littleo(\theta^3), \\
        \E_{\theta}[D^-] &= \mu - \var(D^-) \theta + \bigo(\theta^2), \\
        \text{and} \quad \E_{\theta}[D^+] &= \mu - \cov(D^-, D^+) \theta + \bigo(\theta^2).
    \end{align*}
\end{lemma}
\begin{proof}
    Since $\E\left[\abs{D^-}^3\right] < \infty$ and $D^-$ is non-negative, by the dominated convergence theorem
    \begin{equation}
        \E \left[ (D^-)^3 \exp(-\theta D^-) \right] = \E \left[ (D^-)^3 \right] + \littleo(1)
        \label{eq:mgf-start}
    \end{equation}
    as $\theta \downarrow 0$. Integrating \cref{eq:mgf-start} with respect to $\theta$ and applying Fubini's theorem to exchange the order of the expectation and integral gives
    \begin{equation*}
        \E \left[ \int_0^{\theta} (D^-)^3 e^{-\theta' D^-} \dif \theta' \right]
        = \E \left[ \int_0^{\theta} \left\{ (D^-)^3 + \littleo(1) \right\} \dif \theta' \right]
        = \E \left[ (D^-)^3 \right] \theta + \littleo(\theta).
    \end{equation*}
    Evaluating the integral with respect to $\theta'$ on the left hand side and rearranging gives that
    \begin{equation*}
        \E \left[ (D^-)^2 e^{-\theta D^-} \right]
        = \E \left[ (D^-)^2 \right] - \E \left[ (D^-)^3 \right] \theta + \littleo(\theta).
    \end{equation*}
    Repeating this method yields
    \begin{align}
        \E \left[ D^- e^{-\theta D^-} \right]
        &= \mu - \E \left[ (D^-)^2 \right] \theta + \tfrac{1}{2} \E \left[ (D^-)^3 \right] \theta^2 + \littleo(\theta^2), \label{eq:asymin} \\
        \text{and} \quad \E \left[ e^{-\theta D^-} \right]
        &= 1 - \mu \theta + \tfrac{1}{2} \E \left[ (D^-)^2 \right] \theta^2 - \tfrac{1}{6} \E \left[ (D^-)^3 \right] \theta^3 + \littleo(\theta^3).
        \label{eq:asym00}
    \end{align}
    Similarly integrating the equation
    \begin{equation*}
        \E \left[ (D^-)^2 D^+ \exp(-\theta D^-) \right] = \E \left[ (D^-)^2 D^+ \right] + \littleo(1)
    \end{equation*}
    twice gives
    \begin{equation}
        \E \left[ D^+ e^{-\theta D^-} \right]
        = \mu \theta - \E \left[ D^- D^+ \right] \theta + \tfrac{1}{2} \E \left[ (D^-)^2 D^+ \right] \theta^2 + \littleo(\theta^2).
        \label{eq:asymout}
    \end{equation}
    \cref{eq:asym00} gives the small-$\theta$ expansion of the normalising constant of the measure change. Combining this with \cref{eq:asymin} and \cref{eq:asymout} yields the expansions for $\E_{\theta}[D^-]$ and $\E_{\theta}[D^+]$ respectively. Taking the logarithm of \cref{eq:asym00} gives the expansion of the cumulant generating function $\alpha(\theta)$.
\end{proof}

To achieve the recentering of $\Delta_{n-m}$ we desire, let us define a sequence of tilted measures $\P_n$ defined by their Radon--Nikodym derivative
\begin{equation}
    \diff{\P_n}{\P} = \exp \left( - \theta_n \Xi^-_{n-m} - (n - m) \alpha(\theta_n) \right),
\end{equation}
where $\theta_n = \frac{m}{\mu n}$. This factorises and so $\vD_1, \ldots, \vD_n$ remain i.i.d.\ under this tilting, each having the law of $\vD$ under $\P_{\theta_n}$. Applying \cref{lem:asym-expansions}, we can compute that
\begin{align*}
    \E_n[\Delta_{n-m}] 
    &= m(\lambda_+ - \lambda_-) + \bigo(n^{1/3}).
\end{align*}
Hence,
\begin{align*}
    \textstyle \sum_{i=1}^m (k_i^+ - k_i^-) - \E_n[\Delta_{n-m}] 
    &= s^-(m) - s^+(m) + \Big[ m(\lambda_+ - \lambda_-) - \E_n[\Delta_{n-m}] \Big] \\
    &= \bigo(n^{1/3} \log n),
\end{align*}
which is within the $\bigo(n^{1/2})$ range from the mean required for a typical deviation. This justifies our choice of $\theta_n = \frac{m}{\mu n}$. 

\subsubsection{Expansion of the numerator}

Remarkably the same tilting to apply the local limit theorem also correctly recenters $\Xi_{n-m}^-$ such that the expectation in \cref{eq:measure-change-numerator} is dominated by the typical behaviour of $\Xi_{n-m}^-$ under $\P_n$. Using \cref{lem:asym-expansions}, we have that
\begin{equation*}
    \E_n[\Xi_{n-m}^-] =  \mu n - \lambda_- m + \bigo(n^{1/3})
\end{equation*}
under the tilting. Thus we will expand the numerator under the event
\begin{equation*}
    \omegaEvent_n = \left\{ 
        \abs{\Xi^-_{n-m} - \mu n + \lambda_- m} \leq n^{1/2} \log(n)
     \right\}.
\end{equation*}
This event is saying that $\Xi_{n-m}^-$ is at `typical fluctations' from its tilted mean. The next lemma then expands the numerator of $\phi^n_m$ on the event $\omegaEvent_n$.
\begin{lemma}
    \label{lem:measure-change-numerator}
    We have that
    \begin{align*}
        \E\left[ 
            \one_{\lltEvent_n \cap \omegaEvent_n}
            \prod_{i=1}^m \frac{(n - i + 1) \mu}{\sum_{j=1}^m k_j^- + \Xi_{n-m}^-}
        \right]
        & =
        \left\{ \exp\left(
            \frac{1}{\mu n} \sum_{i=0}^m (s^-(i) - s^-(m)) - \frac{\sigma_-}{6 \mu^2} \frac{m^3}{n^2}
        \right) + \littleo(1) \right\} \\
        & \hspace{18em} \times \P_n(\lltEvent_n \cap \omegaEvent_n)
    \end{align*}
    where the $\littleo(1)$ term is bounded independently of $\vk_1, \ldots, \vk_m$ satisfying the assumption in \vref{eq:s-condition}.
\end{lemma}
\begin{proof}
    Firstly,
    \begin{equation*}
        \prod_{i=1}^m \frac{(n-i+1)\mu}{\sum_{j=i}^m k_j^- + \Xi^-_{n-m}} = \exp(X_n - Y_n),
    \end{equation*}
    where
    \begin{equation*}
        X_n = \sum_{i=1}^m \log\left( 1 - \frac{i-1}{n} \right)
        \quad \text{and} \quad
        Y_n = \sum_{i=1}^m \log\left( \frac{\sum_{j=i}^m k_j^- + \Xi^-_{n-m}}{\mu n} \right).
    \end{equation*}
    Note that
    \begin{equation*}
        \sum_{j=i}^m k^-_j = s^-(m) - s^-(i - 1) + (m - i + 1) \lambda_-.
    \end{equation*}
    For convenience, define
    \begin{equation*}
        \Omega^-_n = \Xi^-_{n-m} - \mu n + \lambda_- m
    \end{equation*}
    such that $\omegaEvent = \{ \abs{\Omega^-_n} < n^{1/2} \log n \}$. Then we have
    \begin{align*}
        Y_n 
        &= \sum_{i=1}^m \log \left( \frac{s^-(m) - s^-(i-1) + (m - i + 1) \lambda_- + \Omega^-_n + \mu n - \lambda_- m}{\mu n} \right) \\
        &= \sum_{i=1}^m \log \left( 1 + A_{i, n} + B_{i, n} \right)
    \end{align*}
    where
    \begin{align*}
        A_{i, n} = \frac{1}{\mu n} \left\{ 
            \Omega_n^- -\left[ s^-(i-1) - s^-(m) \right] 
        \right\}, \quad
        B_{i, n} = - \frac{\lambda_-}{\mu n} (i-1).
    \end{align*}
    Then on the event $\omegaEvent_n$,
    \begin{equation*}
        \max_{i=1, \ldots, m} \abs{A_n^i} = \bigo(n^{-1/2} \log n)
        \quad \text{and} \quad
        \max_{i=1, \ldots, m} \abs{B_n^i} = \bigo(n^{-1/3}).
    \end{equation*}
    where the $\bigo$ bounds are uniform for $\vk_1, \ldots, \vk_m$ satifying \cref{eq:s-condition}. There are $m = \theta(n^{2/3})$ terms in the summation. Thus to keep all terms of order $\Omega(1)$, we keep terms of order $\Omega(n^{-2/3})$, uniformly in $i$, when expanding $\log(1 + A_{i, n} + B_{i, n})$. The only such terms are $A_{i, n}, B_{i, n}$ and $B_{i, n}^2$. Moreover,
    \begin{equation*}
        \sum_{i=1}^m B_n^i = - \frac{\lambda_-}{2 \mu} \frac{m^2}{n} + \littleo(1) \quad \text{and} \quad
        \sum_{i=1}^m (B_n^i)^2 = \frac{\lambda_-^2}{3 \mu^2} \frac{m^3}{n^2} + \littleo(1).
    \end{equation*}
    Therefore,
    \begin{align*}
        Y_n
        &= \sum_{i=1}^m (A_{i, n} + B_{i, n} - \tfrac{1}{2} B_{i, n}^2) + \littleo(1) \\
        &= - \frac{1}{\mu n} \sum_{i=0}^m \left( s^-(i) - s^-(m) \right)
        + \frac{m}{\mu n} \Omega_n^-
        - \frac{\lambda_-}{2\mu} \frac{m^2}{n} - \frac{\lambda_-^2}{6 \mu^2} \frac{m^3}{n^2} + \littleo(1),
    \end{align*}
    where we use that $\sum_{i=1}^m \left( s^-(i-1) - s^-(m) \right) = \sum_{i=0}^m \left( s^-(i) - s^-(m) \right)$.

    Similarly we can expand $X_n$ as
    \begin{equation*}
        X_n = - \frac{m}{2 n} - \frac{m^3}{3 n^2} + \littleo(1).
    \end{equation*}

    Thus,
    \begin{align*}
        &\one_{\lltEvent_n \cap \omegaEvent_n} \prod_{i=1}^m \frac{(n-i+1)\mu}{\sum_{j=i}^m k_j^- + \Xi^-_{n-m}} \\
        = &\exp \Bigg( \frac{1}{\mu n} \sum_{i=1}^m (s^-(i) - s^-(m)) \nonumber
         - \frac{m}{\mu n} \Omega^-_n + \frac{(\lambda_- - \mu)}{2 \mu} \frac{m^2}{n} + \frac{(\lambda_-^2 - \mu^2)}{6 \mu^2} \frac{m^3}{n^2} + \littleo(1) \Bigg) \one_{\lltEvent_n \cap \omegaEvent_n}.
    \end{align*}
    In addition, using \cref{lem:asym-expansions}, the measure change can be expanded as
    \begin{equation*}
        \diff{\P_n}{\P} = \exp \left( 
            - \frac{m}{\mu n} \Omega_n^- + \frac{(\lambda_- - \mu)}{2\mu} \frac{m^2}{n}
            + \frac{(\lambda_-^2 - \mu^2)}{6 \mu^2} \frac{m^3}{n^2} + \frac{\sigma_-}{6 \mu^2} \frac{m^3}{n^2} + \littleo(1)
        \right).
    \end{equation*}
    Hence,
    \begin{align*}
        &\E\left[ 
            \one_{\lltEvent_n \cap \omegaEvent_n}
            \prod_{i=1}^m \frac{(n-i+1)\mu}{\sum_{j=i}^m k^-_j + \Xi^-_{n-m}}
        \right] \\
        =& \E_n\left[ 
            \frac{\dif \P}{\dif \P_n}
            \one_{\lltEvent_n \cap \omegaEvent_n}
            \prod_{i=1}^m \frac{(n-i+1)\mu}{\sum_{j=i}^m k^-_j + \Xi^-_{n-m}}
        \right] \\
        =& \E_n \left[ 
            \one_{\lltEvent_n \cap \omegaEvent_n}
            \exp \left( 
                \frac{1}{\mu n} \sum_{i=0}^m (s^-(i) - s^-(m)) - \frac{\sigma_-}{6 \mu^2} \frac{m^3}{n^2} + \littleo(1)
            \right)
        \right] \\
        =& \left\{  
            \exp \left( 
                \frac{1}{\mu n} \sum_{i=0}^m (s^-(i) - s^-(m)) - \frac{\sigma_-}{6 \mu^2} \frac{m^3}{n^2}
            \right) + \littleo(1)
         \right\} \P_n(\lltEvent_n \cap \omegaEvent_n)
    \end{align*}
    as required.
\end{proof}

\subsubsection{Multivariate Local Limit Theorem}

To complete the proof of \cref{prop:measure-change-approx} we need to understand the asymptotic behaviour of $\P_n(\lltEvent_n \cap \omegaEvent_n)$. Recall an effect of the tilting was to center $\Delta_{n-m}$ in such a way that the probability of the event
\begin{equation*}
    \lltEvent_n = \left\{ 
        \Delta_{n-m} = \textstyle \sum_{i=1}^m (k_i^+ - k_i^-)
    \right\}
\end{equation*}
can be addressed by the local limit theorem. However, due to the tilting, $\P_n$ changes with $n$. In effect, $\Delta_n$ under $\P_n$ has the same distribution as $\sum_{i=1}^{n-m} X_{n, i}$ where $(X_{n, i})_{i=1}^n$ has the same joint distribution as $(D_i^- - D_i^+)_{i=1}^n$ under $\P_n$. Then $X_{n, 1}, \ldots, X_{n, n}$ are i.i.d.\ but the distribution of $X_{n, 1}$ can change with $n$. A collection of random variables $(X_{n, 1}, \ldots, X_{n, n})_{n = 1}^{\infty}$ satisfying this property is a \emph{row-wise i.i.d.\ triangular array}. Thus we require a generalisation of the discrete local limit theorem which can deal with such arrays. In addition, to deal with the event
\begin{equation*}
    \omegaEvent_n = \left\{ 
        \abs*{\Xi_{n-m}^- - \mu n + \lambda_- m} \leq n^{1/2} \log n
    \right\},
\end{equation*}
we will prove a multivariate local limit theorem applicable to $(\Delta_{n-m}, \Xi^-_{n-m})$ under $\P_n$ and then sum over the possible values of $\Xi^-_n$.

Before we state the result we use, we first define some terminology regarding lattices in $\R^d$. A set of points in $\R^d$ is a \emph{lattice} if there exists a basis $\va_1, \ldots, \va_d$ of $\R^d$ such that
\begin{equation*}
    \lattice = \left\{ 
        \textstyle \sum_{i=1}^d n_i \va_i : n_i \in \Z \text{ for $i = 1, \ldots, d$}
    \right\}.
\end{equation*}
We say $\lattice$ is generated by $\va_1, \ldots, \va_d$. We can summarise the basis by a $n \times n$ matrix $A$ whose columns are $\va_1, \ldots, \va_n$. In other words $A_{ij} = \va_j^{(i)}$. The choice of basis generating a lattice is not unique, and the following lemma adapted from \cite[Corollary 4.3a]{schrijverTheoryLinearInteger1998} characterises when two basis generate the same lattice.
\begin{lemma}
    Let $A$ and $B$ be $n \times n$ matrices of full rank. Then the columns of $A$ and $B$ generate the same matrix if and only if there exists a matrix $U$ such that $U$ has integer entries, $\det(U) = \pm 1$ and $A = UB$.
\end{lemma}
Therefore we can define $\det(\lattice)$ to be $\abs{\det(A)}$ for any matrix $A$ whose columns generate $\lattice$, and this definition is independent of the choice of $A$.

For integer lattices, we can obtain a canonical choice of the basis generating the lattice. We say a $d \times d$ matrix $A$ is in \emph{Hermite normal form} if $A$ is lower triangular with entries 
\begin{equation*}
    A = \begin{pmatrix}
        a_{1, 1} &        & 0 \\
        \vdots   & \ddots &   \\
        a_{d, 1} & \cdots & a_{d, d}
    \end{pmatrix}
\end{equation*}
satisfying
\begin{enumerate}
    \item $a_{i, j}$ is a non-negative integer for all $i = 1, \ldots, d$ and $j \geq i$,
    \item $a_{i, i} > 0$ for all $i = 1, \ldots, d$, and
    \item $a_{j, i} < a_{i, i}$ for all $j > i$.
\end{enumerate}
Then the following lemma, adapted from \cite[Corollary 1]{casselsIntroductionGeometryNumbers1997}, gives existence of a canonical choice of basis generating an integer lattice.
\begin{lemma}
    \label{lem:hnf-basis}
    Suppose $\lattice \subset \Z^d$ is a lattice. Then there exists a $d \times d$ matrix $A$ in Hermite normal form such that the columns of $A$ form a basis which generates $\lattice$.
\end{lemma}

An $\R^d$-valued random variable $\vX$ is \emph{non-degenerate} if it is not supported on an affine hyperplane of $\R^d$. $\vX$ is \emph{lattice} if it is non-degenerate and supported on a translation of a lattice. To avoid dealing with translations, it is convenient to work with the \emph{symmetrisation} of $\vX$. This is the random variable $\vX^* = \vX_1 - \vX_2$ where $\vX_1$ and $\vX_2$ are independent copies of $\vX$. For each lattice $\lattice$, $\vX$ is supported on a translation of $\lattice$ if and only if $\vX^*$ is supported on $\lattice$ without translation.

If $\vX$ is lattice, the \emph{main lattice} $\lattice(\vX)$ of $\vX$ is the intersection of all lattices containing the support of $\vX^*$. This is in itself a lattice, and is explicitly given by
\begin{equation*}
    \lattice(\vX) = \bigcup_{k=1}^{\infty} \left\{ 
        \textstyle \sum_{i=1}^k n_i \vx^*_i : \text{$n_i \in \Z$ and $\vx^*_i \in \supp(\vX^*)$ for $i = 1, \ldots, k$}
    \right\}.
\end{equation*}
It will turn out that if $\vX$ is an $\R^d$-valued lattice random variable with main lattice $\lattice$, then $\det(\lattice)$ can be seen as a generalisation of the span of an $\R$-valued random variable.

To deal with the triangular array, we recall the exponential tilt is given by
\begin{equation*}
    \frac{\dif \P_n}{\dif \P} = \exp(
        - \theta_n \Xi^-_{n-m} - (n - m) \alpha(\theta_n)
    )
\end{equation*}
where $\theta_n = \frac{m}{\mu n}$. Since $\theta_n \to 0$, the distribution of $\vD_i$ under $\P_n$ is converging to that of $\vD_i$ under $\P$ as $n \to \infty$. This allows us to ignore the tilting in the limit.

\begin{restatable}{theorem}{llt}
    \label{thm:multi-triangular-llt}
    For each $n \geq 1$ let $\vX_n$ be an $\R^d$ valued random variable and
    \begin{equation*}
        \vX_{n, 1}, \vX_{n, 2}, \ldots, \vX_{n, n}  
    \end{equation*}
    be i.i.d.\ copies of $\vX_n$. Assume that the following holds:
    \begin{enumerate}
        \item There exists a random variable $\vX$ such that $\vX_n \todist \vX$ as $n \to \infty$.
        \item $(\norm{\vX_n}^2)_{n \geq 1}$ is a uniformly integrable sequence of random variables. Explicitly
            \begin{equation}
                \label{eq:ui-condition}
                \lim_{L \to \infty} \sup_n \E\left[
                    \norm{\vX_n}^2
                    \one \left\{ \norm{\vX_n}^2 > L \right\}
                \right] = 0.
            \end{equation}
        \item For all $n$, $\vX_n$ and $\vX$ are lattice with common main lattice $\lattice$.
    \end{enumerate}
    Then $\vX$ has finite second moment. Further, for each $n$ let $\vc_n$ be an arbitrary element in the support of $\sum_{i=1}^n \vX_{n, i}$. Then uniformly for $\vy \in \vc_n + \lattice$,
    \begin{equation*}
        \P\Big(
            \textstyle \sum_{i=1}^n \vX_{n, i} = \vy
        \Big)
        = n^{-d/2} \det(\lattice) f\left( \vx_n(\vy) \right) + \littleo \left( n^{-d/2} \right)
        \quad \text{where} \quad
        \vx_n(\vy) = \frac{\vy - n \E[\vX_n]}{\sqrt{n}}
    \end{equation*}
    and $f$ is the density of a $N(0, \cov(\vX))$ distribution. This means that
    \begin{equation*}
        \lim_{n \to \infty} \sup_{\vy \in \vc_n + \lattice} \abs*{
            n^{d/2} \P\Big( \textstyle \sum_{i=1}^n \vX_{n, i} = \vy \Big)
            - \det(\lattice) f(\vx_n(\vy))
        } = 0.
    \end{equation*}
\end{restatable}

We defer the proof of this to \cref{sec:llt} in the appendix, and instead make a few remarks. Firstly $\vX$ is assumed to be lattice and thus non-degenerate. Hence $\cov(\vX)$ is invertible, ensuring $N(0, \cov(\vX))$ has a valid density $f$, which is explicitly given by
\begin{equation*}
    f(\vx) = \frac{1}{\sqrt{(2 \pi)^d \det(\cov(\vX))}}
    \exp \Bigg(
        -\frac{1}{2} \vx \cdot \cov(\vX)^{-1} \vx
        \Bigg).
\end{equation*}

Secondly, since the $\vX_1, \vX_2, \ldots$ do not necessarily live in the same probability space we should not technically refer to the sequence $(\norm{\vX_n}^2)_{n \geq 1}$ as uniformly integrable. However the condition in \cref{eq:ui-condition} is still well defined.

We apply \Cref{thm:multi-triangular-llt} to $(\Xi_{n-m}^-, \Delta_{n-m})$. Suppose $(D^- - D^+, D^-)$ is non-degenerate and let $\lattice$ be its main lattice. By \cref{lem:hnf-basis}, $\lattice$ is generated by the columns of a matrix $A$ in Hermite normal formal. Since $D^- - D^+$ has span 1, it must be the case that $A_{1, 1} = 1$. Thus there is some positive integer $q$ such that
\begin{equation*}
    A = \begin{pmatrix}
        1 & 0 \\ 
        0 & q
    \end{pmatrix}.
\end{equation*}
Finally let $\Sigma$ be the covariance matrix of $(D^- - D^+, D^-)$. With this notation, the following lemma holds:
\begin{lemma}
    \label{lem:bivar-llt}
    Suppose $(D^- - D^+, D^-)$ is non-degenerate. For each $n$, let $\vc_n$ be in the support of $(\Delta_{n-m}, \Xi^-_{n-m})$. Then uniformly for $(x, y) \in \vc_n + \Lambda$,
    \begin{align*}
        &\P_n\left(
            \Delta_{n-m} = \E \big[ \Delta_{n-m} \big] + x, \ 
            \Xi^-_{n-m} = \E \big[ \Xi^-_{n-m} \big] + y
        \right)  \nonumber \\
        &\hspace{7em} = \frac{q}{2\pi \det(\Sigma)^{1/2} \: n} \exp\left( 
            \frac{-1}{2n}
            \begin{pmatrix}
                x & y
            \end{pmatrix}
            \Sigma
            \begin{pmatrix}
                x \\ y
            \end{pmatrix}
        \right)
        + \littleo( n^{-1} )
    \end{align*}
    as $n \to \infty$.
\end{lemma}
\begin{proof}
    Let $\vX = (D^- - D^+, D^-)$. For each $n$, let $\vX_n$ be distributed as $(D^- - D^+, D^-)$ under $\P_{\theta_n}$. Then
    \begin{equation*}
        \vectwo{D_1^- - D_1^+}{D_1^-},
        \ldots,
        \vectwo{D_n^- - D_n^+}{D_n^-}
    \end{equation*}
    under $\P_n$ can be seen as $n$ i.i.d.\ copies of $\vX_n$. Since $\theta_n \to 0$, we have that $\vX_n \todist \vX$ as $n \to \infty$. 

    For any $L > 0$,
    \begin{align*}
        \sup_n \E \left[ 
            \norm{\vX_n}^2 \one\{\norm{\vX_n}^2 > L\}
        \right]
        &=
        \sup_n \E \left[ 
            e^{- \theta_n D^- - \alpha(\theta_n)}
            \norm{\vX}^2 \one\{\norm{\vX}^2 > L\}
        \right] \\
        &\leq
        \left(\sup_n e^{- \alpha(\theta_n)}\right) 
        \E \left[ 
            \norm{\vX}^2 \one\{\norm{\vX}^2 > L\}
        \right]
    \end{align*}
    since $\theta_n$ and $D^-_n$ are non-negative. Since $\theta_n$ is convergent,
    \begin{equation*}
        \sup_n e^{-\alpha(\theta_n)} < \infty.
    \end{equation*}
    Moreover $\E \left[ \norm{\vX}^2 \one\{\norm{\vX}^2 > L\} \right] \to 0$ as $L \to \infty$ as $\vX$ has finite second moment. Thus $(\norm{\vX_n}^2)_{n \geq 1}$ satisfies the uniform integrability condition in \cref{eq:ui-condition}.

    Finally the exponential tilt does not change the support of the random variables. Thus $\vX$ and $\vX_n$ share a common main lattice $\lattice$. In addition, $\det(\lattice) = q$.

    Hence the result follows by \cref{thm:multi-triangular-llt}. There is a small change in that we are considering a sum of $n - m$ random variables rather than $n$. However since $m = \littleo(n)$, the same asymptotic result holds.
\end{proof}

Now we show $\P(\lltEvent_n, \omegaEvent_n)$ has the same asymptotic behaviour as $\P(\Delta_n = 0)$. We only prove a lower bound, but this is sufficient for proving \cref{prop:measure-change-approx}.
\begin{lemma}
    \label{lem:mod-dev-local}
    Under the assumptions of \cref{prop:measure-change-approx},
    \begin{equation*}
        \mathbb P_n \left(
            \Delta_{n-m} = \sum_{i=1}^m (k_i^+ - k_i^-), \ 
            \abs{\Xi^-_{n-m} - \E_n[\Xi^-_{n-m}]} \leq n^{\frac{1}{2}} \log n
        \right)
        \geq \frac{1}{\sqrt{2 \pi \sigma^2 n}} (1 + \littleo(1)).
    \end{equation*}
\end{lemma}
\begin{proof}
    For convenience let
    \begin{equation*}
        P_n = \mathbb P_n \left(
            \Delta_{n-m} = \sum_{i=1}^m (k_i^+ - k_i^-), \ 
            \abs{\Xi^-_{n-m} - \E_n[\Xi^-_{n-m}]} \leq n^{\frac{1}{2}} \log n
        \right).
    \end{equation*}
    Firstly, suppose $(D^- - D^+, D^-)$ is degenerate. Then since we assume that $D^- - D^+$ is non-constant, it must be the case that either $D^-$ or $D^+$ is constant. Either way, it becomes the case that
    \begin{equation*}
        \left\{
            \Delta_{n-m} = \sum_{i=1}^m (k_i^+ - k_i^-), \ 
            \abs{\Xi^-_{n-m} - \E_n[\Xi^-_{n-m}]} \leq n^{\frac{1}{2}} \log n
        \right\}
        =
        \left\{
            \Delta_{n-m} = \sum_{i=1}^m (k_i^+ - k_i^-)
        \right\}.
    \end{equation*}
    Then applying \cref{thm:multi-triangular-llt}, as we did in the proof of \cref{lem:bivar-llt}, shows that
    \begin{equation*}
        \P\left( 
            \Delta_{n-m} = \textstyle \sum_{i=1}^m (k_i^+ - k_i^-)
        \right) = \frac{1}{\sqrt{2 \pi \sigma^2 n}} (1 + \littleo(1)).
    \end{equation*}
    Otherwise assume that $(D^- - D^+, D^-)$ is non-degenerate. Define
    \begin{equation*}
        a_n = \sum_{i=1}^m (k_i^+ - k_i^-) - \E_n[\Delta_{n-m}].
    \end{equation*}
    Also let
    \begin{equation*}
        L_n = \Big\{
            y : \bigg( \textstyle\sum_{i=1}^m (k_i^+ - k_i^-), y \bigg) \in \vc_n + \lattice
            \Big\}.
    \end{equation*}
    $L_n$ has a simpler representation. Fix any $y_0 \in L_n$. Then if $\lattice$ is generated by the columns of
    \begin{equation*}
        \begin{pmatrix}
            1 & 0 \\
            0 & q
        \end{pmatrix}
    \end{equation*}
    we must have $L_n = y_0 + q\Z$. Fix an arbitrary $M > 0$. Then
    \begin{align*}
        P_n &= \sum_{\substack{y \in L_n \\ \abs{y} \leq n^{1/2} \log n}} \mathbb P_n \left( \Delta_{n-m} = \E_n[\Delta_{n-m}] + a_n, \ \Xi^-_{n-m} = \E_n[\Xi^-_{n-m}] + y \right) \\
        &\geq \sum_{\substack{y \in L_n \\ \abs{y} \leq M n^{1/2}}} \mathbb P_n \left( \Delta_{n-m} = \E_n[\Delta_{n-m}] + a_n, \ \Xi^-_{n-m} = \E_n[\Xi^-_{n-m}] + y \right)
    \end{align*}
    for all $n$ sufficiently large. By \cref{lem:bivar-llt}, using that the error is uniform, we have that
    \begin{equation*}
        P_n \geq \sum_{\substack{y \in L_n \\ \abs{y} \leq M n^{1/2}}} 
         \frac{q}{2\pi \det(\Sigma)^{1/2} \: n} \exp\left( 
            \frac{-1}{2n} \vectwo{a_n}{y} \cdot \Sigma^{-1} \vectwo{a_n}{y}
         \right)
         + \littleo( n^{-1/2} )
    \end{equation*}

    We wish to factorise the summand. To this end, we make a change of variables. There exists $c \in \R$ such that
    \begin{equation*}
        \cov(D^- - c(D^- - D^+), D^- - D^+) = 0.
    \end{equation*}
    Let $\tau^2$ be the variance of $D^- - c(D^- - D^+)$. Then
    \begin{align*}
         &\frac{q}{2\pi \det(\Sigma)^{1/2} \: n} \exp\left( 
            \frac{1}{2n} \vectwo{a_n}{y} \cdot \Sigma^{-1} \vectwo{a_n}{y}
         \right) \nonumber \\
         & \hspace{5em} =
         \frac{1}{\sqrt{2 \pi \sigma^2 n}} \exp \left( - \frac{1}{2 \sigma^2} \frac{a_n^2}{n} \right)
         \frac{q}{\sqrt{2 \pi \tau^2 n}} \exp \left( - \frac{1}{2 \tau^2} \frac{(y - ca_n)^2}{n} \right).
    \end{align*}
    We now examine the asymptotic behaviour of $a_n$. By \cref{lem:asym-expansions},
    \begin{align*}
        \E_n[\Delta_{n-m}]
        &= (n - m) \E_{\theta_n}[D^- - D^+] \\ 
        &= -(\lambda_- - \lambda_+)m + \bigo(n^{1/3}).
    \end{align*}
    Therefore 
    \begin{equation*}
        a_n = s_+(m) - s_-(m) + \bigo(n^{1/3}) = \bigo(n^{1/3} \log n),
    \end{equation*}
    by the assumption in \cref{eq:s-condition}. Thus
    \begin{equation*}
        P_n \geq
        \frac{1}{\sqrt{2 \pi \sigma^2 n}}(1 + \littleo(1))
        \sum_{\substack{y \in L_n \\ \abs{y} \leq M n^{1/2}}}
        \frac{q}{\sqrt{2 \pi \tau^2 n}} \exp \left( - \frac{1}{2 \tau^2} \frac{(y - ca_n)^2}{n} \right)
        + \littleo(n^{-1/2})
    \end{equation*}
    Note that
    \begin{equation*}
        \sum_{\substack{y \in L_n \\ \abs{y} \leq M n^{1/2}}}
        \frac{q}{\sqrt{2 \pi \tau^2 n}} \exp \left( - \frac{1}{2 \tau^2} \frac{(y - ca_n)^2}{n} \right)
        = \sum_{\substack{y \in L_n \\ \abs{y} \leq M n^{1/2}}}
        \frac{q}{\sqrt{n}}\ g\left( \frac{y - ca_n}{\sqrt{n}} \right)
    \end{equation*}
    where
    \begin{equation*}
        g(z) = \frac{1}{\sqrt{2 \pi \tau^2}} \exp\left( \frac{-z^2}{2 \tau^2} \right).
    \end{equation*}
    Since $a_n = \bigo(n^{1/3 + \epsilon})$, for $n$ sufficiently large
    \begin{align}
        \sum_{\substack{y \in L_n \\ \abs{y} \leq M n^{1/2}}}
        \frac{q}{\sqrt{2 \pi \tau^2 n}} \exp \left( - \frac{1}{2 \tau^2} \frac{(y - ca_n)^2}{n} \right)
        &\geq \sum_{\substack{z \in L_n - ca_n\\ \abs{z} \leq \frac{1}{2} M n^{1/2}}} 
        \frac{q}{\sqrt{n}}\ g \left( \frac{z}{\sqrt{n}} \right) \\
        &= \sum_{\substack{z \in \tilde{L}_n \\ \abs{z} \leq \frac{1}{2} M }}
        \frac{q}{\sqrt{n}}\ g(z) \label{eq:riemann-sum}
    \end{align}
    where
    \begin{equation*}
        \tilde{L}_n = \frac{L_n - ca_n}{\sqrt{n}}.
    \end{equation*}
    Then $\tilde{L}_n \cap [-\frac{1}{2}M, \frac{1}{2}M]$ is a partition of $[-\frac{1}{2}M, \frac{1}{2}M]$ where adjacent points are distance $q/\sqrt{n}$ apart from each other. Thus \cref{eq:riemann-sum} is a Riemann sum approximation of an integral. Hence
    \begin{equation*}
        \sum_{\substack{y \in L_n \\ \abs{y} \leq M n^{1/2}}}
        \frac{q}{\sqrt{2 \pi \tau^2 n}} \exp \left( - \frac{1}{2 \tau^2} \frac{(y - ca_n)^2}{n} \right)
        \geq (1 + \littleo(1)) \int_{-\frac{1}{2} M}^{\frac{1}{2}M} g(z) \dif z.
    \end{equation*}
    Thus
    \begin{equation*}
        P_n \geq \frac{1}{\sqrt{2 \pi \sigma^2 n}} (1 + \littleo(1)) \int_{-\frac{1}{2} M}^{\frac{1}{2}M} g(z) \dif z.
    \end{equation*}
    This holds for all $M > 0$, and $\int_{-\infty}^{\infty} g(z) \dif z = 1$. Therefore,
    \begin{equation*}
        P_n \geq \frac{1}{\sqrt{2 \pi \sigma^2 n}} \left( 1 + \littleo(1) \right),
    \end{equation*}
    as required.
\end{proof}

\subsection{Proof of lower bound}

Now we are ready to prove \cref{prop:measure-change-approx}.

\begin{proof}[Proof of \cref{prop:measure-change-approx}]
    By \cref{lem:exact-measure-change} and \cref{lem:measure-change-numerator} we have that
    \begin{align*}
        \phi(\vk_1, \ldots, \vk_m)
        &\geq \left\{ \exp\left(
            \frac{1}{\mu n} \sum_{i=0}^m (s^-(i) - s^-(m)) - \frac{\sigma_-}{6 \mu^2} \frac{m^3}{n^2}
        \right) + \littleo(1) \right\} \frac{\P_n(\lltEvent_n \cap \omegaEvent_n)}{\P(\Delta_n = 0)}
    \end{align*}
    where the $\littleo(1)$ term is independent of $\vk_1, \ldots, \vk_m$ satisfying our assumptions. Then by \cref{lem:bivar-llt} and \cref{cor:measure-change-denominator-control} we have that
    \begin{equation*}
        \frac{\P_n(\lltEvent_n \cap \omegaEvent_n)}{\P(\Delta_n = 0)} \geq 1 + \littleo(1)
    \end{equation*}
    where the $\littleo(1)$ term is independent of $\vk_1, \ldots, \vk_m$ satisfying our assumptions. Thus
    \begin{equation*}
        \phi(\vk_1, \ldots, \vk_m) \geq
        \exp\left(
            \frac{1}{\mu n} \sum_{i=0}^m (s^-(i) - s^-(m)) - \frac{\sigma_-}{6 \mu^2} \frac{m^3}{n^2}
        \right) + \littleo(1)
    \end{equation*}
    as required.
\end{proof}

\subsection{Convergence of the measure change}

We are now ready to prove the main result of this section.

\begin{proof}[Proof of \cref{prop:measure-change-no-crit}]
    The existence of the measure change is covered by \cref{lem:exact-measure-change}. Define
    \begin{equation*}
        \Gamma(n, m) = \exp \left( 
            \frac{1}{\mu n} \sum_{i=0}^m \left( 
                V^-(i) - V^-(m)
             \right)
             - \frac{\sigma_-}{6 \mu^2} \frac{m^3}{n^2}
         \right).
    \end{equation*}
    Then by Donsker's invariance principle,
    \begin{equation*}
        \left( 
            n^{-1/3} V^-(\floor{tn^{2/3}}),
            n^{-1/3} V^+(\floor{tn^{2/3}})
         \right)_{t \geq 0} \todist
         \left( 
             \sigma_- W^-_t,
             \sigma_+ W^+_t
          \right)_{t \geq 0}
    \end{equation*}
    in $\D\left([0, \infty), \R^2\right)$, where $(W^-_t, W^+_t)_{t \geq 0}$ are a pair of correlated standard Brownian motions with correlation $\corr(Z^-_1, Z^+_1)$. We can write
    \begin{align*}
        \frac{1}{n} \sum_{i=0}^{\floor{T n^{2/3}}} V^-(i)
        &= n^{-2/3} \int_0^{\floor{T n^{2/3}} + 1} n^{-1/3} V^-(\floor{u}) \dif u \\
        &= \int_0^{n^{-2/3}(\floor{T n^{2/3} + 1})} n^{-1/3} V^-(\floor{s n^{2/3}}) \dif s.
    \end{align*}
    Thus, by the continuous mapping theorem,
    \begin{align*}
        \frac{1}{n} \sum_{i=0}^{\floor{T n^{2/3}}} \left(V^-(i) - V^-(m) \right)
        &\todist \int_0^T \left(W^-_s - W^-_T \right) \dif s
        = - \int_0^T s \dif W^-_s.
    \end{align*}
    Hence,
    \begin{align*}
        &\left( 
            \Gamma(n, \floor{T n^{2/3}}),
            \left(
                n^{-1/3} V^-(\floor{t n^{2/3}}),
                n^{-1/3} V^+(\floor{t n^{2/3}})
            \right)_{t \in [0, T]}
         \right) \\
         &\hspace{22em} \todist
         \left( 
             \Phi(T), (\sigma_- W^-_t, \sigma_+ W^+_t)_{t \in [0, T]}
          \right)
    \end{align*}
    in $\R \times \D([0, T], \R)$, as $n \to \infty$. Recall the event
    \begin{equation*}
        \walkDeviationEvent_m = \left\{ 
            \max_{i=1, \ldots, m} \abs{V^-(i)} \leq m^{1/2} \log m
            \quad \text{and} \quad
            \max_{i=1, \ldots, m} \abs{V^+(i)} \leq m^{1/2} \log m
         \right\}
    \end{equation*}
    By \cref{prop:measure-change-approx}, it is the case that
    \begin{equation*}
        \Phi(n, m) \geq (\Gamma(n, m) + \littleo(1)) \one_{\walkDeviationEvent_m}.
    \end{equation*}
    The processes $(V^{\pm}(n))_{n \geq 0}$ are discrete martingales. Therefore, by Doob's maximal inequality,
    \begin{equation*}
        \P\left(\max_{i=1, \ldots, m} \abs{V^{\pm}(i)} > m^{1/2} \log(m) \right)
        \leq \frac{\E[(V^{\pm}(m))^2]}{m (\log m)^2} = \frac{\sigma_{\pm}^2}{(\log m)^2} \to 0
    \end{equation*}
    as $m \to \infty$. Thus $\P(\walkDeviationEvent_m) \to 1$ as $m \to \infty$. Hence, we still have that
    \begin{align*}
        &\left( 
            (\Gamma(n, \floor{T n^{2/3}}) + \littleo(1))\one_{E_{\floor{T n^{2/3}}}} ,
            \left(
                n^{-1/3} V^-(\floor{t n^{2/3}}),
                n^{-1/3} V^+(\floor{t n^{2/3}})
            \right)_{t \in [0, T]}
         \right) \\
         & \hspace{22em} \todist
         \left( 
             \Phi(T), (\sigma_- W^-_t, \sigma_+ W^+_t)_{t \in [0, T]}
          \right).
    \end{align*}
    We have $\E[\Phi(T)] = 1$ by a standard stochastic calculus calculation. Therefore, by \cref{lem:sandwiching-lemma}, we get the desired result that
    \begin{align*}
        &\left( 
            \Phi(n, \floor{T n^{2/3}}) ,
            \left(n^{-1/3} V^-(\floor{t n^{2/3}}), n^{-1/3} V^+(\floor{t n^{2/3}}) \right)_{t \in [0, T]}
         \right) \\
        & \hspace{22em} \todist
         \left( 
             \Phi(T), (\sigma_- W^-_t, \sigma_+ W^+_t)_{t \in [0, T]}
          \right),
    \end{align*}
    and that $(\Phi(n, \floor{T n^{2/3}}))_{n \geq 1}$ is a uniformly integrable sequence.
\end{proof}
\section{Convergence of the out-forest}\label{sec.convoutforest}
 Fix $T>0$. In this section we will show that the \L ukasiewicz path and height process corresponding to the out-forest converge under rescaling up to time $\lfloor T n^{2/3}\rfloor$. Note that the out-forest will contain at least $n$ vertices, so for $n$ large enough, $\lfloor T n^{2/3}\rfloor\leq n$ and the encoding processes are well-defined up to time $\lfloor T n^{2/3}\rfloor$. 
 
We will show that the convergence under rescaling of the \L ukasiewicz path and height process $(\hat{S}^{+}_n(k),\hat{H}_n(k),k\leq \lfloor Tn^{2/3}\rfloor)$ occurs jointly with convergence in distribution under rescaling of $(\hat{S}^-_n(k),\hat{P}_n(k), k\leq \lfloor Tn^{2/3}\rfloor)$, for $\hat{S}^-_n(k)$ the number of unpaired in-half-edges of vertices that have been discovered at time $k$, and $\hat{P}_n(k)$ the number of dummy leaves added in the first $k$ time-steps.  \\
We let $(B_t)_{t\geq 0}$ be a Brownian motion, and define
$$(\hat{B}_t,t\geq 0):=\left( B_t-\frac{\sigma_{-+}+\nu_-}{2\sigma_+ \mu}t^2, t\geq 0\right).$$ 
We define the reflected process $$(\hat{R}_t,t\geq 0)= \left(\hat{B}_t-\inf\left\{\hat{B}_s: s\leq t\right\},t\geq 0\right).$$

The main result of this section is as follows. 

\begin{proposition}\label{prop:convoutforest}
It holds that
\begin{align*}\left(n^{-1/3}\hat{S}^{+}_n\left(\lfloor n^{2/3}t\rfloor \right),n^{-1/3}\hat{H}_{n}\left(\lfloor n^{2/3}t\rfloor \right), t\leq T \right)
\todist\left(\sigma_+ \hat{B}_t, \frac{2}{\sigma_+} \hat{R}_t, t\leq T \right)\end{align*}
in $\D([0,T],\R)^2$, and 
\begin{align*}\left( n^{-2/3}\hat{S}_n^-\left(\lfloor n^{2/3}t\rfloor \right), n^{-1/3}\hat{P}_n\left(\lfloor n^{2/3}t\rfloor \right), t\leq T \right)\toprob\left(\nu_-t,  \frac{\nu_-}{2\mu} t^2, t\leq T \right)\end{align*}
in $\D([0,T],\R)^2$ as $n\to \infty$. 
\end{proposition}
We prove \cref{prop:convoutforest} by studying two other forests that are related to the out-forest via a change of measure.  \\
The proof is structured as follows.
\begin{enumerate}
    \item \label{item.measurechangeexists} Recall that $(\mathbf{\hat{D}}_{n,1},\dots,\mathbf{\hat{D}}_{n,n})$ are the degree pairs of the vertices in order of discovery. Also recall $\mathbf{Z}_1, \mathbf{Z}_2, \ldots$ in an i.i.d.\ sequence of $\N\times \N$-valued random variables, $\mathbf{Z}_i:=(Z_i^-,Z_i^+)$, such that 
    $$\P(Z_i^-=k^-, Z_i^+=k^+)=\frac{k^-\P(D^-=k^-,D^+=k^+)}{\mu}.$$
    In \cref{sec:measure-change}, we showed that the law of $(\mathbf{\hat{D}}_{n,1},\dots,\mathbf{\hat{D}}_{n,m})$ conditional on $\sum_{i=1}^n D_i^-=\sum_{i=1}^n D_i^+$ and $m \leq R_n$ is absolutely continuous with respect to that of $(\mathbf{Z}_1,\dots, \mathbf{Z}_m)$, and we showed the convergence under rescaling of the Radon-Nikodym derivative $\phi_m^n$ for $m=\lfloor T n^{2/3}\rfloor$. 
    \item Point \ref{item.measurechangeexists} motivates us to study a Bienaymé forest with offspring distributed as $Z_1^+$. The convergence of the \L ukasiewicz path of this forest under rescaling follows from Donsker's theorem.
    \item In Subsection \ref{subsec.purpleleavesGWforest}, we modify the Bienaymé forest in order to include dummy leaves. We add extra randomness, approximating the procedure described in \cref{prop:sampleoutforest}, in such a way that at some time-steps, a dummy leaf is added. We call the resulting forest \emph{the forest with dummy leaves}. We respect the order of the degrees in the Bienaymé forest, in the sense that for any $k$, the $k$th true vertex in the forest with dummy leaves has the same number of children as the $k$th vertex in the Bienaymé forest. The law of the forest with dummy leaves depends on $n$, because the probability of finding a dummy leaf depends on $n$. We then show that the \L ukasiewicz path and height process of the forest with dummy leaves converge under rescaling, jointly with the convergence of the \L ukasiewicz path and height process of the Bienaymé forest under rescaling up to time $\lfloor T n^{2/3}\rfloor$.
    \item We show convergence under rescaling of the out-forest up to time $\lfloor T n^{2/3}\rfloor$ by applying the measure change to the forest with dummy leaves and showing that the resulting forest is a good approximation of the out-forest. 
\end{enumerate}

\subsection{Convergence before adding the dummy leaves}
We define the two processes
$$ \hat{Y}^{\pm}(k)=\sum\limits_{i=1}^k (\hat{D}^{\pm}_{n,i}-1), $$
for $1\leq k\leq n$, which encode the degrees in order of discovery.

We will study these processes via the measure change that we defined in \cref{sec:measure-change}. Let
\begin{equation*}
  Y^{\pm}(k) = \sum_{i=1}^k (Z_i^{\pm} - 1)
\end{equation*}
be the corresponding walks for $(\vZ_i)_{i=1}^{\infty}$. Then, in the critical case, these are related to the centered random walks $V^{\pm}$ by
\begin{equation*}
  Y^+(k) = V^+(k)
  \quad \text{and} \quad
  Y^-(k) = V^-(k) - (\lambda_- - 1) k = V^-(k) - \nu_- k.
\end{equation*}
Therefore, we obtain the following corollary of \cref{prop:measure-change-no-crit}.
\begin{corollary}
  \label{cor:measure-change}
  Suppose we are in the setting of \cref{prop:measure-change-no-crit} and that the criticality condition holds. Then for all $T > 0$,
  \begin{align*}
      &\left( 
          \Phi(n, \floor{n^{2/3} T}),
          \left(
              n^{-2/3} Y^-\left( \floor{n^{2/3} t} \right),
              n^{-1/3} Y^+\left( \floor{n^{2/3} t} \right)
          \right)_{t \in [0, T]}
      \right) \\
      & \hspace{23em} \todist \left( 
          \Phi(T),
          (\nu_- t, \sigma_+ W^+_t)_{t \in [0, T]}
      \right)
  \end{align*}
  in $\R \times \D([0, T], \R^2)$ as $n \to \infty$.
\end{corollary}
Let $(\hat{B}_t,t\geq 0)$ be distributed as follows. For $F$ a suitable test function, and for $(B_t)_{t\geq 0}$ a Brownian motion,
\begin{align*} &\E\left[F(\sigma_+ \hat{B}_t,0\leq t \leq T)\right]\\&=\E\left[\exp\left(-\frac{\sigma_{-+}}{\sigma_+ \mu} \int_0^T s dB_s -\frac{\sigma_{-+}^2 T^3}{6\sigma_+^2 \mu^2}\right)F(\sigma_+ B_t,   0\leq t \leq T)\right].\end{align*}

%  Recall that 
%  $$ {Y}^+(k)=\sum\limits_{i=1}^k (Z^+_i-1),$$ 
%  and 
%  $$ {Y}^-(k)=\sum\limits_{i=1}^k (Z^-_i-1).$$ 
%  Then, Donsker's theorem and the law of large numbers imply the following straightforward lemma.
% \begin{lemma}
% \label{lem.jointconvergenceinout}
%  We have $$ \left(n^{-2/3}{Y}^-\left(\lfloor n^{2/3}t\rfloor\right), t\geq 0\right)
%  \toprob 
%  \left(\nu_-t, t\geq 0 \right)$$
%  in $\D(\R_+,\R)^2$ as $n\to \infty$.  Moreover,
%  \begin{align*} &\left(n^{-1/3}\left({Y}^-\left(\lfloor n^{2/3}t\rfloor\right)-n^{2/3}\nu_-t\right),n^{-1/3}Y^+\left(\lfloor n^{2/3}t\rfloor\right), t\geq 0\right)\\
%  &\todist 
%  \left(\mathbf{B}_t, t\geq 0 \right)\end{align*}
%  in $\D(\R_+,\R)^2$, as $n\to \infty$, with $(\mathbf{B}_t,t\geq 0)$ a Gaussian process with covariance matrix
%  $$\begin{pmatrix} \sigma_-^2  & \sigma_{-+} \\ \sigma_{-+}  & \sigma_+^2  \end{pmatrix}t.$$
% \end{lemma} 

\begin{proposition}\label{prop:convaftermeasurechange} 
 We have that
$$\left(n^{-2/3}\hat{Y}^-\left(\lfloor n^{2/3} t\rfloor\right), n^{-1/3}\hat{Y}^+\left(\lfloor  n^{2/3} t\rfloor\right), 0\leq t \leq T \right) \todist \left( \nu_- t, \sigma_+\hat{B}_t, 0 \leq t \leq T \right)$$
in the Skorokhod topology as $n\to \infty$.
\end{proposition}
\begin{proof}
 We recall from the statement of \cref{cor:measure-change} that $(W^-,W^+)$ is a pair of correlated standard Brownian motions with correlation $\operatorname{Corr}(Z_1^-,Z_1^+)$. 
Let $(B^1_t,t\geq 0)$ and $(B^2_t,t\geq 0)$ be two independent Brownian motions, so that we may define $$(\sigma_-W^-_t,\sigma_+W^+_t,t\geq 0)=\left(\frac{\sigma_{-+}}{\sigma_+}B_t^1+\left(\sigma_-^2-\frac{\sigma_{-+}^2}{\sigma_+^2}\right)^{1/2} B_t^2, \sigma_+ B_t^1, t\geq 0\right).$$ 
 Then, \cref{cor:measure-change} implies that for $F$ a continuous, bounded test function, 
 \begin{align*}&\E\left[F\left(n^{-1/3}\hat{Y}^+\left(\lfloor n^{2/3}t\rfloor\right), 0\leq t \leq T \right) \right]\\
 &=\E\left[F\left(n^{-1/3}\hat{Y}^+\left(\lfloor n^{2/3}t\rfloor\right), 0\leq t \leq T \right)\one_{\lfloor Tn^{2/3}\rfloor \leq R_n} \right]+o(1)
 \\\to& \E\left[\exp\left(-\frac{1}{\mu}\int_0^Tsd\left(\frac{\sigma_{-+}}{\sigma_+}B_s^1+\left(\sigma_-^2-\frac{\sigma_{-+}^2}{\sigma_+^2}\right)^{1/2}B_s^2\right) - \frac{T^3 \sigma_-^2}{6\mu^2}\right) F\left(\sigma_+ B_t^1, 0\leq t \leq T\right)\right]\\
 &=\E\left[\exp\left(-\frac{\sigma_{-+}}{\sigma_+ \mu} \int_0^T s dB^1_s -\frac{\sigma_{-+}^2 T^3}{6\sigma_+^2 \mu^2}\right)F(\sigma_+ B^1_t,   0\leq t \leq T)\right].\end{align*}
 For the details of the argument, we refer the reader to the proof of Theorem 4.1 in \cite{conchon--kerjanStableGraphMetric2021}. Then, the fact that $(Y(k),k\geq 1)$ is a random walk with steps of mean $\nu_-$ implies that
 $$\left(n^{-2/3}Y^-\left(\lfloor n^{2/3}t\rfloor\right),t\geq 0\right)\toprob\left(\nu_- t,t\geq 0\right),$$
 and then, by \cref{cor:measure-change}, also 
 $$\left(n^{-2/3}\hat{Y}^-\left(\lfloor n^{2/3}t\rfloor\right),t\geq 0\right)\toprob\left(\nu_- t,t\geq 0\right),$$
 which proves the statement.
 \end{proof}
 The following proposition characterises the distribution of $(\hat{B}_t, {0\leq t\leq T})$.
 \begin{proposition}
 \label{prop:characterizelimitprocess}
 We have that
 $$(\sigma_+ \hat{B}_t, {0\leq t\leq T})\overset{d}{=}\left(\sigma_+ B_t-\frac{\sigma_{-+}}{2\mu}t^2, {0\leq t\leq T}\right),$$
 where $(B_t)_{t\geq 0}$ is a standard Brownian motion.
 \end{proposition}
 \begin{proof}
 Firstly, we have that for any $t\in [0,T]$ and $\theta>0$,
 \begin{align*} \E\left[\exp(-\theta \sigma_+ \hat{B}_t)\right] &= \E \left[ \exp\left(-\frac{\sigma_{-+}}{\sigma_+ \mu}\int_0^t s dB_s-\frac{\sigma_{-+}^2 t^3}{6\sigma_+^2\mu^2}-\theta\sigma_+ B_t  \right)\right]\\
 &=\E\left[\exp \left( -\frac{\sigma_{-+}}{\sigma_+ \mu}\int_0^t \left(s+\frac{\sigma_+^2\theta \mu}{\sigma_{-+}}\right) dB_s -\frac{\sigma_{-+}^2 t^3}{6\sigma_+^2\mu^2}\right)\right] \\
 &= \exp\left(-\frac{\sigma_{-+}^2}{2\sigma_+^2 \mu^2}\int_0^t \left(s+\frac{\sigma_+^2\theta \mu}{\sigma_{-+}}\right)^2 ds -\frac{\sigma_{-+}^2 t^3}{6\sigma_+^2\mu^2}\right)\\
 &=\exp\left(\frac{\sigma_+^2 t}{2}\theta^2+\frac{\sigma_{-+} t^2}{2\mu}\theta \right)\\
 &= \E \left[\exp\left(-\theta\left(\sigma_+ B_t - \frac{\sigma_{-+}}{2\mu} t^2\right)\right)\right].
 \end{align*}
 Then, more generally, for $m>0$, $0=t_0\leq t_1\leq \cdots \leq t_m=T$, and $\theta_1, \dots, \theta_m \in \R_+$, 
 \begin{align*}
     &\E\left[\exp\left(-\sum_{i=1}^m \theta_i(\sigma_+\hat{B}_{t_i}-\sigma_+\hat{B}_{t_{i-1}})\right)\right]\\
    %  &=\E \left[ \exp\left(-\frac{\sigma}{\mu} \sum_{i=1}^m \int _{t_{i-1}}^{t_i} s dB_s^1-\frac{\sigma_-^2}{2\mu^2}\sum_{i=1}^m \int_{t_{i-1}}^{t_i}s^2ds-\frac{ \sigma_{-+}}{\sigma} \sum_{i=1}^m \theta_i(B_{t_i}^1-B_{t_{i-1}}^1)- 
    %  \left(\sigma_+^2-\frac{\sigma_{-+}^2}{\sigma_-^2}\right)^{1/2}\sum_{i=1}^m \theta_i(B_{t_i}^2-B_{t_{i-1}}^2)\right)\right]\\
     &= \prod_{i=1}^m \E\left[\exp\left(-\frac{\sigma_{-+}}{\sigma_+\mu} \int _{t_{i-1}}^{t_i} s dB_s-\frac{\sigma_{-+}^2 (t_i^3-t_{i-1}^3)}{6\sigma_+^2 \mu^2}-\theta_i \sigma_+  (B_{t_i}-B_{t_{i-1}})\right)\right]\\
     &= \prod_{i=1}^m  \exp\left( -\frac{\sigma_{-+}^2}{2\sigma_+^2 \mu^2}\int_{t_{i-1}}^{t_i}\left(s+\frac{\sigma_+^2\theta_i\mu}{\sigma_{-+}}\right)^2 ds - \frac{\sigma_{-+}^2 (t_i^3-t_{i-1}^3)}{6\sigma_+^2 \mu^2} \right)\\
     &=  \prod_{i=1}^m \exp\left(\frac{ \sigma_+^2 (t_i-t_{i-1}) }{2}\theta_i^2+\frac{\sigma_{-+} (t_i^2-t_{i-1}^2)}{2\mu}\theta_i \right)\\
     &=\E \left[\exp\left(-\sum_{i=1}^m\theta_i\left(\sigma_+ (B_t-B_{t_i}) - \frac{\sigma_{-+}}{2\mu} (t_i^2-t_{i-1}^2)\right)\right)\right],
 \end{align*}
 which proves the result.
 \end{proof}
\subsection{Adding dummy leaves to a Bienaymé forest}\label{subsec.purpleleavesGWforest}
We would like to add dummy leaves to the forest encoded by $(Y^+(l),1\leq l \leq k)$. However, in the absence of a true stack of in-edges, we need to approximate the probability of adding a dummy leaf. We do this by approximating the stack size by its mean $\mu n$. 
% \begin{lemma}
% Consider an eDFS of a configuration model on $n$ vertices, with the total number of in-half-edges equal to $\mu n$. Suppose the number of unpaired in-half-edges of discovered vertices at step $k$ in the exploration is equal to $S_n^{-}(k)$, suppose $(S_n^{+}(l),1\leq l\leq k)$ encodes the \L ukasiewicz path of the out-forest up to time $k$, and set $$I_n^{+}(k)=\inf\left\{S_n^{+}(l),1\leq l\leq k\right\}.$$
% Then, the probability that, in the $(k+1)$th time-step, we sample a surplus edge is given by
% $$p_{k+1}:=\frac{S_n^{-}(k)}{\mu n - k -I^{+}(k)+1}\one_{\{I^{+}(k)=I^{+}(k-1)\}}.$$
% \end{lemma}
% \begin{proof}
% This is a slight adaptation of Lemma \ref{lemma.sampleoutforest}, with the sequence $(\mathbf{\hat{D}}_{n,1},\dots,\mathbf{\hat{D}}_{n,m})$ replaced by $(\mathbf{Z}_1,\dots, \mathbf{Z}_m)$, and the total number of in-edges replaced by its mean $\mu n$.
% \end{proof}
We use this idea to define the forest with dummy leaves and its \L ukasiewicz path $(S_n^{+}(k), k\geq 1)$ as a function of $(Y^-(k), Y^+(k) ,k\geq 1)$ and some extra randomness to decide at which time-steps we add a dummy leaf.
\begin{enumerate} 
    \item Set $P_n(1)=0$, $S_n^{+}(1)=Z_1^+-1$, $S_n^{-}(1)=Z_1^-$. 
    \item Suppose we are given $(P_n(l),S_n^{+}(l),S_n^{-}(l), 1\leq l \leq k)$. Define 
    $I^{+}(k)=\min\{S_n^{+}(l), l\leq k\}$. Then, with probability $$p_{k+1}:=\frac{S_n^{-}(k)}{\mu n - k -I^{+}(k)+1}\one_{\{I^{+}(k)=I^{+}(k-1)\}},$$ independent from everything else, set $P_n(k+1)=P_n(k)+1$. Otherwise, set $P_n(k+1)=P_n(k)$. 
    \item Set $$S_n^{+}(k+1)=Y^+(k+1-P_n(k+1))-P_n(k+1),$$ and $$S_n^{-}(k+1)=Y^-(k+1-P_n(k+1))-P_n(k+1)-I^{+}(k)+1.$$
\end{enumerate}
Let the forest with dummy leaves be the forest with \L ukasiewicz path $(S_n^{+}(k), k\geq 1)$ in which the $k$th vertex is a dummy leaf if and only if $P^n(k)-P^n(k-1)=1$. 
\subsubsection{Convergence of the \L ukasiewicz path}
To show the convergence of the \L ukasiewicz path corresponding to the forest with dummy leaves, we will first examine the limit of $(P_n(k), k\geq 1)$ under rescaling. We will first prove tightness, after which we will show convergence.

\begin{lemma}\label{lemma.tightnesssurplusedges}
 We have that, $$\left(n^{-1/3}P_n\left(\lfloor  n^{2/3}t\rfloor \right) \right)_{n\geq 1}$$ 
 is tight for all $t>0$.
 \end{lemma}
 \begin{proof}
Set $m=\lfloor  n^{2/3}t\rfloor$ and fix $\epsilon>0$. It is trivial that for any $k\leq m$, $$S^{-}(k)\leq \sum_{i=1}^k Z^-_i=Y^-(k)+k.$$ Moreover, $\mu n - k -I^{+}(l)+1>\mu n-k$.  Therefore, $$p_{k+1}\leq \frac{Y^-(k)+k}{\mu n - k}.$$
This upper bound is increasing in $k$. Consequently, conditional on $(Y^+(j),Y^-(j),j\geq 1),$ $n^{-1/3}P_n(m)$ is stochastically dominated by a binomial random variable with parameters  $m$ and $$\frac{Y^-(m)+m}{\mu n - m}\wedge 1.$$
Since $(Y^-(k)+k,k\geq 1)$ is a random walk with steps of finite mean, $\left(n^{-2/3}(Y^-(m)+m)\right)_{n\geq 1}$ is tight. Therefore,
$$\left(n^{1/3} \frac{Y^-(m)+m}{\mu n - m}\right)_{n\geq 1}$$ is tight, which implies that a binomial random variable with parameters  $m$ and $$\frac{Y^-(m)+m}{\mu n - m}\wedge 1$$
is tight. The statement follows.
\end{proof}
\begin{lemma}\label{lemma.convergenceQandP}
We have  
$$\left(n^{-1/3}P_n(\lfloor n^{2/3}t\rfloor), t \geq 0 \right)\toprob \left(\frac{\nu_-}{2\mu} t^2, t\geq 0 \right)$$
in $D(\R_+,\R)$ as $n\to \infty$.

\end{lemma}
\begin{proof}
Recall that
$$p_{k+1}=\frac{S_n^{-}(k)}{\mu n - k -I^{+}(k)+1}\one_{\{I^{+}(k)=I^{+}(k-1)\}}.$$
Define $M^+(k)=\min\{Y^+(l):l\leq k\}$ so that $0\geq I^{+}(k)\geq M^+(k)-P_n(k)$.  Then, by Lemma \ref{lemma.tightnesssurplusedges}, the convergence under rescaling of $Y^+$ shown in \cref{cor:measure-change}, and the continuous mapping theorem, $\left(n^{-1/3}I^+(\lfloor n^{2/3} t \rfloor)\right)_{n\geq 1}$ is tight for all $t\geq 0$.
We will now argue that the indicator, which ensures that the roots are never dummy leaves, does not have an effect on $(P_n(k),k\leq m)$ on the scale of interest. Let $m=\lfloor n^{2/3}t\rfloor$. Define
\begin{align*}\begin{split}
E^p(m)&:=\sum_{k=0}^{m-1}\frac{S_n^{-}(k)}{\mu n - k -I^{+}(k)+1}\one_{\{I^{+}(k)\neq I^{+}(k-1)\}}\\
&\leq -I^{+}(m) \frac{Y^-(m)+m}{\mu n - m},\end{split}\end{align*}
so since $I^{+}(m)$ is of order $n^{1/3}$ and $\frac{Y^{-}(m)+m}{\mu n - m}$ is of order $n^{-1/3}$, $(E^p(m))_{n\geq 1}$ is tight.  This means that if we allow the roots to be dummy leaves, with high probability, we would only sample $O(1)$ roots that are dummy leaves up to time $O(n^{2/3})$. This does not affect $(P_n(k),k\leq m)$ on the scale of interest. \\
 Then, the convergence under rescaling of $Y^-$ and $Y^+$ shown in \cref{cor:measure-change}, the tightness of $\left(n^{-1/3}I^{+}(\lfloor n^{2/3} t \rfloor)\right)_{n\geq 1}$ and Lemma \ref{lemma.tightnesssurplusedges} imply that
\begin{align}\begin{split}\label{eq.convergenceprob}
  &\left(n^{1/3}\frac{S_n^{-}\left(\lfloor n^{2/3} t \rfloor\right)}{\mu n - \lfloor n^{2/3} t \rfloor -I^{p,+}\left(\lfloor n^{2/3} t \rfloor\right)+1},t\geq 0\right)\\
 &=\left(n^{1/3}\frac{Y^-\left(\lfloor n^{2/3} t \rfloor-P_n\left(\lfloor n^{2/3} t \rfloor\right)\right)-P_n\left(\lfloor n^{2/3} t \rfloor\right)-I^{+}\left(\lfloor n^{2/3} t \rfloor\right)+1}{\mu n - \lfloor n^{2/3} t \rfloor -I^{+}\left(\lfloor n^{2/3} t \rfloor\right)+1},t\geq 0\right)\\
 &\toprob \left(\frac{\nu_-}{\mu}t,t\geq 0\right)\end{split}\end{align}
in $D(\R_+,\R)$ as $n\to \infty$. 
Then, by the continuous mapping theorem and the tightness of $(E^p(m))_{n\geq 1}$,
$$\left(n^{-1/3}\sum_{i=0}^{\lfloor n^{2/3}t \rfloor} p_k , t \geq 0\right)\toprob \left(\frac{\nu_-}{2\mu}t^2,t\geq 0\right)$$
in $D(\R_+,\R)$ as $n\to \infty$. \\
Let $\cG=(\cG_k,k\geq 1)$ denote the filtration such that $\cG_{k}$ contains the information on the shape of the forest until time $k$, including which of the first $k$ vertices are dummy vertices. Then, 
$$M_n(k):=\sum_{i=1}^k (\one_{\{P_n(i)-P_n(i-1)=1\}}-p_i)$$ is a $\cG$-martingale. We claim that $(n^{-1/3}M_n(\lfloor n^{2/3} t\rfloor ), t\geq 0)$ converges to $0$ in probability in $D(\R_+,\R)$. Indeed, for any $t\geq 0$,
\begin{align*}\E[n^{-2/3}M_n(\lfloor n^{2/3} t\rfloor )^2]&=n^{-2/3}\sum_{i=1}^{\lfloor n^{2/3} t\rfloor} \E[\E[(\one_{\{P_n(i)-P_n(i-1)=1\}}-p_i)^2|\cG_{i-1}]]\\&=n^{-2/3}\sum_{i=1}^{\lfloor n^{2/3} t\rfloor} \E[p_i-p_i^2]\to 0.\end{align*}
Hence, since for all $t\geq 0$,
\begin{align*}n^{-1/3}P_n(\lfloor n^{2/3}t\rfloor)&=n^{-1/3}\sum_{i=1}^{\lfloor n^{2/3}t\rfloor}  \one_{\{P_n(i)-P_n(i-1)=1\}}\\
&= n^{-1/3}\sum_{i=0}^{\lfloor n^{2/3}t \rfloor} p_k + n^{-1/3} M_n\left(\lfloor n^{2/3}t\rfloor\right),
\end{align*}
we have
$$\left(n^{-1/3}P_n(\lfloor n^{2/3}t\rfloor),t\geq 0\right)\todist  \left( \frac{\nu_-}{2\mu}t^2, t\geq 0 \right),$$
 which proves the statement.

\end{proof}

The convergence of $P_n$ under rescaling implies the convergence of $S^{+}_n$ and $S^{-}_n$ under rescaling, which is the content of the following lemma. Let $(B_t,t\geq 0)$ be a Brownian motion, and define 
$$({B}^{\mathrm{d}}_t,t\geq 0)=\left(B_t-\frac{\nu_-}{2\mu\sigma_+}t^2,t\geq 0\right).$$ 

\begin{lemma}
\label{lem:lukasiewiczpathpurplevertices}
 We have 
 \begin{align*}\left(n^{-1/3}Y^{+}\left(\lfloor n^{2/3}t\rfloor\right), n^{-1/3}S^{+}_n\left(\lfloor n^{2/3}t\rfloor\right), t\geq 0\right)\todist\left(\sigma_+B_t,\sigma_+B^{\mathrm{d}}_t,  t\geq 0\right)\end{align*}
 in $\D(\R_+,\R)^2$  and 
 $$\left(n^{-2/3}S^{-}_n\left(\lfloor n^{2/3}t\rfloor\right),t\geq 0\right)\toprob\left(\nu_- t,t\geq 0\right)$$
 in $\D(\R_+,\R)$ as $n\to\infty$.
\end{lemma}
\begin{proof}
 This follows from the convergence under rescaling of $Y^+$ and $Y^-$ shown in \cref{cor:measure-change} and Lemma \ref{lemma.convergenceQandP}, and the expressions 
 $$S_n^{+}(k+1)=Y^+\left(k+1-P_n(k+1)\right)-P_n(k+1),$$ and $$S_n^{-}(k+1)=Y^-\left(k+1-P_n(k+1)\right)-P_n(k+1)-I^{+}(k)+1.$$
\end{proof}

\subsubsection{Convergence of the height process}\label{subsubsec.convheightprocess}
In this subsection, we will extend \cref{lem:lukasiewiczpathpurplevertices}. We will show that, under rescaling, the height process of the forest with dummy leaves converges jointly with the other encoding processes of the forest with dummy leaves. 
Let $(H^+_n(k),k\geq 1)$ be the height process corresponding to the forest with dummy leaves. Set $$(R^{\mathrm{d}}_t,t\geq 0)=\left({B}^{\mathrm{d}}_t-\inf\left\{{B}^{\mathrm{d}}_s: s\leq t\right\},t\geq 0 \right).$$
\begin{proposition}
\label{prop:convheightprocesspurple}
We have that 
\begin{align*}&\left(n^{-1/3}Y^{+}\left(\lfloor n^{2/3}t\rfloor\right), n^{-1/3}S^{+}_n\left(\lfloor n^{2/3}t\rfloor\right),n^{-1/3}H^{+}_n\left(\lfloor n^{2/3}t\rfloor\right), t\geq 0\right) \\
&\qquad \todist\left(\sigma_+B_t,\sigma_+{B}^{\mathrm{d}}_t, \frac{2}{\sigma_+} R^{\mathrm{d}}_t,  t\geq 0\right)\end{align*}
 in $\D(\R_+,\R)^3$, and 
 $$\left(n^{-2/3}S^{-}_n\left(\lfloor n^{2/3}t\rfloor\right),t\geq 0\right)\toprob\left(\nu_- t,t\geq 0\right)$$
 in $\D(\R_+,\R)$ as $n\to\infty$.
\end{proposition}

\begin{figure}
    \centering
    \includegraphics[scale=1]{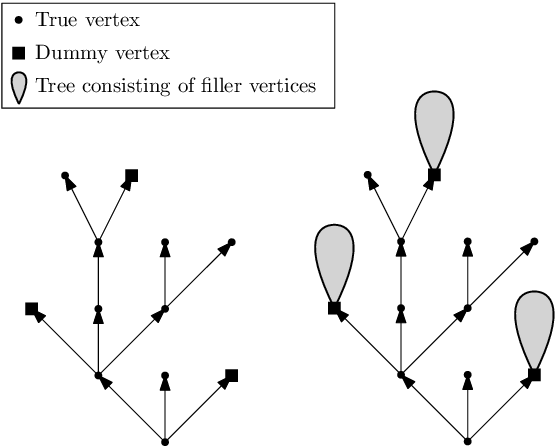}
    \caption{Given a component of the forest with dummy vertices (left), we modify it by sampling independent Bienaymé trees with offspring distributed as $Z^+$ consisting of filler vertices and identifying each dummy leaf with a root of such a tree. The resulting tree (right) is a Bienaymé tree, and the resulting forest is a Bienaymé forest.}
    \label{fig.blackpurpleredforest}
\end{figure}

The difficulty in proving this proposition is the fact that the forest with dummy leaves is not a Bienaymé forest, because the probability of sampling a dummy leaf changes as the exploration is performed. The theory of convergence of height processes under rescaling is well-developed for Bienaymé processes (see e.g. \citet{AST_2002__281__R1_0}), but this is not the case for more general processes.  We will adapt a technique that Broutin, Duquesne and Wang developed in \cite{broutinLimitsMultiplicativeInhomogeneous2021} to show the convergence of the height process of an inhomogeneous random graph under rescaling. The key idea is that  the forest with dummy leaves itself is not a Bienaymé forest, but we can embed it in a Bienaymé forest that does not depend on $n$. We call the extra vertices \emph{filler vertices} and call the resulting forest \emph{the forest with dummy and filler vertices}. We then show convergence under rescaling of the height process corresponding to the forest with dummy and filler vertices, and use this to obtain height process convergence for the forest with dummy leaves. \\
We start by defining the forest with dummy and filler vertices. Informally, we obtain it by modifying the forest with dummy leaves in such a way that a sub-tree consisting of the descendants of a dummy vertex has the same law as a sub-tree consisting of the descendants of a true vertex. We do this by sampling extra Bienaymé trees with offspring distributed as $Z^+$, whose vertices are all filler vertices, and then identifying their roots with the dummy leaves. The resulting forest is a Bienaymé forest containing true, dummy and filler vertices, in which the forest with true vertices and dummy leaves is embedded. This is illustrated in Figure \ref{fig.blackpurpleredforest}. 

The formal procedure is as follows. Suppose we are given $(Y^+(k),S^{+}_n(k),P_n(k),k\geq 1)$, which encodes the forest with dummy leaves.
\begin{enumerate}
    \item Let $(Y^{\mathrm{f}}(k),k\geq 1)$ be an independent copy of $(Y^+(k),k\geq 1)$, which will encode the pendant subtrees that consist of filler vertices.
    \item Define $\theta_n(k)=k-P_n(k-1)+\min\{j: Y^{\mathrm{f}}(j)=-P_n(k-1)\}$. 
    \item Set $\Lambda_n(k)=\max\{j:\theta_n(j)\leq k\}-P_n(\max\{j:\theta_n(j)\leq k\})$. 
    \item We now define \begin{equation}\label{eq.definitionYdf}(Y^{\mathrm{df}}(k),k\geq 1)=(Y^+(\Lambda_n(k))+Y^{\mathrm{f}}(k-\Lambda_n(k)),k\geq 1)\end{equation}
    and we let \emph{the forest with dummy and filler vertices} be the forest with \L ukasiewicz path $(Y^{\mathrm{df}}(k),k\geq 1)$, in which $P_n(\max\{j:\theta_n(j)\leq k\})$ of the first $k$ vertices are dummy vertices, $\Lambda_n(k)$ of the first $k$ vertices are true vertices, and the rest are filler vertices. We let $(H^{\mathrm{df}}(k),k\geq 1)$ be the height process corresponding to the forest with dummy and filler vertices.
\end{enumerate}
By removing the filler vertices from the forest with dummy and filler vertices, we obtain the original forest with dummy leaves. We make the following observations.
\begin{enumerate}
    \item We claim that $\theta_n(k)$ is equal to index in depth first order of the $k$th true or dummy vertex in the forest with dummy and filler vertices. Indeed, note that $\min\{j: Y^{\mathrm{f}}(j)=-P_n(k-1)\}$ is equal to the number of vertices in the first $P_n(k-1)$ trees in the forest encoded by $Y^{\mathrm{f}}$, so that $$\min\{j: Y^{\mathrm{f}}(j)=-P_n(k-1)\}-P_n(k-1)$$ is equal to the number of filler vertices in depth-first order until the $k$th true or dummy vertex. 
    \item Note that $\Lambda_n(k)$ is the number of true vertices amongst the first $k$ vertices. This follows from the fact that $\max\{j:\theta_n(j)\leq k\}$ is the number of true or dummy vertices amongst the first $k$ vertices. 
    \item By the previous remark, $(\Lambda_n(k),k\geq 1)$ only takes steps of size $0$ or $1$. Both $(Y^+(k),k\geq 1)$ and $(Y^{\mathrm{f}}(k),k\geq 1)$ are random walks with steps distributed as $Z^+-1$, so, by construction, $(Y^{\mathrm{df}}(k),k\geq 1)$ is a random walk with steps distributed as $Z^+-1$, so the forest with dummy and filler vertices is a Bienaymé forest with offspring distributed as $Z^+$.
    \item By construction, $(H^{\mathrm{df}}(\theta_n(k)),k\geq 1)$ is the height process corresponding to the forest with dummy vertices. Moreover,
   \begin{equation}\label{eq.constructionSp}(S^{+}_n(k),k\geq 1)=(Y^{\mathrm{df}}(\theta_n(k))-E(\theta_n(k)),k\geq 1),\end{equation}
    where 
    $E(k)$ counts the number of children of the $k$th vertex in the forest with dummy and filler vertices that are filler vertices.
\end{enumerate}
In order to prove \cref{prop:convheightprocesspurple}, considering the construction above and \cref{lem:lukasiewiczpathpurplevertices}, it is sufficient to prove the following lemma.
\begin{lemma}\label{lem:heightprocessblackpurplered}
There exists a process $(D_t,t\geq 0)$ such that 
\begin{align*}
    &\left(n^{-1/3}\left[Y^{\mathrm{df}}\left(\theta_n\left(\lfloor n^{2/3}t\rfloor \right)\right)-E\left(\lfloor n^{2/3}t\rfloor \right)\right], n^{-1/3}H^{\mathrm{df}}\left(\theta_n\left(\lfloor n^{2/3}t\rfloor \right)\right),t\geq 0\right)\\
    &\hspace{15em}\todist\left(\sigma_+D_t,\frac{2}{\sigma_+}\left(D_t-\inf\left\{D_s,s\leq t\right\}\right),t\geq 0\right)
\end{align*}
in $\D(\R_+,\R)^2$ as $n\to \infty$ and $\left(\frac{2}{\sigma_+}\left(D_t-\inf\left\{D_s,s\leq t\right\}\right),t\geq 0\right)$ is the height process corresponding to $(\sigma_+D_t,t\geq 0)$.
\end{lemma} 

The next lemma show that the pathwise construction of $(Y^{\mathrm{df}}(k),H^{\mathrm{df}}(k),k\geq 1)$ converges to its continuous counterpart.

Let $(B_t, t \geq 0)$ and $(B^{\mathrm{f}}_t, t\geq 0)$ be two independent Brownian motions and let $$\theta(t):=t+\inf\left\{s\geq 0 : \sigma_+ B^{\mathrm{f}}_s< -\frac{\nu_-}{2\mu} t^2\right\},$$ and $\Lambda(t)=\inf\{s\geq 0:\theta(s)> t\}$. Define \begin{equation}\label{eq.definitionBpr}\left(B^{\mathrm{df}}_t,t \geq 0\right):=\left( B_{\Lambda(t)}+ B^{\mathrm{f}}_{t-\Lambda(t)}, t\geq 0\right)\end{equation}
and set
$$(R^{\mathrm{df}}_t, t\geq 0):=\left(B^{\mathrm{df}}_t-\inf\{B^{\mathrm{df}}_s,s\leq t\},t\geq 0\right).$$ 

\begin{lemma}\label{lemma.convergenceX+}
We have that $\left((2/\sigma_+)R^{\mathrm{df}}_t, t\geq 0\right)$
is the height process corresponding to $\left(\sigma_+ B^{\mathrm{df}}_t,t \geq 0\right)$.
Moreover,

\begin{equation}\label{eq.convergenceYpr} \left(n^{-1/3}Y^{\mathrm{df}}\left( \lfloor n^{2/3}t \rfloor \right), n^{-1/3}H^{\mathrm{df}}\left( \lfloor n^{2/3}t \rfloor \right),t\geq 0\right)\todist\left( \sigma_+ B^{\mathrm{df}}_{t} ,\frac{2}{\sigma_+}R^{\mathrm{df}}_t, t\geq 0\right)\end{equation}
in $D(\R_+,\R)^2$, jointly with 
$$\left(n^{-1/3}Y^+\left(\lfloor n^{2/3}t \rfloor \right), n^{-1/3}Y^{\mathrm{f}}\left(\lfloor n^{2/3}t \rfloor \right), t\geq 0\right) \todist\left(\sigma_+ B_t,\sigma_+ B^{\mathrm{f}}_t, t\geq 0\right)$$
in $D(\R_+,\R)^2$ and 
$$\left(n^{-2/3} \Lambda_n\left(\lfloor n^{2/3}t\rfloor \right), n^{-2/3}\theta_n\left(\lfloor n^{2/3}t\rfloor \right),t\geq 0\right)\todist\left(\Lambda(t),\theta(t), t\geq 0 \right)$$
in $D(\R_+,\R)^2$ as $n\to \infty$. 
Moreover,   
\begin{equation}\label{eq.convergencecompSprandtheta}\left(n^{-1/3}Y^{\mathrm{df}}\left(\theta_n \left(\lfloor n^{2/3}t\rfloor \right)\right), n^{-1/3}H^{\mathrm{df}}\left(\theta_n\left(\lfloor n^{2/3}t\rfloor \right) \right),t\geq 0 \right) \todist \left(\sigma_+ B^{\mathrm{df}}_{\theta(t)}, \frac{2}{\sigma_+}R^{\mathrm{df}}_{\theta(t)},t\geq 0\right)\end{equation}
in $D(\R_+,\R)^2$ as $n\to \infty$ jointly with the other convergences.
\end{lemma}
In the proof of Lemma \ref{lemma.convergenceX+} we use the following straightforward technical result that follows immediately from the characterization of convergence in the Skorokhod topology given in Ethier and Kurtz \cite[Proposition 3.6.5]{ethierMarkovProcessesCharacterization1986}, .

\begin{lemma}\label{lemma.technicalcomposedfunctions}
If $h_n\to h$ and $f_n\to f$ in $\D(\R_+,\R)$ as $n\to\infty$, and $h_n$ and $h$ are monotone non-decreasing and $h$ is continuous, then 
$$h_n\circ f_n \to h\circ f$$
and 
$$f_n\circ h_n \to f\circ h$$
in $\D(\R_+,\R)$ as $n\to\infty$.
\end{lemma}

We also use the following technical result, that is proved in Appendix \ref{app.technical}.

\begin{lemma}\label{lem.technicalhittingtimes}
If $f_n\to f$ in $\D(\R_+,\R)$ as $n\to\infty$, and $f$ is a continuous function that is not bounded from above, with $f(0)=0$ and  with unique local maxima, then 
$$\left(\inf\{t:f_n(t)>s\},s>0\right)\to \left(\inf\{t:f(t)>s\},s>0\right)$$
in $\D(\R_+,\R)$ as $n\to\infty$.
\end{lemma}

\begin{proof}[Proof of Lemma \ref{lemma.convergenceX+}]
Firstly, note that since $(Y^{\mathrm{df}}(k),k\geq 1)$ encodes a critical Galton-Watson forest with offspring variance $\sigma_+^2$, the proof of Theorem 1.8 in \citet{legallRandomTreesApplications2005} gives us that for $(B^*_s,s\geq 0)$ a Brownian motion,
\begin{align}\label{eq.convergenceX}
&\left(n^{-1/3}Y^{\mathrm{df}}\left(\lfloor n^{2/3}s\rfloor \right), n^{-1/3}H^{\mathrm{df}}\left(\lfloor n^{2/3}s\rfloor \right), s\geq 0 \right) \nonumber \\
& \hspace{15em} \todist \left(\sigma_+ B^*_s,\frac{2}{\sigma_+} \left(B^*_s-\inf\{B^*_u:u\leq s\}\right),  s\geq 0\right)
\end{align} 
 in $D(\R_+,\R)^2$ as $n\to \infty$, and that $\left(\frac{2}{\sigma_+}(B^*_s-\inf\{B^*_u,u\leq s\}),s\geq 0\right)$ is the height process corresponding to $\left(\sigma_+ B^*_s,s \geq 0\right)$. Then, we note that since $(Y^+(k),k\geq 1)\overset{d}{=}(Y^{\mathrm{df}}(k),k\geq 1)$, so that also
 $$\left(n^{-1/3}Y^+\left(\lfloor n^{2/3} t\rfloor \right),t\geq 0\right)\todist\left(\sigma_+ B_t, t\geq 0\right)$$
in $D(\R_+,\R)$ as $n\to \infty$. Then, since also $(Y^+(k),k\geq 1)\overset{d}{=}(Y^{\mathrm{f}}(k),k\geq 1)$ and by Lemma \ref{lem.technicalhittingtimes} and the almost sure uniqueness of the local minima of Brownian, we get that
\begin{align}\begin{split}\label{eq.convergencebychaumont}&\left(n^{-1/3}Y^{\mathrm{f}}\left(\lfloor n^{2/3} s \rfloor \right), n^{-2/3}\inf\left\{k:n^{-1/3}Y^{\mathrm{f}}(k) \leq -x\right\}, s \geq 0, x\geq 0 \right)\\
&\hspace{15em}\todist\left( \sigma_+ B^{\mathrm{f}}_s, \inf\left\{u:\sigma_+ B^{\mathrm{f}}_u < -x\right\}, s\geq 0, x \geq 0\right)\end{split}\end{align}
in $D(\R_+,\R)^2$ as $n\to \infty$.
 
Since $(P_n(k),k\geq 1)$ is non-decreasing, applying Lemma \ref{lemma.technicalcomposedfunctions}, and combining the convergence in \cref{eq.convergencebychaumont} with Lemma \ref{lemma.convergenceQandP} gives that also
$$\left(n^{-2/3}\inf\left\{k:Y^{\mathrm{f}}(k) \leq - P_n\left(\lfloor n^{2/3} t \rfloor -1\right)\right\},t\geq 0\right)\todist\left(\inf\left\{u:\sigma_+ B^{\mathrm{f}}_u< -\frac{\nu_-}{2\mu} t^2\right\},t\geq 0\right)$$
  in $D(\R_+,\R)$ as $n\to \infty$ jointly with the convergence in \cref{eq.convergencebychaumont}. Therefore, 
 \begin{equation}\label{eq.convergencetheta}\left(n^{-2/3}\theta_n\left(\lfloor n^{2/3}t\rfloor \right),t\geq 0 \right) \todist \left(\theta(t),t\geq 0\right)\end{equation}
  in $D(\R_+,\R)$ as $n\to \infty$ jointly with the convergence in \cref{eq.convergencebychaumont}.
Recall that 
$$\Lambda_n(k)=\max\{j:\theta_n(j)\leq k\}-P_n(\max\{j:\theta_n(j)\leq k\}). $$ By definition, for all $n$, $(\theta_n(k),k\geq 1)$ and $(\theta(t),t\geq 0)$ are strictly increasing, so
$$\left(n^{-2/3}\max\{j:\theta_n(j)\leq \lfloor n^{2/3} t \rfloor\} ),t\geq 0\right)\todist\left( \Lambda(t),t\geq0 \right)$$
in $D(\R_+,\R)$ as $n\to \infty$ jointly with the convergence in \cref{eq.convergencebychaumont} and \cref{eq.convergencetheta}. Since $\max\{j:\theta_n(j)\leq \lfloor n^{2/3} t \rfloor\}$ is of order $n^{2/3}$, and, by Lemma \ref{lemma.convergenceQandP}, $P_n(\lfloor n^{2/3}t\rfloor)$ is of order $n^{1/3}$, we get that 
$$\left(n^{-2/3}\Lambda_n\left(\lfloor n^{2/3} t \rfloor\} \right),t\geq 0\right)\todist\left( \Lambda(t),t\geq0 \right)$$
in $D(\R_+,\R)$ as $n\to \infty$ jointly with the convergence in \cref{eq.convergencebychaumont} and \cref{eq.convergencetheta}.\\
To finish the proof, we examine the construction of $(Y^{\mathrm{df}}(k),k\geq 1)$ in \cref{eq.definitionYdf} and the construction of $(B^{\mathrm{df}}_s,s\geq 0)$ in \cref{eq.definitionBpr}. 
Note that $\Lambda_n(k)$ and $k-\Lambda_n(k)$ are non-decreasing. Again, by Lemma \ref{lemma.technicalcomposedfunctions}, this implies that 
$$\left(n^{-1/3}Y^{\mathrm{df}}\left( \lfloor n^{2/3} t \rfloor \right), t\geq 0 \right)\todist \left( B^{\mathrm{df}}_{t}, t\geq 0\right)$$
in $D(\R_+,\R)$ as $n\to \infty$ jointly with all earlier mentioned convergences. Combining this with the convergence in \cref{eq.convergenceX} proves \cref{eq.convergenceYpr}. The fact that $(\theta_n(k),k\geq 1)$ is non-decreasing and Lemma \ref{lemma.technicalcomposedfunctions} then imply \cref{eq.convergencecompSprandtheta}. 
\end{proof}

\begin{lemma}\label{lemma.subtracterrorconverges}
We have that 
\begin{align*}
& \left(n^{-1/3}S^{+}\left(\lfloor n^{2/3}t\rfloor \right)), n^{-1/3}H^{+}\left(\lfloor n^{2/3}t\rfloor \right) ,t\geq 0 \right) \\
& \hspace{15em} \todist \left(\sigma_+ B^{\mathrm{df}}_{\theta (t)},\frac{2}{\sigma_+} \left(B^{\mathrm{df}}_{\theta (t)}-\inf\{B^{\mathrm{df}}_{s}:s\leq \theta(t)\}\right) ,t\geq 0 \right)
\end{align*}
in $\D(\R_+,\R)^2$ as $n\to \infty$. 
\end{lemma}

\begin{proof}
By \cref{eq.constructionSp}, and by Lemma \ref{lemma.convergenceX+}, it is sufficient to show that for any $t>0$,
$$n^{-1/3}\max_{k\leq \lfloor n^{2/3}t\rfloor}E(k)\toprob0.$$
We remind the reader that $E(k)$ counts the number children of the $k$th vertex in the forest with dummy and filler vertices that are filler vertices, so
$$n^{-1/3}\max_{k\leq \lfloor n^{2/3}t\rfloor}E(k)\leq n^{-1/3}\max_{k\leq \theta_n(\lfloor n^{2/3}t\rfloor)}(Y^{\mathrm{f}}(k)-Y^{\mathrm{f}}(k-1)+1),$$
which converges to $0$ by tightness of $\left(n^{-2/3}\theta^{n}(\lfloor n^{2/3}t\rfloor)\right)_{n\geq 1}$ and the fact that $$\left(n^{-1/3}Y^{\mathrm{f}}\left(\lfloor n^{2/3}s\rfloor\right),s\geq 0\right)$$ converges in distribution to a continuous process in $D(\R_+,\R)$ as $n\to\infty$.
\end{proof}

The following lemma is the last ingredient in the proof of \cref{lem:heightprocessblackpurplered}.
\begin{lemma}\label{lemma.heightprocesstimechange}
We have that with probability $1$, $$\left(\frac{2}{\sigma_+} \left(B^{\mathrm{df}}_{\theta (t)}-\inf\{B^{\mathrm{df}}_{s}:s\leq \theta(t)\}\right), t\leq T \right)=\left(\frac{2}{\sigma_+} \left(B^{\mathrm{df}}_{\theta (t)}-\inf\{B^{\mathrm{df}}_{\theta(s)}:s\leq t\}\right), t\leq T \right),$$ which is continuous, and it is the height process corresponding to $\left(\sigma_+ B^{\mathrm{df}}_{\theta (t)},t\leq T\right)$. 
\end{lemma}
\begin{proof}
From \cite{legallRandomTreesApplications2005}, we know that $\left(\frac{2}{\sigma_+}R^{\mathrm{df}}_t,t\geq 0\right)$ is the height process corresponding to $\left(\sigma_+ B^{\mathrm{df}}_{t},t\geq 0\right)$. By definition of the height process, it is sufficient to show that, firstly, with probability $1$, $(B^{\mathrm{df}}_{\theta(t)},t\geq 0)$ is continuous, and, secondly, for all $t\geq 0$, and all $s$ such that $\theta(t-)<s<\theta(t)$, we have $B^{\mathrm{df}}_s > B^{\mathrm{df}}_{\theta(t)}$.

Recall that $(B_t, t \geq 0)$ and $(B^{\mathrm{f}}_t, t\geq 0)$ are two independent Brownian motions, $$\theta(t)=t+\inf\left\{s\geq 0 : \sigma_+ B^{\mathrm{f}}_s< -\frac{\nu_-}{2\mu} t^2\right\},$$ we have $\Lambda(t)=\inf\{s\geq 0:\theta(s)> t\}$, and 
\begin{equation*}
  \left(B^{\mathrm{df}}_t,t \geq 0\right):=\left( B_{\Lambda(t)}+ B^{\mathrm{f}}_{t-\Lambda(t)}, t\geq 0\right).
\end{equation*}

Firstly, note that the jumps of $\theta$ correspond to excursions above the infimum of $B^{\mathrm{f}}$.  With probability $1$, for each of these excursions, the minimum on the excursion is only attained at the endpoints. This can be seen by the almost sure uniqueness of local minima of Brownian motion. We will work on this event of probability 1.

Now fix $t$ such that $\theta(t-)\neq \theta(t)$ and let $s\in (\theta(t-),\theta(t))$. Observe that $\Lambda$ is equal to $t$ on $[\theta(t-),\theta(t)]$. For $[\theta(t-),\theta(t))$ this follows by definition of $\Lambda$, and for $\theta(t)$ it follows since $(\theta(u):u\geq 0)$ is strictly increasing. This implies that
\begin{equation*}
  s-\Lambda(s)<\theta(t)-\Lambda(\theta(t))=\inf\left\{ u\geq 0: \sigma_+ B_u^{\mathrm{f}}<-\frac{\nu_-}{2\mu} t^2\right\}. 
\end{equation*}
By our assumption on the minima of the excursions above the infimum of $B^{\mathrm{f}}$, this implies that
\begin{equation*}
  B^{\mathrm{f}}_{s-\Lambda(s)}>-\frac{\nu_-}{2\mu} t^2=B^{\mathrm{f}}_{\theta(t)-\Lambda(\theta(t))}
\end{equation*}
where the last equality follows from continuity of $B^{\mathrm{f}}$. Combining this with $\Lambda(s)=\Lambda(\theta(t))$ implies that
$B^{\mathrm{df}}_s>B^{\mathrm{df}}_{\theta(t)}$.

Finally, 
$$B^{\mathrm{df}}_{\theta(t-)}=B_{\Lambda(\theta(t-))}+B^{\mathrm{f}}_{\theta(t-)-\Lambda(\theta(t-))}
=B_{t}+B^{\mathrm{f}}_{\theta(t-)-t}$$
and by continuity of $(B^{\mathrm{f}}_s,s\geq 0)$,
\begin{align*}B^{\mathrm{f}}_{\theta(t-)-t}&=B^{\mathrm{f}}\left({\lim_{s\uparrow t}\inf\{u:B^{\mathrm{f}}_u<-\frac{\nu_-}{2\mu} s^2\}}\right)\\&=\lim_{s\uparrow t} B^{\mathrm{f}}\left({\inf\left\{u:B^{\mathrm{f}}_u<-\frac{\nu_-}{2\mu} s^2\right\}}\right)\\&= -\frac{\nu_-}{2\mu^2}t^2\\
&=B^{\mathrm{f}}_{\theta(t)-t}, \end{align*}
so 
$B^{\mathrm{df}}_{\theta(t-)}=B^{\mathrm{df}}_{\theta(t)}.$
\end{proof}

\subsection{Proof of Proposition \ref{prop:convoutforest}}\label{subsubsec.convaftermeasurechange}
We will now combine the convergence of the measure change under rescaling, which is the content of \cref{cor:measure-change}, and the convergence of the encoding processes of the forest with dummy leaves, which is the content of \cref{prop:convheightprocesspurple}, in order to prove \cref{prop:convoutforest}.

\begin{proof}[Proof of \cref{prop:convoutforest}]
Recall that $\hat{P}_n(k)$ denotes the number of dummy leaves amongst the first $k$ vertices in the forest with dummy leaves. Then, as shown in \cref{prop:sampleoutforest}, the probability that the $(k+1)$th vertex in the out-forest is purple, given the degrees in order of discovery and the dummy leaves amongst the first $k$ vertices is equal to
$$q_{k+1}:=\frac{\hat{S}^-_n(k)}{\sum_{i=1}^n D^-_i-k-\hat{I}^+_n(k)}\one_{\left\{\hat{I}^+_n(k-1)= \hat{I}^+_n(k)\right\}},$$
where $\hat{I}^+_n(k)=\min\{\hat{S}^{+}_n(l):l\leq k\}$.
In order to use the results on the forest with dummy leaves, we need to replace the term $\sum_{i=1}^n D^-_i$ in the denominator by $\mu n$. Therefore, define a new forest, \emph{the approximate out-forest}, in which the degrees in order of discovery are the same as in the out-forest. However, in this forest, the probability that the $(k+1)$th vertex is a dummy leaf, given the degrees in order of discovery and the dummy leaves amongst the first $k$ vertices, is equal to
$$\tilde{q}_{k+1}:=\frac{\tilde{S}^-_n(k)}{\mu n-k-\tilde{I}_n(k)}\one_{\left\{\tilde{I}_n(k-1)=\tilde{I}_n(k)\right\}},$$
where $\tilde{S}^-_n(k)$ is the number of unused in-edges of previously discovered vertices in the approximate out-forest up to time $k$ and $-\tilde{I}^+_n(k)$ is the number of components in the approximate out-forest up to time $k$. We let $\tilde{P}_n(k)$ denote the number of dummy leaves amongst the first $k$ vertices in the approximate out-forest. 
We claim that there exists a coupling such that
$$\sum_{i=1}^{\lfloor n^{2/3}T\rfloor }|q_i-\tilde{q}_i|\toprob0$$
as $n\to \infty$. 
Indeed, by the convergence in \cref{prop:convaftermeasurechange}, 
$$\left(n^{-2/3}\sum_{i=1}^{\lfloor n^{2/3}T\rfloor} \hat{D}^n_i\right)_{n>0}$$ is tight. Moreover, with a slight adaptation to the proof of Lemma \ref{lemma.tightnesssurplusedges}, we can show that $\left(n^{-1/3}\tilde{P}_n\left(\lfloor n^{2/3}T\rfloor \right)\right)_{n>0}$ is tight. This, combined with the convergence under rescaling of $(\hat{Y}^+_n(k),k\geq 1)$, implies that also $\left(n^{-1/3}\tilde{I}^+_n\left(\lfloor n^{2/3}T\rfloor \right)\right)_{n>0}$ is tight.  Since $D^-_1,\dots,D^-_n$ are i.i.d.\ random variables with mean $\mu$ and finite variance,
$\left(n^{-1/2}\left(\sum_{i=1}^{n}D^-_i-\mu n\right)\right)_{n>0}$ is tight. By using the trivial identity $a/b-c/d=(b(a-c)-c(d-b))/bd$, this implies that
$\left(n^{2/3}\max_{k\leq \lfloor n^{2/3}T\rfloor }|q_k-\tilde{q}_k|\right)_{n>0}$ is tight, which implies that there exists a coupling such that $\left(\max_{k\leq \lfloor n^{2/3}T\rfloor } |\hat{P}_n(k)-\tilde{P}_n(k)|\right)_{n>1}$ and $\left(\max_{k\leq \lfloor n^{2/3}T\rfloor } |\hat{I}^+_n(k)-\tilde{I}^+_n(k)|\right)_{n>1}$ are tight, which implies that, again by $a/b-c/d=(b(a-c)-c(d-b))/bd$, 
$\left(n^{5/6}\max_{k\leq \lfloor n^{2/3}T\rfloor }|q_k-\tilde{q}_k|\right)_{n>0}$ is tight, which implies that 
$$\sum_{i=1}^{\lfloor n^{2/3}T\rfloor }|q_i-\tilde{q}_i|\toprob0$$
as $n\to \infty$. 
Therefore, under the right coupling, 
$$\P\left(\max_{k\leq \lfloor n^{2/3} T \rfloor}|\hat{P}_n(k)-\tilde{P}_n(k)|>0\right)\to 0.$$
In other words, we can couple the out-forest and the approximate out-forest in such a way that we do not see any difference on the scale of interest. Therefore, we can show convergence under rescaling of the encoding processes of the approximate out-forest instead. To avoid further complicating notation, we will from now on refer to its encoding processes as $$(\hat{S}^{+}_n(k),\hat{H}_n, \hat{S}^-_n(k), \hat{P}_n(k),k\leq \lfloor n^{2/3}T\rfloor).$$ Then, these processes are constructed out of sample paths of $(\hat{Y}^+(k),\hat{Y}^-(k), k\leq \lfloor n^{2/3}T\rfloor )$ and independent randomness in exactly the same way as the sample paths of $$({S}_n^{+}(k),{H}_n^+(k),{S}_n^-(k),P_n(k), k \leq \lfloor n^{2/3}T\rfloor )$$  (corresponding to the forest with dummy vertices) are constructed out of sample paths of $(Y^+(k),Y^-(k), k\leq \lfloor n^{2/3}T\rfloor )$ and independent randomness. 
We will use the following notation:\begin{align*}
    \hat{S}^{+}_{(n)}&:=\left(n^{-1/3}\hat{S}^{+}_n\left(\lfloor n^{2/3} t \rfloor\right),0\leq t \leq T\right)\\
    \hat{H}_{(n)}&:=\left(n^{-1/3}\hat{H}_n\left(\lfloor n^{2/3} t \rfloor\right),0\leq t \leq T\right)\\
    \hat{Y}^+_{(n)}&:=\left(n^{-1/3}\hat{Y}^+\left(\lfloor n^{2/3} t \rfloor\right),0\leq t \leq T\right)\\
     {S}^{+}_{(n)}&:=\left(n^{-1/3}{S}^{+}_n\left(\lfloor n^{2/3} t \rfloor\right),0\leq t \leq T\right)\\
    {H}^+_{(n)}&:=\left(n^{-1/3}{H}^+_n\left(\lfloor n^{2/3} t \rfloor\right),0\leq t \leq T\right)\\
    {Y}^+_{(n)}&:=\left(n^{-1/3}{Y}^+\left(\lfloor n^{2/3} t \rfloor\right),0\leq t \leq T\right).
\end{align*}
Let $f:D([0,T],\R)^3\to \R$ be a bounded, continuous test-function. Then, for $m=\lfloor n^{2/3}T\rfloor$
\begin{align*}\E\left[f\left(\hat{Y}^+_{(n)}, \hat{S}^+_{(n)},  \hat{H}_{(n)}\right) \right]&= \E\left[f\left(\hat{Y}^+_{(n)}, \hat{S}^+_{(n)},  \hat{H}_{(n)}\right) \one_{R_n\geq m} \right]+o(1)\\
&=\E\left[ \E\left[\left. f\left(\hat{Y}^+_{(n)},\hat{S}^+_{(n)},  \hat{H}_{(n)}\right)\right|\hat{\mathbf{D}}_{n,1},\dots, \hat{\mathbf{D}}_{n,m} \right]\one_{R_n\geq m} \right]+o(1)
 \\&=\E\left[ \Phi(n,m)\E\left[\left.f\left(Y^+_{(n)}, S^+_{(n)},  H^+_{(n)}\right)\right| \mathbf{Z}_1,\dots, \mathbf{Z}_m\right]\right]+o(1)\\&=\E\left[ \Phi(n,m)f\left(Y^+_{(n)},S^+_{(n)},  H^+_{(n)}\right)\right]+o(1),
\end{align*}
where we use that $\E\left[\left. f\left(\hat{Y}^+_{(n)},\hat{S}^+_{(n)},  \hat{H}_{(n)}\right)\right| \hat{\mathbf{D}}_{n,1},\dots, \hat{\mathbf{D}}_{n,m}\right]$ and $\one_{R_n\geq m}$ are bounded, adapted functions of $\hat{\mathbf{D}}_{n,1},\dots, \hat{\mathbf{D}}_{n,m}$, and that $\Phi(n,m)$ is the measure change from  $(\hat{\mathbf{Z}}_{1},\dots, \hat{\mathbf{Z}}_{m})$ to $(\hat{\mathbf{D}}_{n,1},\dots, \hat{\mathbf{D}}_{n,m})$. Then, using \cref{cor:measure-change} and \cref{prop:convheightprocesspurple}, following the proof of \cite[Theorem 4.1]{conchon--kerjanStableGraphMetric2021} gives us that 
\begin{align*}
    &\E\left[f\left(\hat{Y}^+_{(n)},\hat{S}^+_{(n)},  \hat{H}_{(n)}\right) \right]\\
    &\to \E\left[\Phi(T)f\left(\sigma_+ B_t,\sigma_+ B^+_t,\frac{2}{\sigma_+}R^+_t,0\leq t \leq T\right)\right].
\end{align*}
Since $$(B^+_t,t\geq 0)=\left(B_t-\frac{\nu_-}{2\sigma_+ \mu}t^2,t\geq 0\right),$$
 \cref{prop:characterizelimitprocess} implies the convergence under rescaling of $(\hat{S}^+_n(k),\hat{H}_n(k),k\geq 0)$. By \cref{prop:convheightprocesspurple}, $S^{-}_n$ converges in distribution under rescaling to a deterministic process, which will not be affected by the measure change. This completes the proof. 
\end{proof}
\subsection{The convergence of the out-forest holds conditionally on the multigraph being simple}
We will now show that the parts of the directed multigraph we observe far beyond the timescale of interest are with high probability simple. We will then use an argument by Joseph \cite{josephComponentSizesCritical2014} to show that this implies that \cref{prop:convoutforest} holds conditional on the resulting multigraph being simple. We let $B_n(k)$ be the number of `bad edges' up to time $k$; to be precise, it equals be the number of self-loops and edges created parallel to an existing edge in the same direction as that edge, up until discovery of the $k$th vertex of the out-forest. Following \cite{conchon--kerjanStableGraphMetric2021}, we call these anomalous edges. 
\begin{proposition}\label{prop.anomalousedges}
Suppose $\beta<1$. Then we have
$$\P\left(B_n(\lfloor n^\beta \rfloor)>0\right)\to 0$$
as $n\to \infty$.
\end{proposition}
\begin{remark}
We adapt the proof of \cite[Lemma 7.1]{josephComponentSizesCritical2014} and of \cite[Proposition 5.3]{conchon--kerjanStableGraphMetric2021} to the directed setting. A significant complication is caused by the conditioning on $$\left\{\sum_{i=1}^n D^-_i=\sum_{i=1}^n D^+_i\right\}.$$ We observe that in both papers, the proof of the aforementioned result is not fully correct, because the authors use the wrong expression for the probability of sampling an anomalous edge. However, the argument below can be adapted to the setting of \cite{josephComponentSizesCritical2014} and \cite{conchon--kerjanStableGraphMetric2021} to yield a correct proof.
\end{remark}
\begin{proof}
We distinguish between the following types of anomalous edges.\\
Self-loops occur when the out-half-edge of a vertex is paired to an in-half-edge of the same vertex.  Let $B^1_n(k)$ be the number of self-loops that are found up to time $k$. For $v$ explored up to time $\lfloor n^\beta\rfloor$, a vertex with in-degree $d^-_v$ and out-degree $d^+_v$, there are $d^-_v d^+_v$ possible combinations of an in-half-edge and an out-half-edge that form a self-loop connected to $v$. Any of these combinations of half-edges is paired with probability bounded above by 
$$\frac{1}{\sum_{i=\lfloor n^\beta \rfloor+1}^n\hat{D}^-_i}.$$
Parallel edges occur when an out-half-edge of a vertex is paired to an in-half-edge of one of its previously explored children. Let $B^2_n(k)$ be the number of parallel edges that are found up to time $k$. For any vertex $v$ with in-degree $d^-_v$, and a parent $p(v)$ with out-degree $d^+_{p(v)}$, there are at most $d^-_v d^+_{p(v)}$ possible combinations of an in-half-edge and an out-half-edge that form a parallel edge from $p(v)$ to $v$. Again, any of these combinations of half-edges is paired with probability bounded above by 
$$\frac{1}{\sum_{i=\lfloor n^\beta \rfloor+1}^n \hat{D}^-_i}.$$
The last type of anomalous edges is a surplus edge with multiplicity greater than 1. Let $B^3_n(k)$ be the number of surplus edges with multiplicity greater than 1 that are found up to time $k$. For a vertex $w$ with out-degree $d^+_w$ and a vertex $v$ with in-degree $d^-_v$, a multiple surplus edge from $w$ to $v$ can only occur if $v$ is discovered before $w$. In that case, there are at most $(d^+_w)^2(d^-_v)^2$ possible pairs of combinations of half-edges, and each of these pairs appears with probability bounded above by
$$\left(\frac{1}{\sum_{i=\lfloor n^\beta \rfloor+1}^n \hat{D}^-_i}\right)^2.$$
Let $p(i)$ denote the index of the parent of the vertex with index $i$. Also, denote $$\cG^n=\sigma\left(\hat{D}^-_1,\hat{D}^+_1,\dots,\hat{D}^-_n,\hat{D}^+_n \right).$$ Then, by the conditional version of Markov's inequality, 

\begin{align*}\P\left(\left.B^1_n(\lfloor n^\beta \rfloor)>0\right| \cG^n \right)&\leq \frac{\sum_{i=1}^{\lfloor n^\beta \rfloor} \hat{D}^-_i\hat{D}^+_i}{\sum_{i=\lfloor n^\beta \rfloor+1}^n \hat{D}^-_i}\wedge 1,\\
\P\left(\left.B^2_n(\lfloor n^\beta \rfloor)>0\right| \cG^n \right)&\leq \frac{\sum_{i=1}^{\lfloor n^\beta \rfloor} \hat{D}^-_i\E\left[\left.\hat{D}^+_{p(i)}\right|\cG^n\right]}{\sum_{i=\lfloor n^\beta \rfloor+1}^n \hat{D}^-_i}\wedge 1,\\
\P\left(\left.B^3_n(\lfloor n^\beta \rfloor)>0\right| \cG^n \right)&\leq \frac{\sum_{i=1}^{\lfloor n^\beta \rfloor}\sum_{j<i} (\hat{D}^+_i)^2 (\hat{D}^-_j)^2 }{\left(\sum_{i=\lfloor n^\beta \rfloor+1}^n \hat{D}^-_i\right)^2 }\wedge 1,\end{align*}
where we note that $p(i)$ is not adapted to $\cG^n$, because ancestral relations in the tree also depend on the surplus edges. However, we observe that by the Cauchy-Schwarz inequality,
\begin{align*}\sum_{i=1}^{\lfloor n^\beta \rfloor} \hat{D}^-_i\E\left[\left.\hat{D}^+_{p(i)}\right|\cG^n\right]&\leq \left(\sum_{i=1}^{\lfloor n^\beta \rfloor} (\hat{D}^-_i)^2\right)^{1/2}\left(\sum_{i=1}^{\lfloor n^\beta \rfloor} \E\left[\left.\hat{D}^+_{p(i)}\right|\cG^n\right]^2\right)^{1/2}\\
&= \left(\sum_{i=1}^{\lfloor n^\beta \rfloor} (\hat{D}^-_i)^2\right)^{1/2}\left(\sum_{j=1}^{\lfloor n^\beta \rfloor} (\hat{D}^+_{j})^2\sum_{i=1}^{\lfloor n^\beta \rfloor}\E\left[\left.\one_{j=p(i)}\right| \cG^n\right]\right)^{1/2}\\ 
&\leq\left(\sum_{i=1}^{\lfloor n^\beta \rfloor} (\hat{D}^-_i)^2\right)^{1/2}\left(\sum_{i=1}^{\lfloor n^\beta \rfloor} (\hat{D}^+_i)^3\right)^{1/2}.\end{align*}
We will show that \begin{equation}\label{eq.conditionalprobanamolousedges}\P\left(\left.B^1_n(\lfloor n^\beta \rfloor)+B^2_n(\lfloor n^\beta \rfloor)+B^3_n(\lfloor n^\beta \rfloor)>0\right| \cG^n \right)\overset{p}{\to}0\end{equation} as $n\to\infty$. We note that 
$$\sum_{i=\lfloor n^\beta \rfloor+1}^n \hat{D}^-_i=\sum_{i=1}^n D^-_i-\sum_{i=1}^{\lfloor n^\beta \rfloor -1}\hat{D}^-_i,$$
and by the weak law of large numbers, $\frac{1}{n}\sum_{i=1}^n D^-_i\overset{p}{\to} \mu n$, so \cref{eq.conditionalprobanamolousedges} follows if we show that 
\begin{enumerate}
    \item $\frac{1}{n}\sum_{i=1}^{\lfloor n^\beta \rfloor} \hat{D}_i^-\overset{p}{\to}0 $, 
    \item $\frac{1}{n}\sum_{i=1}^{\lfloor n^\beta \rfloor}\hat{D}_i^- \hat{D}_i^+\overset{p}{\to}0$,
    \item $\frac{1}{n}\sum_{i=1}^{\lfloor n^\beta \rfloor}(\hat{D}_i^-)^2\overset{p}{\to}0$, and 
    \item $\frac{1}{n}\sum_{i=1}^{\lfloor n^\beta \rfloor}(\hat{D}_i^+)^3\overset{p}{\to}0$
\end{enumerate}
as $n\to \infty$. The proposition will then follow from the bounded convergence theorem.

Note that we can only show the convergence of the Radon-Nikodym derivative $\Phi(n,m)$ under rescaling for $m=O(n^{2/3})$, so it is not straightforward to use the measure change to prove results on the time scale $O(n^\beta)$ for $\beta>2/3$, such as the convergences above. Therefore, instead, we will use \emph{Poissonization} to sample $(\mathbf{\hat{D}}_{n,1},\dots,\mathbf{\hat{D}}_{R_n,n})$. This technique was also used by Joseph in \cite{josephComponentSizesCritical2014}. 

 Let $R_n$ be as before, and, conditional on $R_n$, let $D^{0,+}_1,\dots,D^{0,+}_{n-R_n}$ i.i.d.\ random variables with the law of $D^+$ conditional on the event $\{D^-=0\}$, and set $S_n=\sum_{i=1}^{n-R_n}D^{0,+}_i$. Suppose $R_n=r$ and $S_n=s$. 
Let
$$\pi_0(dt,k_1,k_2)=r\P(D^-=k_1,D^+=k_2|D^->0)k_1\exp(-k_1 t)dt$$
be a measure on $\R_+\times \N^2$, and let $\Pi_0$ be a Poisson point process with intensity measure $\pi_0$ conditional on $\Pi_0(\R,\N,\N)=r$. We view the first coordinate as the time coordinate, and refer to the second and third coordinate as the \emph{point}. Then, the points in $\Pi_0$ ordered by time have the same law as $(\mathbf{\hat{D}}_{n,1},\dots,\mathbf{\hat{D}}_{r,n})$ (before conditioning on the event $\{\sum_{i=1}^nD^-_i=\sum_{i=1}^nD^+_i\}$).  
The intensity of this process is not constant in $t$, so we perform a time change. Define
$$\cL_{\mathbf{D}}(x,y)=\E\left[\left.\exp(-xD^--yD^+)\right|D^->0 \right],$$
and set 
$$\psi(t)=\left(1-\cL_{\mathbf{D}}(\cdot,0)\right)^{-1},$$
so that, by a trivial adaptation of \cite[Lemma 4.1]{josephComponentSizesCritical2014}, for 
$$\pi_r(dt,k_1,k_2):=\P(D^-=k_1,D^+=k_2|D^->0)k_1\exp\left(-k_1 \psi(t/r)\right)\psi'(t/r)dt$$
on $(0,r)\times \N^2$, we have that for $t\in (0,r)$, there exists a probability measure $P_t$ on $\N^2$ such that
$$\pi_r(dt,k_1,k_2)=P_t(D^-=k_1,D^+=k_2)dt.$$
Let ${\Pi}^r$ be a Poisson point process with intensity $\pi_r$. Define $N_r= {\Pi}_r((0,r),\N,\N)$ and $\Delta_r=\int_{(0,r)\times \N^2}(k_1-k_2){\Pi}^r(dt,k_1,k_2)=s$. Then, let ${\Pi}^{r,s}$ have the law of ${\Pi}_r$ conditional on the events $\{N_r=r\}$ and $\{\Delta_r=s\}$. Then, the points of ${\Pi}^{r,s}$ ordered by time are distributed as $(\mathbf{\hat{D}}_{n,1},\dots,\mathbf{\hat{D}}_{n,R_n})$ conditional on the events $\left\{\sum_{i=1}^nD^-_i=\sum_{i=1}^nD^+_i\right\}$, $\{R_n=r\}$ and $\{S_n=s\}$. Let ${\lambda}^{r,s}_t$ be the marginal density of ${\Pi}^{r,s}$ in $t$, so that there exists a probability distribution ${P}^{r,s}_t(k_1,k_2)$ on $\N^2$ such that for ${\pi}^{r,s}_t(k_1,k_2)$ the marginal intensity measure on $\N^2$ of ${\Pi}^{r,s}$ in $t$, 
$${\pi}^{r,s}_t(k_1,k_2)={\lambda}^{r,s}_t{P}^{r,s}_t(k_1,k_2)$$
for all $k_1,k_2\in \N$.

For any $L>0$, define
$$\cE_L=\left\{|R_n-\E[R_n]|\leq Ln^{1/2},  |S_n-\E[S_n]|\leq Ln^{1/2}\right\}.$$
Then, note that 
\begin{align*}\P\left(\frac{1}{n}\sum_{i=1}^{\lfloor n^\beta\rfloor} \hat{D}_i^-\hat{D}_i^+>\epsilon \right)\leq & \P(\cE_L^c) + \P\left(\left.\Pi_{R_n,S_n}\left((0,2n^\beta),\N^2\right)<n^\beta\right|\cE_L \right)\\
&+\P\left(\left.\frac{1}{n}\int_{(0,2n^\beta)\times \N^2}k_1k_2 \Pi_{R_n,S_n}(dt,k_1,k_2)>\epsilon  \right|\cE_L \right)\end{align*}
Fix $\epsilon>0$. By the central limit theorem, we can pick an $L$ such that $\P(\cE_L^c)<\epsilon$ for all $n$. We condition on $\cE_L$. Suppose $R_n=r$ and $S_n=s$. Then, for $P$ a Poisson random variable with rate $2n^\beta$,
 $$\P\left(\Pi_{r,s}\left((0,2n^\beta),\N^2\right)<n^\beta\right)\leq \frac{\P\left(P<n^\beta\right)}{\P(\Delta_r=s, N_r=r)}$$
 We note that the numerator is the probability of a large-deviation event and decreases exponentially fast in $n^\beta$, while the local limit theorem yields that the denominator is of order $n^{-1/2}$ uniformly in all $r$ and $s$ that we consider on $\cE_L$. This implies that $$\P\left(\left.\Pi_{R_n,S_n}\left((0,2n^\beta),\N^2\right)<n^\beta\right|\cE_L\right)\to 0$$
as $n\to \infty$.  
Now, note that for $E^{r,s}_t$ denoting the expectation with respect to $P_t^{r,s}$,
$$\E\left[\frac{1}{n}\int_{(0,2n^\beta)\times \N^2}k_1k_2 \Pi_{r,s}(dt,k_1,k_2) \right]=\frac{1}{n}\int_{(0,2n^\beta)}\lambda^{r,s}_t E^{r,s}_t[D^- D^+]dt,$$
so we start by bounding $E^{r,s}_t[D^- D^+]$. 
We note that
$$E^{r,s}_t\left[{D}^-{D}^+\right]={E}^r_t\left[\left.{D}^-{D}^+\right| \Delta_r=s, N_r=r\right]={E}^r_t\left[{D}^-{D}^+\frac{\P\left[\left. \Delta_n=s, N_r=r \right|\Pi_r(t,D^-,D^+)=1\right]}{\P\left[ \Delta_r=s, N_r=r\right]}\right].$$
By the fact that $\Pi_r$ is a point process, we have that for $k_1$, $k_2$ in $\N$, 
$$\P\left[ \Delta_r=s, N_r=r \left| \Pi_r(t,k_1,k_2)=1\right.\right]=\P\left[ \Delta_r=s+k_2-k_1, N_r=r-1 \right],$$
so that, since $N_r\sim \operatorname{Poisson}(r)$, and since on the event $\{N_r=r-1\}$ (resp. $\{N_r=r\}$),  $\Delta_r-s$ is the sum of $r-1$ (resp. $r$) i.i.d. random variables with finite variance and mean at most $O(n^{-1/2})$, we observe that, by the local limit theorem,
\begin{align*}
    \P\left[ \Delta_r=s, N_r=r \left| \hat{D}^-_t=k_1,\hat{D}^+_t=k_2\right.\right]&=O(n^{-1/2})\text{, and}\\
    \P\left[ \Delta_r=s, N_r=r \right]&=\Theta(n^{-1/2})
\end{align*} for any $k_1$ and $k_2$, and any $r$ and $s$ that we consider on $\cE_L$. Therefore, there exists a $c_1$ such that
$$\frac{\P\left[ \Delta_r=s, N_r=r \left| \hat{D}^-_t=k_1,\hat{D}^+_t=k_2\right.\right]}{\P\left[ \Delta_r=s, N_r=r\right]}<c_1$$
for any $k_1$, $k_2$, $t$ and $n$, and any $r$ and $s$ that we consider on $\cE_L$. If we show that for some $c_2$ $${E}^r_t\left[\hat{D}^-\hat{D}^+\right]<c_2$$ for all $r$ in the interval that we consider and all $t<2n^\beta$, it follows that there is a $c_3$ such that
$${E}^{r,s}_t\left[\hat{D}^-\hat{D}^+\right]<c_3$$ 
for any $k_1$, $k_2$, $t$ and $n$, and any $r$ and $s$ that we consider on $\cE_L$.
We note that by definition of $\pi_r(dt,k_1,k_1)$, 
$${E}^r_t\left[\hat{D}^-\hat{D}^+\right]=\frac{\frac{d^3}{dx^2 dy}\cL_{\mathbf{D}}(x,y)|_{(\psi(t/r),0)}}{\frac{d}{dx}\cL_{\mathbf{D}}(x,y)|_{(\psi(t/r),0)}}.$$
Careful analysis of $\cL_{\mathbf{D}}(x,y)$ and $\psi(s)$ implies that this quantity is bounded uniformly for all $n$, all $r$ in the interval that we consider and all $t\in(0,2n^\beta)$. We refer the reader to the proof of   \cite[Lemma A.1]{josephComponentSizesCritical2014} for the details of a similar argument in the undirected setting.
This implies that 
$$\E\left[\frac{1}{n}\int_{(0,2n^\beta)\times \N^2}k_1k_2 \Pi_{r,s}(dt,k_1,k_2) \right]\leq \frac{C}{n}\E\left[\Pi_{r,s}\left((0,2n^\beta),\N,\N\right)\right].$$
Then, we note that for any $x>0$, for $P$ a Poisson random variable with rate $2n^\beta$,
$$\P\left(\Pi_{r,s}\left((0,2n^\beta),\N,\N\right)>(x+1)2n^\beta\right)\leq \frac{\P\left[P>(x+1)2n^\beta \right]}{\P\left[ \Delta_r=s, N_r=r\right]}.$$
Then, by the local limit theorem and the exponential tail of the Poisson distribution, we obtain that there exist $c_4,c_5>0$ such that for all $n$, all $r$ and $s$ in the interval of interest and all $x>1$,
$$\P\left(\Pi_{r,s}\left((0,2n^\beta),\N,\N\right)>(x+1)2n^\beta\right)\leq c_4\exp(-c_5xn^\beta).$$
This implies that there is a constant $c_6$ such that 
$$\E\left[\Pi_{r,s}\left((0,2n^\beta),\N,\N\right)\right]\leq c_6 n^\beta$$
for all $n$ and all $r$ and $s$ that we consider under $\cE_L$. 
It then follows that 
$$\E\left[\frac{1}{n}\int_{(0,2n^\beta)\times \N^2}k_1k_2 \Pi_{r,s}(dt,k_1,k_2) \right]\to 0$$
as $n\to \infty$ uniformly in all $r$ and $s$ of interest, so for $n$ large enough,
$$\P\left(\left.\frac{1}{n}\int_{(0,2n^\beta)\times \N^2}k_1k_2 \Pi_{R_n,S_n}(dt,k_1,k_2)>\epsilon  \right|\cE_L \right)<\epsilon.$$
This implies that
$$\frac{1}{n}\sum_{i=1}^{\lfloor n^\beta \rfloor}\hat{D}_i^- \hat{D}_i^+\overset{p}{\to}0.$$
The other convergences are proved similarly, and the result follows. 
\end{proof}

\begin{proposition}
 \cref{prop:convoutforest} holds conditionally on the resulting multigraph being simple. 
\end{proposition}
\begin{proof}
Let $\rho(n)=\inf\{k\geq 1:B_n(k)>0\}$, and note that the event that the multigraph formed by the configuration model on $n$ vertices is simple is equal to $\{\rho(n)=\infty\}$. Proposition \ref{prop.anomalousedges} shows that we do not observe any anomalous edges far beyond the timescale in which we explore the largest components of the out-forest. This allows us to conclude that all of the results we prove using the exploration up to time $O(n^{2/3})$ are also true conditioned on $\{\rho(n)=\infty\}$. This follows from the proof of Theorem 3.2 in \cite{josephComponentSizesCritical2014}.
\end{proof}
The results that follow are all obtained by studying the exploration up to time $O(n^{2/3})$, so will also be true conditional on the resulting directed multigraph being simple.

\section{Convergence of the SCCs under rescaling}\label{sec.convSCCs}

In this section, we will use the convergence of the out-forest that we obtained in Section \ref{sec.convoutforest} to show that the SCCs ordered by decreasing number of edges converge under rescaling in the $d_{\vec{\cG}}$-product topology. 

\subsection{Convergence of the out-components that contain an ancestral surplus edge}\label{subsec.ancestral}
In this subsection, we will prove that the out-components that are explored up to time $O(n^{2/3})$ that contain an ancestral surplus edge converge under rescaling. Recall the definition of $(A_n(k),k\geq 1)$ from Subsection \ref{subsubsec.samplecandidates}, and recall that the out-components that contain a non-trivial SCC are the out-components on which $(A_n(k),k\geq 1)$ increases. Moreover, if $(A_n(k),k\geq 1)$ increases on a component, the law of the first increase time corresponds to the position of the tail of the first ancestral surplus edge in the component. \\
We first study the convergence of $(\hat{H}_n^\ell(k),k\geq 1)$ under rescaling. This is an extension of \cref{prop:convoutforest}. Recall that for $(B_t, t\geq 0)$ a standard  Brownian motion, we defined
$$(\hat{B}_t,t\geq 0)=\left( B_t-\frac{\sigma_{-+}+\nu_-}{2\sigma_+ \mu}t^2, t\geq 0\right),$$ 
and its reflected process
$$(\hat{R}_t,t\geq 0)=\left(\hat{B}_t-\inf\left\{\hat{B}_s: s\leq t\right\},t\geq 0\right).$$ 
\begin{proposition}
\label{prop:heightprocesswithlengths}
We have that
\begin{align*}
    &\left(n^{-1/3}\hat{S}^{+}_n\left(\lfloor n^{2/3}t\rfloor \right),n^{-1/3}\hat{H}_{n}\left(\lfloor n^{2/3}t\rfloor \right),n^{-1/3}\hat{H}^\ell_{n}\left(\lfloor n^{2/3}t\rfloor \right),  t\leq T\right)\\
    &\hspace{17em} \todist\left(\sigma_+ \hat{B}_t, \frac{2}{\sigma_+} \hat{R}_t,\frac{2(\sigma_{-+}+\nu_-)}{\sigma_+\mu} \hat{R}_t, t\leq T\right)
\end{align*}
in $\D([0,T],\R)^3$,
jointly with 
$$\left(n^{-2/3}\hat{S}_n^-\left(\lfloor n^{2/3}t\rfloor \right), n^{-1/3}\hat{P}_n\left(\lfloor n^{2/3}t\rfloor \right),t\leq T\right)\overset{p}{\to}\left(\nu_-t,  \frac{\nu_-}{2\mu} t^2, t\leq T\right)$$
in $\D([0,T],\R)^2$ as $n\to \infty$.
\end{proposition}
\begin{proof}
We use \citet[Theorem 1]{deraphelisScalingLimitMultitype2017}, which states the convergence of the height process of a Bienaymé forest with edge-lengths under a few conditions on the degree and edge length distribution. We will apply this result to the Bienaymé forest with dummy and filler vertices, as defined in Subsection \ref{subsubsec.convheightprocess}. 

We equip this forest with edge lengths similarly to how we equipped the out-forest with edge-lengths when we described how to sample the candidates in Subsection \ref{subsubsec.samplecandidates}. We do this as follows. For a dummy or filler vertex with out degree $d^+$, sample its in-degree with law $Z^-$ for $\mathbf{Z}$ conditional on the event $\{Z^+=d^+\}$. The in-degree of the true vertices is encoded by $(Y^{-}(k),k\geq 1)$.  Then, for a vertex with in-degree $d^-$, let the edges connecting it to its children have length $d^--1$ (unless it is the root of the component, then let the edges connecting it to its children have length $d^-$). Let $(H^{\mathrm{df},\ell}(k),k\geq 1)$ be the height process of the resulting forest.\\
We will translate the conditions of Theorem 1 in \cite{deraphelisScalingLimitMultitype2017} to our setting and check them. The conditions are as follows.
\begin{enumerate}
    \item $\E[Z^+]=1$
    \item $1<\E[(Z^+)^2]<\infty$
    \item $\E\left[Z^+\one_{\{Z^->x\}}\right]=o(x^{-2})$ as $x\to \infty$.
\end{enumerate}
Under these conditions, using the notation from Subsection \ref{subsubsec.convheightprocess},
\begin{align}
    \label{eq.convmodifiedheightprocess}
    &\left(n^{-1/3}Y^{\mathrm{df}}\left(\lfloor t n^{2/3}\rfloor \right),n^{-1/3}H^{\mathrm{df}}\left(\lfloor t n^{2/3}\rfloor \right), n^{-1/3}H^{\mathrm{df},\ell}\left(\lfloor t n^{2/3}\rfloor \right),t\geq 0\right) \nonumber \\
    &\hspace{17em} \todist\left(\sigma_+ B_s, \frac{2}{\sigma_+}R_s, \frac{2(\sigma_{+-}+\nu_-)}{\mu\sigma_+}R_s,t\geq 0\right)
\end{align}
in $D(\R_+,\R)^3$ as $n\to \infty$. 
Then, we observe that the rest of the argument in Subsections \ref{subsubsec.convheightprocess} and  \ref{subsubsec.convaftermeasurechange} can be extended to include the height process with edge lengths. This yields the result.\\
Therefore, to finish the proof, we need the conditions of Theorem 1 in \cite{deraphelisScalingLimitMultitype2017} to hold. The conditions are equivalent to 
\begin{enumerate}
    \item $\E[D^+D^-]=\E[D^-]$
    \item $1<\frac{\E[(D^+)^2D^-]}{\E[D^-]}<\infty$
    \item $\E\left[D^+D^-\one_{D^->x}\right]=o(x^{-2})$ as $x\to \infty$. 
\end{enumerate}
Note that the first and second conditions follow directly from the assumptions, and the third condition is implied by $\E[D^+(D^-)^3]<\infty$.
\end{proof}

\begin{proposition}\label{prop.convergenceancestraledges}
We have, jointly with the convergence in \cref{prop:heightprocesswithlengths},
\begin{align*}\left(A_n\left(\lfloor tn^{2/3}\rfloor\right),t\leq T\right)\todist\left(A_t,t\leq T\right),\end{align*}
as $n\to \infty$, where $(A_t,t\geq 0)$ is a Cox process of intensity $$\frac{2(\sigma_{-+}+\nu_-)}{\sigma_+\mu^2} \hat{R}_t$$ at time $t$. The convergence is in $D([0,T],\R)$.
\end{proposition}

\begin{proof}
By definition, $(A_n(k),k\geq 1)$ is a counting process with compensator 
\begin{align*}
    A_n^{comp}(k)&=\sum_{i=1}^k \frac{\hat{H}^\ell(i)}{\hat{S}^-_n(i)}\one_{\{\hat{P}_n(i)-\hat{P}_n(i-1)=1\}}\\
    &=\sum_{j=1}^{\hat{P}_n(k)}\frac{\hat{H}^\ell(\min\{l:\hat{P}_n(l)\geq k\})}{\hat{S}^-_n(\min\{l:\hat{P}_n(l)\geq k\})}.
\end{align*}
 By \citet[Theorem 14.2.VIII]{daleyIntroductionTheoryPoint2008}, the claimed convergence under rescaling of $(A_n(k),k\geq 1)$ follows if we show that 
\begin{equation}\label{eq.convergencecompensator}
    \left(A_n^{comp}\left(\lfloor tn^{2/3}\rfloor \right), t\geq 0\right)\todist\left(\frac{2(\sigma_{-+}+\nu_-)}{\sigma_+\mu^2} \int_0^t\hat{R}_v dv, t \geq 0\right)
\end{equation}
in $\D(\R_+,\R)$ as $n\to \infty$ jointly with the convergence in \cref{prop:heightprocesswithlengths}. Therefore, we will now prove that \cref{eq.convergencecompensator} holds. Since
$$\left(n^{-1/3}\hat{P}_n\left(\lfloor n^{2/3}t\rfloor \right),t\geq 0\right)\overset{p}{\to}\left(\frac{\nu_-}{2\mu}t^2,t\geq 0\right)$$
in $\D(\R_+,\R)$ as $n\to \infty$,
we get that
\begin{align*}\left(n^{-2/3}\min\{l\geq 1:n^{-1/3}\hat{P}_n(l)\geq t\},t\geq 0\right)&\overset{p}{\to}\left(\min\left\{s>0: \frac{\nu_-}{2\mu}s^2>t\right \}, t\geq 0\right)\\
&=:\left(\tau(t),t\geq 0\right) \end{align*}
in $\D(\R_+,\R)$ as $n\to \infty$, because $\left(\frac{\nu_-}{2\mu}t^2,t\geq 0\right)$ is strictly increasing. Then, \cref{prop:heightprocesswithlengths}, Lemma \ref{lemma.technicalcomposedfunctions}, Slutsky's lemma and the continuous mapping theorem imply that 
\begin{align*}&\left(\sum_{j=1}^{\lfloor n^{1/3}t\rfloor}\frac{\hat{H}^\ell(\min\{l:\hat{P}_n(l)\geq k\})}{\hat{S}^-_n(\min\{l:\hat{P}_n(l)\geq k\})},t\geq 0\right)\todist \left( \frac{2(\sigma_{-+}+\nu_-)}{\sigma_+\mu} \int_0^t \frac{\hat{R}_{\tau(s)}}{\nu_- \tau(s)}ds,t\geq 0 \right)
\end{align*}
in $D(\R_+,\R)$ as $n\to\infty$. If we combine this with the convergence under rescaling of $(P_n(k),k\geq 1)$ from Lemma \ref{lemma.convergenceQandP} and apply Lemma \ref{lemma.technicalcomposedfunctions}, some simple analysis then yields \cref{eq.convergencecompensator}, which proves the statement.
\end{proof}

\subsection{Finding the important components in the out-forest}\label{subsec.componentswithancestral}

In this subsection, we will show that, conditional on the convergence under rescaling in Proposition \ref{prop.convergenceancestraledges}, the sequence of intervals that encode the trees with ancestral surplus edges sampled up to time $\lfloor Tn^{2/3}/2\rfloor $ converges as well under rescaling. We want all of the trees that contain such an ancestral surplus edge to be fully explored by time $\lfloor Tn^{2/3}\rfloor$, so we let $T$ be large enough so that this is likely. To be precise, fix $\epsilon>0$ and, from now on, let $T$ be large enough such that $\inf\{\hat{B}_t, t\leq T\}<\inf\{\hat{B}_t, t\leq T/2\}$ with probability at least $1-\epsilon$.

Lemma \ref{lemma.extractexcursions} is a statement about extracting excursion intervals from deterministic functions with marks, which we will apply to the sample paths of $(\hat{S}_n^{+}(k),k\geq 1)$ with the increase times of $(A_n(k),k\geq 1)$ playing the rôle of the marks. The lemma tells us that if the sample paths and increase times converge under rescaling, then the beginnings and endpoints of the excursions above the running infimum that contain the increase times converge as well. 

Let $(f_n(t), t\leq T )$ for $n\geq 1$, and $(f(t),t\leq  T)$ be functions in $\D(\R_+,\R)$, such that 
$$(f_n(t), t\leq T)\to (f(t),t\leq  T)$$ in $\D([0,T],\R)$ as $n\to \infty$. Assume that $(f(t),t\leq T )$ is continuous and that the local minima of $(f(t),t\geq 0)$ are unique. Moreover, let $(x_i^n)_{1\leq i\leq m}$, for $n\geq 1$, and $(x_i)_{1\leq i\leq m}$ be elements of $[0,T]^{m}$ such that for all $i\in [m]$, $x_i^n\to x_i$ in $[0,T]$ as $n\to \infty$, and such that $f(x_i)-\inf\{f(s):s\leq x_i\}>0$ for all $i\in [m]$. Moreover, assume that $\inf\{f(t):t\leq T\}<\inf\{f(t):t\leq x_m\}$ and that $\inf\{f_n(t):t\leq T\}<\inf\{f_n(t):t\leq x^n_m\}$ . For $i \in [m], n \geq 1$, let $g_i^n$ be the left endpoint of the excursion of $f_n$ above its running infimum that contains $x_i^n$, and let $\sigma_i^n$ be the length of this excursion, i.e. 
\begin{align*}
    g_i^n&=\inf\left\{t\geq 0:f_n(t)=\inf\{f_n(s):s\leq x_i^n\}\right\}, \\
    \sigma_i^n&=\inf\left\{ t\geq 0: \inf\{f_n(s):s\leq g_i^n+t\} < \inf\{f_n(s):s\leq x_i^n\}\right\}.
\end{align*}
Similarly, let $g_i$ be the left endpoint of the excursion of $f$ above its running infimum that contains $x_i$, and let $\sigma_i$ be the length of this excursion, i.e. 
\begin{align*}
    g_i&=\inf\left\{t\geq 0:f(t)=\inf\{f(s):s\leq x_i\}\right\}, \\
    \sigma_i&=\inf\left\{ t\geq 0: \inf\{f(s):s\leq g_i+ t\} < \inf\{f(s):s\leq x_i\}\right\}.
\end{align*}
For $S=\{(a_i,b_i), i\in [m]\}$, let $\operatorname{ord}(S)$ be a sequence consisting of the elements of $S$ put in decreasing order of $a_i$, with ties broken arbitrarily, and concatenated with $(0,0)_{i\geq 1}$ so that $\operatorname{ord}(S)\in (\R^2)^\infty$.

\begin{lemma}\label{lemma.extractexcursions}
We have that 
$$\operatorname{ord}\left(\left\{(g_i^n,\sigma_i^n):1\leq i \leq m\right\}\right)\to \operatorname{ord}\left(\left\{(g_i,\sigma_i):1\leq i \leq m\right\}\right)$$
in $(\R^{2})^\infty$ equipped with the product topology as $n\to \infty$. 
\end{lemma}
Note that if a given excursion of $f$ above its running infimum contains multiple marks, only one instance of its left endpoint and excursion length will appear in $\operatorname{ord}\left(\left\{(g_i^n,\sigma_i^n):1\leq i \leq m\right\}\right)$. Therefore, the number of non-zero entries of $\operatorname{ord}\left(\left\{(g_i^n,\sigma_i^n):1\leq i \leq m\right\}\right)$ can vary as $n$ varies, which is why we work in  $(\R^{2})^\infty$. This lemma is proved in Appendix \ref{app.technical}. 

We now apply this result to our process to extract the excursion intervals that contain the marks representing ancestral backedges that are sampled up to time $\lfloor Tn^{2/3}/2\rfloor$.  We recall the following definitions from Subsection  \ref{subsubsec.samplecandidates}. We have that $G_i^n$ is the left endpoint of the excursion of $\hat{S}_n^+$ above its running infimum that encodes the out-component that contains the $i$th ancestral surplus edge, and $\Sigma_i^n$ is the length of this excursion.  Moreover, $G_i$ and $\Sigma_i$ are their continuous counterparts. Formally, for $i\in \left\{1,\dots, A_n\left(\lfloor T n^{2/3}/2\rfloor\right)\right\}$, 
\begin{align*}G_i^n&=\min\left\{k\geq 1:\hat{S}^{+}_n(k)=\min\{\hat{S}^{p,+}_n(l):l\leq X_i^n\}\right\}\text{ and}\\
\Sigma_i^n&=\min\left\{k \geq 1: \min\left\{\hat{S}^{p,+}_n(l):l\leq G_i^n+k\right\} < \min\left\{\hat{S}^{p,+}_n(l):l\leq X_i^n\right\}\right\},
\end{align*}
and for $i\in \left\{1,\dots, A\left(T/2\right)\right\}$, 
\begin{align*}
G_i&=\inf\left\{t\geq 0:\sigma_+\hat{B}_t=\inf\{\sigma_+\hat{B}_s:s\leq X_i\}\right\}\text{ and}\\
\Sigma_i&=\inf\left\{ t\geq 0: \inf\{\sigma_+\hat{B}_s:s\leq G_i+t\} < \inf\{\sigma_+\hat{B}_s:s\leq X_i\}\right\}.
\end{align*}
We recall that the function $\operatorname{ord}$ sorts a set of elements by decreasing second coordinate and appends an infinite sequence of zeroes; the formal definition was given before the statement of Lemma \ref{lemma.extractexcursions}. 
\begin{proposition}\label{prop.extractexcursions}
It holds that
$$\operatorname{ord}\left(\left\{\left(n^{-2/3}G_i^n,n^{-2/3}\Sigma_i^n\right):1\leq i \leq A_n\left(\lfloor T n^{2/3}/2\rfloor\right)\right\}\right)\todist \operatorname{ord}\left(\left\{(G_i,\Sigma_i):1\leq i \leq A(T/2)\right\}\right)$$
in the product topology on $(\R^2)^\infty$ as $n\to \infty$, jointly with the convergence in Proposition \ref{prop.convergenceancestraledges}. 
\end{proposition}
\begin{proof}
 By Skorokhod's representation theorem, we may work on a probability space where the convergence in Proposition \ref{prop.convergenceancestraledges} holds almost surely. We only consider the event on which the convergence holds and $\inf\{\hat{B}_t, t\leq T\}<\inf\{\hat{B}_t, t\leq T/2\}$ holds and claim that we can apply Lemma \ref{lemma.extractexcursions} to the sample paths of $\left(n^{-1/3}\hat{S}^{+}_n\left(\lfloor n^{2/3}t\rfloor\right),t \leq T\right)$ with marks $$\left(n^{-2/3}X_n^i\right)_{1\leq i\leq A_n\left(\lfloor T n^{2/3}/2\rfloor\right)},$$ where we observe that by the convergence, for $n$ large enough, also $$\inf\left\{\hat{S}^{+}_n\left(\lfloor n^{2/3}t\rfloor\right), t\leq T\right\}<\inf\left\{\hat{S}^{+}_n\left(\lfloor n^{2/3}t\rfloor\right), t\leq T/2\right\}$$ holds.  We check the conditions.
Firstly, note that by $A_n\left(\lfloor T n^{2/3}/2\rfloor\right)\to A\left(T/2\right)$ as $n\to \infty$, we can pick $n$ large enough such that $A_n\left(\lfloor T n^{2/3}/2\rfloor\right)=A\left(T/2\right)$. By the local absolute continuity of $(\hat{B}_t,t\geq 0)$ to a Brownian motion, its local minima are almost surely unique. Since
$$\left(A_n\left(\lfloor t n^{2/3}\rfloor\right), t\leq T/2\right) \overset{a.s.}{\to}\left(A\left(t\right),t\leq T/2\right)$$
in $\D(\R_+,\R)$ as $n\to \infty$, we observe that for all $i\in \{1,\dots,A(T/2)\}$, $n^{-2/3}X_i^n\to X_i$ almost surely in $\R$ as $n\to \infty$. The intensity of the Cox process $(A_t,t\geq 0)$ at time $t$ is proportional to $\hat{R}_t$, so $\hat{R}_{X_i}>0$ for all $i$ almost surely. This allows us to apply Lemma \ref{lemma.extractexcursions}, and the convergence follows.
\end{proof}

\subsection{Convergence of the set of candidates}
By \cref{lemma.extractexcursions}, we know that the intervals that encode the out-components that contain an ancestral surplus edge converge under rescaling. This convergence holds jointly with the convergence under rescaling of the first time-step at which an ancestral surplus edge is found in each of these components. We will show that the positions of the other candidates in a component converge as well under rescaling. Recall the procedure to sample candidates that is described in Subsection \ref{subsubsec.samplecandidates}. 
% \begin{lemma}\label{lemma.samplingprocedure}
% Suppose we are exploring the component of $(\hat{F}_n(k),k\geq 1)$ that contains vertex $k$, and suppose $k$ is purple. Denote this component by $\cT$, with root $l$. Moreover, suppose the tails of all important surplus edges in $\cT$ that are discovered up to time $k$ are contained in $C_k\subset\{g+1,\dots,k-1\}$. Then, for $S$ a subset of the vertices of $\cT$, let $\cT(S)$ be the subtree of $\cT$ spanned by $S$. Then, $k$ is the tail of an important surplus edge only if the surplus edge corresponding to $k$ has its head in $\cT(C_k\cup\{g,k
% \})$. 
% \end{lemma}
% \begin{proof}
% This is a direct consequence of Lemma \ref{lemma.whatispartofscc}.\ref{item.factsonsccs2} and \ref{lemma.whatispartofscc}.\ref{item.factsonsccs4}. 
% \end{proof}

We will now show convergence under rescaling of the sequence of candidates in a particular component of $(\hat{F}_n(k),k\geq 1)$ . 

By Skorokhod's representation theorem, we may work on a probability space where the convergence in Propositions \ref{prop.convergenceancestraledges} and \ref{prop.extractexcursions} holds almost surely. Let $(g,\sigma)\in \left\{(G_i,\Sigma_i):i\leq A(T/2)\right\}$, so that, for each $n$ large enough, we can find $(g_n,\sigma_n)\in\left\{(G_i^n,\Sigma_i^n):i\leq A_n\left(\lfloor Tn^{2/3}/2\rfloor\right)\right\}$ such that $(g_n,\sigma_n)\to (g,\sigma)$. Set $V_1=\inf\{t\in [g,g+\sigma]:A(t)=A(g)+1\}$, and similarly, set $V_1^n=\min\{g_n<k\leq g_n+\sigma_n:A_n(k)=A_n(g_n)+1\}$, which are well-defined by definition of $g$, $\sigma$, $g_n$ and $\sigma_n$. By construction, $\{g_n+1,\dots,g_n+\sigma_n\}$ encodes an out-component. Call this component $T^n_{g_n}$. We apply the procedure defined in \cref{prop:samplecandidates} to find the candidates in $T^n_{g_n}$. Let $\mathbf{V}_n(g_n)$ denote the sequence of candidates in $T^n_{g_n}$. Similarly, $[g,g+\sigma]$ encodes a component of the out-$\R$-forest. Call this component $\cT_g$, and apply the procedure in Subsection \ref{subsubsec.samplecontinuousobject} to find the candidates in $\cT_g$. Denote the sequence of candidates by $\mathbf{V}(g)$. 

\begin{proposition}\label{prop.convergencestartingpointscandidates}
Jointly with the convergence in Proposition \ref{prop.extractexcursions}, 
$$n^{-2/3}\mathbf{V}_n(g_n)\todist\mathbf{V}(g)$$
in the product topology.
\end{proposition}
\begin{proof}
We will find a coupling such that $n^{-2/3}\mathbf{V}_n(g_n)\overset{a.s.}{\to}\mathbf{V}(g).$ By the convergence in Propositions \ref{prop.convergenceancestraledges} and \ref{prop.extractexcursions}, $n^{-2/3}V_1^n\overset{a.s.}{\to}V_1$. In general, let $V_m^n$ denote the $m^{th}$ candidate that is found in $T^n_{g_n}$, and let $V_m$ denote the $m^{th}$ candidate that is found in $\cT_g$. Suppose that, for some $m$, we have found a coupling such that 
\begin{equation}\label{eq.convergencemcandidates}n^{-2/3}(V_1^n,\dots,V_m^n)\overset{a.s.}{\to}(V_1,\dots,V_m).\end{equation}
Then, $V_{m+1}^n$ is distributed as the position of the first jump of a counting process $K^n_{m+1}(k)$ on $[0,\infty)$ with compensator 
$$K^n_{comp,m+1}(k)=\sum_{i=V_m^n+1}^k \frac{\ell\left(T^{n,\text{mk}}_{i}\right)-m}{\hat{S}^-(i)}  \one{\left\{P_n(i)=P_n(i-1)+1\right\}}$$
for $k\in [V_m^n+1,g_n+\sigma_n]$ and $0$ otherwise, where $T^{n,\text{mk}}_{i}$ is the subtree of $T^n_{g_n}$ spanned by $\{g_n+1,V^n_1,\dots,V^n_m,i\}$. 
Moreover, for $T_s$ the subtree of $\cT_g$ spanned by $\{g,V_1,\dots,V_m, s\}$, and $|T_s|$ its length as encoded by $\left(\frac{2}{\sigma_+}\hat{R}_t,t\geq 0\right)$, $V_{m+1}$ is the first jump in a counting process $K_{m+1}(t)$ on $[0,\infty)$ with compensator 
$$K_{comp,m+1}(t)= \int_{V_m}^t\frac{\sigma_{-+}+\nu_-}{\mu^2}|T_s|ds$$
for $t\in [V_m,g+\sigma]$ and $0$ otherwise. By the convergence under rescaling of $(\hat{H}^\ell_n(k),k\geq 1)$ in \cref{prop:heightprocesswithlengths}, and by Proposition \ref{prop.extractexcursions}, we get that the metric structure of $T^n_{g_n}$ with distances defined by $(\hat{H}^\ell_n(k),k\geq 1)$, and its projection onto $[n^{-2/3}(g_n+1),n^{-2/3}(g_n+\sigma_n)]$, converge under rescaling to the metric structure of $\cT_g$ with distances defined by $$\left(\frac{2(\sigma_{-+}+\nu_-)}{\sigma_+\mu}\hat{R}_t,t\geq 0\right)$$ and its projection onto $[g,g+\sigma]$. This, combined with \cref{eq.convergencemcandidates} implies that 
$$\left(n^{-1/3}\ell\left(T^{n,\text{mk}}_{\lfloor t n^{2/3}\rfloor}\right),V_m\leq t \leq g+\sigma\right)\overset{a.s.}{\to} \left(\frac{\sigma_{-+}+\nu_-}{\mu^2}|T^{\text{mk}}_t|, V_m\leq t \leq l+\sigma\right)$$ in $\D([V_m,g+\sigma],\R_+)$ as $n\to \infty$. Then, a similar argument to that used in the proof of Proposition \ref{prop.convergenceancestraledges} implies that 
$$\left(K^n_{comp,m+1}\left(\lfloor t n^{2/3}\rfloor \right),V_m\leq t \leq g+\sigma\right)\overset{a.s.}{\to}\left(K_{comp,m+1}(t),V_m\leq t \leq g+\sigma\right),$$
$\D(\R_+,\R_+)$ as $n\to\infty$. This implies that 
$$(K^n_{m+1}(\lfloor t n^{2/3} \rfloor ),t\geq 0)\todist (K_{m+1}(t),t\geq 0)$$ in $\D(\R_+,\R_+)$ as $n\to\infty$ and, in particular, we can find a coupling such that $K_m(\infty)>0$ if and only if $K^n_m(\infty)>0$ for all $n$ large enough, and such that on this event,
$$n^{-2/3}V_{m+1}^n\overset{a.s.}{\to}V_{m+1}.$$
If $K_m(\infty)=0$, set $\mathbf{V}(g)=(V_1,\dots,V_m)$, $\mathbf{V}_n(g_n)=(V^n_1,\dots,V^n_m)$, and the statement follows. If $K_m(\infty)>0$, apply the induction step to $(V_1,\dots,V_{m+1})$ and $(V^n_1,\dots,V^n_{m+1})$. The fact that $|\mathbf{V}(g)|<\infty$ almost surely, as shown in Subsection \ref{subsubsec.samplecontinuousobject}, implies that the induction terminates.
\end{proof}

The following proposition shows that also the law of the heads of the surplus edges corresponding to a candidate converges under rescaling. Moreover, we show convergence under rescaling in the pointed Gromov-Hausdorff topology of an out-component with the location of the candidates and the heads of their corresponding surplus edges. 
\begin{proposition}\label{prop.convergenceheadscandidates}
Suppose the convergence in Propositions \ref{prop.convergenceancestraledges}, \ref{prop.extractexcursions} and \ref{prop.convergencestartingpointscandidates} holds almost surely. Then, for $\mathbf{V}_n(g_n)=(V^n_1,\dots, V^n_{N_n})$, $\mathbf{V}(g)=(V_1,\dots, V_{N})$, let $W^n_i$ be the index of the vertex that the surplus edge corresponding to $V^n_i$ connects to. Similarly, let $W_i$ be the index of the vertex that the surplus edge corresponding to $V_i$ connects to. Then, 
\begin{align*}&\left(n^{-1/3}T^n_{g_n}, n^{-2/3}(g_n+1), \left(n^{-2/3}V^n_1,n^{-2/3}W^n_1\right) \dots, \left(n^{-2/3}V^n_{N_n}, n^{-2/3}W^n_{N_n}\right)\right)\\
&\todist\left(\cT_g, l, (V_1,W_1),\dots, (V_{N},W_{N})\right)\end{align*}
in the $(2N+1)$-pointed Gromov-Hausdorff topology. 
\end{proposition}
\begin{proof}
For $S$ a subset of the vertices of $T^n_{g_n}$, let $T^n_{g_n}(S)$ denote the subtree of $T^n_{g_n}$ spanned by $S$. By definition, for $m\leq N_n$, $W^n_m$ is the vertex corresponding to a uniform unpaired in-half-edge of the vertices in $T^n_{g_n}\left(\{g_n+1,V^n_1,\dots,V^n_{m}\}\right)$. By  \cref{prop:heightprocesswithlengths} and Slutsky's lemma,
$$\left(\frac{\hat{H}_n^\ell\left(\lfloor t n^{2/3}\rfloor \right)}{\hat{H}_n\left(\lfloor t n^{2/3}\rfloor \right)},t\geq 0\right)\overset{a.s.}{\to} \left(\frac{\sigma_{-+}+\nu_-}{2\mu},t\geq 0\right)$$
in $\D(\R_+,\R)$ as $n\to \infty$, which implies that the law of $W^n_m$ converges to the law of a uniform vertex in $T^n_{g_n}\left(\{g_n+1,V^n_1,\dots,V^n_{m}\}\right)$. 
Note that, by Theorem \ref{prop:convoutforest},  Propositions \ref{prop.extractexcursions} and \ref{prop.convergencestartingpointscandidates}, we know that the height process of $T^n_{g_n}$ converges under rescaling to the height process of $\cT_g$, jointly with the convergence under rescaling of the positions of the candidates. By the proof of Proposition 5.4 in \cite{goldschmidtScalingLimitCritical2019}, this implies that
$$\left(n^{-1/3}T^n_{g_n},n^{-2/3}g_n+1,n^{-2/3}V^n_1,\dots,n^{-2/3}V^n_{m}\right)\overset{a.s.}{\to}\left(\cT_g, g,V_1,\dots,V_m\right)$$
in the $(m+1)$-pointed Gromov-Hausdorff topology. Since the relation $$\left|T^n_{g_n}\left(\{g_n+1,V^n_1,\dots,V^n_{m}\}\right)\right|=\left|T^n_{g_n}\left(\{g_n+1,V^n_1,\dots,V^n_{m}, W^n_{m}\}\right)\right|$$ passes to the limit, with $|\cdot|$ denoting the length in the tree as encoded by $(\hat{H}_n(k),k\geq 1)$, the limit in distribution of $n^{-2/3}W^n_m$ is a uniform point on the subtree of $\cT_g$ spanned by $\left(g,V_1,\dots,V_m\right)$, which is equal to the law of $W_m$. This proves the statement.
\end{proof}

The proofs of Propositions \ref{prop.convergencestartingpointscandidates} and \ref{prop.convergenceheadscandidates} imply the following proposition.
\begin{proposition}
By Skorokhod's representation theorem, we may work on a probability space where the convergence in Propositions \ref{prop.convergencestartingpointscandidates} and \ref{prop.convergenceheadscandidates} holds almost surely. Let $T^{n,\text{mk}}(g_n)$ be the subtree of $T^n_{g_n}$ spanned by $\{g_n+1,V^n_1,\dots,V^n_{N_n}\}$, and similarly, let $T^{\text{mk}}(g)$ be the subtree of $\cT_g$ spanned by $\{g,V_1,\dots,V_N\}$. Then, also 
\begin{align*}
    &\left(n^{-1/3}T^{n,\text{mk}}(g_n), n^{-2/3}(g_n+1), \left(n^{-2/3}V^n_1,n^{-2/3}W^n_1\right) \dots, \left(n^{-2/3}V^n_{N_n}, n^{-2/3}W^n_{N_n}\right)\right)\\
    &\to \left(T^{\text{mk}}(g), g, (V_1,W_1),\dots, (V_{N},W_{N})\right)
\end{align*}
almost surely in the $(2N+1)$-pointed Gromov-Hausdorff topology as $n\to \infty$. Also the total length in the trees converges, i.e.
$$n^{-1/3}\left|T^{n,\text{mk}}(g_n)\right|\to \left| T^{\text{mk}}(g)\right|$$
almost surely as $n\to\infty$.
\end{proposition}
We now identify the candidates, as described in Subsection \ref{subsubsec.samplecandidates}. In $T^{n,\text{mk}}(g_n)$, set $V_i^n\sim W_i^n$ for each $1\leq i\leq N_n$, and set $M^n_{g_n}:=T^{n,\text{mk}}(g_n)/\sim$. Moreover, in $T^{\text{mk}}(g)$, set $V_i\sim W_i$ for each $1\leq i\leq N$, and set $\cM_g:=T^{\text{mk}}(g)/\sim$. View both as elements of $\vec{\cG}$ in the natural way. To be precise, in  $M^n_{g_n}$, let the vertex set consist of $g_n+1$, $W_i^n$ for $i\leq N_n$, and the branch points $V_i^n\wedge V_j^n$ for $i\neq j\leq N_n$. Similarly, in $\cM_g$, let the vertex set consist of $g$, $W_i$ for $i\leq N$, and the branch points $V_i\wedge V_j$ for $i\neq j\leq N$. Then we have the following proposition.
\begin{proposition}
On the probability space where the convergence in Propositions \ref{prop.convergencestartingpointscandidates} and \ref{prop.convergenceheadscandidates} holds almost surely, 
$n^{-1/3}M^n_{g_n}\overset{a.s.}{\to} \cM_g$
in $\vec{\cG}$.
\end{proposition}
\begin{proof}
The proof is analogous to the proof of Proposition 5.6 in \cite{goldschmidtScalingLimitCritical2019}.
\end{proof}
\begin{proposition}
\label{prop:sccsinonetreeconverge}
On the probability space where the convergence in Propositions \ref{prop.convergencestartingpointscandidates} and \ref{prop.convergenceheadscandidates} holds almost surely, the SCCs in $n^{-1/3}M^n_{g_n}$, listed in decreasing order of length, converge to the SCCs in $\cM_g$, listed in decreasing order of length, in $\vec{\cG}$ almost surely as $n\to \infty$.
\end{proposition}
\begin{proof}
This follows from Proposition 5.3 in \cite{goldschmidtScalingLimitCritical2019}. This proposition requires that the lengths of the SCCs in $\cM_{g}$ have different lengths almost surely, which is the content of Proposition \ref{prop:allengthsaredifferent}. 
\end{proof}

\begin{proposition}
\label{prop:sccordereduptotimeT}
Let $T>0$, and let $(C^T_i(n),i\geq 1)$ be the kernels of the SCCs that contain a candidate with label at most $\lfloor T n^{2/3}/2\rfloor$, ordered by length. Similarly, let $(\cC^T_i,i\geq 1)$ be the kernels of the SCCs obtained from the out-$\R$-forest with a candidate before time $T/2$, ordered by length. Then,
$$\left(n^{-1/3}C^T_i(n), i\geq 1\right) \todist (\cC^T_i,i\geq 1)$$
in the $\vec{\cG}$-product topology, as $n\to \infty$. 
\end{proposition}
\begin{proof}
This follows from Proposition \ref{prop.extractexcursions}, \cref{prop:sccsinonetreeconverge}, and the fact that all SCCs in the limit object have a different length by \cref{prop:allengthsaredifferent}. 
\end{proof}

Finally, we claim that we can choose $T$ large enough such that the SCCs with the highest number of edges are explored before time $\lfloor Tn^{2/3}\rfloor$. This is the content of the following lemma. The proof is in the same spirit as \citet[Lemma 9]{aldous_1991}.
\begin{lemma}\label{lemma.largesccfoundfirst}
For $\delta>0$ and $I$ an interval, let $SCC(n,I,\delta)$ denote the number of SCCs whose vertices have at total of at least $\delta n^{1/3}$ in-edges (including those which are not part of the SCC) and whose time of first discovery is in $n^{2/3}I$. Then,
$$\lim_{s\to \infty}\limsup_{n} \P\left(SCC(n,(s,\infty),\delta)\geq 1\right)=0\text{ for all }\delta>0.$$
\end{lemma}
\begin{proof}
Fix $\delta>0$. Suppose there is an SCC $C$ with $vn^{1/3}$ total in-edges. Conditionally on this fact, the in-edges that are paired up until the time the first in-edge of $C$ is paired are uniform picks (without replacement) from the total set of in-edges. We use $\Xi_n$ to denote the time of discovery of the first in-edge of $C$ multiplied by $n^{-2/3}$. Then, $\Xi_n\todist\operatorname{Exp}(v)$. Fix $\epsilon>0$. We have that, by the memoryless property at time $s$,
$$\P\left(SCC\left(n,(s,2s),\delta\right)=0|SCC\left(n,(s,\infty),\delta\right)\geq 1\right)$$
is asymptotically bounded from above by 
$\exp(-s\delta)$ by the memoryless property at time $s$. So that we can find an $s>0$ such that for all $n$ large enough,
$$\P\left(SCC\left(n,(s,\infty),\delta\right)\geq 1 \text{ and }SCC\left(n,(s,2s),\delta\right)=0\right)<\epsilon.$$
We claim that, by possibly increasing $s$ and $n$, we also get that 
$$\P\left(SCC\left(n,(s,2s),\delta\right)=0\right)>1-\epsilon,$$
which proves the statement.
Firstly, we observe that the ratio of the length of an $SCC$ and its total in-degree are asymptotically equal to $\frac{\sigma_{-+}+\nu_-}{2\mu}$ by the proof of Proposition \ref{prop.convergenceheadscandidates}. Then, note that it is clear from the description of the limit process that, for $s$ large enough, with probability at most $\epsilon/2$, an SCC with total length at least $\frac{\mu}{\sigma_{-+}+\nu_-}\delta$ is discovered after time $s$. By the convergence of the exploration process on compact time intervals, by choosing $n$ large enough, we can then ensure that 
$$\P\left(SCC\left(n,(s,2s),\delta\right)=0\right)>1-\epsilon.$$
We conclude that 
\begin{equation*}
    \P(SCC\left(n,(s,\infty),\delta\right)\geq 1)\leq 2\epsilon. \qedhere
\end{equation*}
\end{proof}
Note that the number of edges in an SCC is bounded from below by the total number of in-edges of vertices in the SCC.

We now show that for any $j$ and any $\epsilon>0$, we can pick $T$ large enough such the $j$ largest components in $(\cC_i,i\geq 1)$ are contained in $(\cC^T_i,i \geq 1)$ with probability at least $1-\epsilon$.

\begin{lemma}\label{lem:largeSCCfoundfirstcont}
For all $j$ holds that 
$$\lim_{T\to\infty}\P\left(\forall i\leq j, \cC_i\in (\cC^T_i,i \geq 1)\right)=1.$$
\end{lemma}
\begin{proof}
Fix $\epsilon>0$. By \cite[Proposition 5.10]{goldschmidtScalingLimitCritical2019} adapted to our limit object, for $k$ large enough, with probability $1-\epsilon/2$, the $j$ largest components of $(\cC_i,i\geq 1)$ are contained in the $k$ largest components of the out-forest with identifications. Moreover, for $T$ large enough, with probability $1-\epsilon/2$, the $k$ largest excursions above the infimum of a Brownian motion with negative parabolic drift occur before time $T$  (see \cite[Section 3]{aldousBrownianExcursionsCritical1997}). This implies the statement.
\end{proof}
\Cref{thm.main} then follows from \cref{prop:sccordereduptotimeT}, \cref{lemma.largesccfoundfirst} and \cref{lem:largeSCCfoundfirstcont}.

\section{Open problems}\label{subsec.openproblems}
Our work contains the first quantitative results on the directed configuration model at criticality, and is the second metric space convergence result for a directed graph model (after the directed Erd\H{o}s-Rényi graph was studied in \cite{goldschmidtScalingLimitCritical2019}), and many interesting unresolved questions remain.
\begin{enumerate}
    \item The law of our limit object is defined by three parameters that are functions of the (mixed) moments of the degree distribution. Does a different choice of parameters always give a different limit distribution? If so, are the laws absolutely continuous to one another? 
    \item Our methods show that the diameter of the configuration model at criticality is  $\Omega(n^{1/3})$ in probability, which is in contrast with the off-critical cases (for deterministic degrees), in which the diameter is $\Theta(\log(n))$ in probability \cite{caiDiameterDirectedConfiguration2020}. We conjecture that the diameter is in fact $\Theta(n^{1/3})$ in probability. Goldschmidt and Maazoun are working on this question for the directed Erd\H{o}s-Rényi graph at criticality. 
    \item In \cite{goldschmidtScalingLimitCritical2019}, the authors show convergence of the sequence of SCCs in the $\ell_1$-sense, which is stronger than the product topology as considered by us. This for example implies that for the directed Erd\H{o}s-Rényi graph, under rescaling, the total length in the SCCs converges in distribution to some finite random variable. Also for undirected configuration models, there are no results that show metric space convergence in a topology on the sequence of components that is stronger than the product topology \cite{bhamidiUniversalityCriticalHeavytailed2020,conchon--kerjanStableGraphMetric2021,bhamidiGlobalLowerMassbound2020}.
     \item We conjecture that, just like the directed Erd\H{o}s-Rényi graph \cite{goldschmidtScalingLimitCritical2019}, the directed configuration model gives rise to a critical window, that in some sense interpolates between subcritical and supercritical models. It would be interesting to adapt our methods to the critical window.
     \item In future work, we plan to extend our understanding of the SCCs by studying the directed graphs in which they are embedded. A first step would be to study all vertices that can be reached from the non-trivial strongly components. This would illuminate connections between the SCCs and expose the fractal structure of the directed graph, which is not observed when only studying the SCCs themselves.
    \item Another natural next step is to study the model under weaker moment conditions. The first condition to eliminate would be $\E\left[(D^-)^i(D^+)^j\right]<\infty$ for $(i,j)=(1,3)$ and $(i,j)=(3,1)$. Removing the former condition would in some sense make the identifications less uniform on the ancestral lines. To be precise, $(\hat{H}^\ell(k)/\hat{H}(k),k\geq 0)$ will not necessarily converge to a constant process under rescaling of time, which means that the in-edges that can be used to form surplus edges are spread out less uniformly on the out-components. We have reason to believe that this would place the model in  a different universality class, but further research is needed to confirm this. Removing the latter condition requires an adaptation of the proof of Proposition \ref{prop.anomalousedges} that does not use the Cauchy-Schwarz inequality. Also the heavy-tailed case is not well-understood, but given our results, it is natural to expect that a potential limit object would be embedded in a tilted stable tree as defined in \cite{conchon--kerjanStableGraphMetric2021}. Moreover, one could define hybrid models by letting the tail-behaviour of the in- and out-degrees be different. 
    \item We conjecture that the rank-1 inhomogeneous directed random graph model under suitable conditions is part of the same universality class as the directed Erd\H{o}s-Rényi graph \cite{goldschmidtScalingLimitCritical2019} and the model we consider in this work. We believe that our methods and the methods of \cite{goldschmidtScalingLimitCritical2019} can be adapted to obtain a metric space scaling limit for the inhomogeneous directed random graph model, and we intend to pursue this in future work. 

\end{enumerate}
\section*{Acknowledgements}
The authors would like to thank Christina Goldschmidt and Robin Stephenson for many fruitful discussions, and for kindly allowing us to use some of their figures, and Igor Kortchemski for his advice on local limit theorems.

\newpage 

\begin{appendices}
\section{Multivariate triangular local limit theorem}
\label{sec:llt}

The goal of this section is to prove \cref{thm:multi-triangular-llt}. This can be deduced from \citet[Corollary 1]{mukhinLocalLimitTheorems1992}. However, Mukhin's result is more general than is needed to prove \cref{thm:multi-triangular-llt}. As a result, the conditions which we need to check in order to apply Mukhin's result are rather complicated. Instead, we offer here an elementary proof.

First, we recall some definitions. An $\R^d$-valued random variable $\vX$ is lattice if it is non-degenerate and is supported on the translation of some lattice. The symmetrisation of $\vX$ is given by $\vX^* = \vX_1 - \vX_2$ where $\vX_1$ and $\vX_2$ are independent copies of $\vX$. If $\vX$ is lattice, the main lattice of $\vX$ is given by
\begin{equation*}
    \lattice = \bigcup_{m=1}^{\infty} \left\{ 
        \textstyle \sum_{i=1}^m n_i \vx^*_i : \text{$n_i \in \Z$ and $\vx^*_i \in \supp(\vX^*)$ for all $i = 1, \ldots, m$}
    \right\}.
\end{equation*}
Now we restate \cref{thm:multi-triangular-llt}.
\llt*

Before we prove \cref{thm:multi-triangular-llt}, we first prove a sequence of lemmas. Our proof of the local limit theorem will use characteristic functions. Let $\vX$ be $\R^d$-valued. We use the convention that the characteristic function of $\vX$ is given by
\begin{equation*}
    \phi(\vu) = \E \left[ 
        e^{i \vu \cdot \vX}
    \right].
\end{equation*}
The following lemma shows the points at which the characteristic function of a lattice random variables attains 1 in absolute value can be precisely characterised when the main lattice is known. This is an adaptation of \cite[P.67, T1]{spitzerPrinciplesRandomWalk1964}.
\begin{lemma}
    \label{lem:cf-periodicity}
    Suppose $\vX$ is lattice with main lattice $\Z^d$ and characteristic function $\phi$. Then $\abs{\phi(\vu)} = 1$ if and only if $\vu \in (2 \pi \Z)^d$.
\end{lemma}
\begin{proof}
    If every coordinate of $\vu$ is a multiple of $2\pi$, then $\vu \cdot \vX$ has support in $t + 2 \pi \Z$ for some $t \in \R$. Therefore $e^{i\vu \cdot \vX}$ is constant and hence $\abs{\phi(\vu)} = 1$.
    
    For the converse, note the characteristic function of the symmetrisation $\vX^*$ satisfies
    \begin{equation*}
        \E\left[e^{i \vu \cdot \vX^*}\right]
        = \E\left[e^{i \vu \vX_1}\right] \E\left[e^{-i \vu \vX_2}\right]
        = \abs*{\E\left[e^{i \vu \vX}\right]}^2 = 1.
    \end{equation*}
    Thus $e^{i \vu \cdot \vx^*} \in 2 \pi \Z$ for all $\vx^*$ in the support of $\vX^*$. Since the fundamental lattice of $\vX$ is $\Z^d$, there exists $\vx^*_1, \ldots, \vx^*_m$ in the support of $\vX^*$ and $k_1, \ldots, k_m \in \Z$ such that
    \begin{equation*}
        \sum_{i=1}^m k_i \vx^*_i = (1, 0, \ldots, 0).
    \end{equation*}
    Therefore,
    \begin{equation*}
        u^{(1)} = \sum_{i=1}^m k_i \vu \cdot \vx^*_i \in 2 \pi \Z.
    \end{equation*}
    Repeating this argument for the other coordinates of $\vu$ shows all coordinates of $\vu$ are multiples of $2 \pi$.
\end{proof}

The next lemma shows convergence of the means and covariance of $\vX_n$ to that of $\vX$, and moreover shows the uniform integrability condition still holds after centering the random variables.
\begin{lemma}
    \label{lem:l2-cvgc-corrs}
    Suppose conditions (1) and (2) of \cref{thm:multi-triangular-llt} hold. Then, as $n \to \infty$, 
    \begin{equation*}
        \E[\vX_n] \to \E[\vX] \quad \text{and} \quad \cov(\vX_n) \to \cov(\vX).
    \end{equation*}
    Further for each $n$, let $\hat{\vX}_n = \vX_n - \E[\vX_n]$, and $\hat{\vX} = \vX - \E[\vX]$. Then the uniform integrability condition in \cref{eq:ui-condition} holds for the centered random variables $(\hat{\vX}_n)_{n \geq 1}$. This means that
    \begin{equation*}
        \label{eq:ui-mean-center}
        \lim_{L \to \infty} \sup_n \E\left[
            \norm{\hat{\vX}_n}^2
            \one \left\{ \norm{\hat{\vX}_n}^2 > L \right\}
        \right] = 0.
    \end{equation*}
\end{lemma}
\begin{proof}
    By Skorokhod's representation theorem, we can assume without loss of generality that $(\vX_n)_{n \geq 1}$ and $\vX$ are in the same probability space and $\vX_n \to \vX$ almost surely as $n \to \infty$. Then, the condition in \cref{eq:ui-condition} gives uniform integrability of $(\norm{\vX_n}_2^2)_{n \geq 1}$. Thus, by Vitali's convergence theorem, $\vX_n \to \vX$ in $L^2$ as $n \to \infty$. Therefore, $\vX$ has finite second moment and the mean and covariance of $\vX_n$ converge to that of $\vX$.

    Since the means converge, the centerings $\hat{\vX}_n \to \hat{\vX}$ in $L^2$ as $n \to \infty$ also. Thus, $(\norm{\hat{\vX}_n}_2^2)_{n \geq 1}$ is uniformly integrable by the converse statement in Vitali's theorem, as required.
\end{proof}

The following lemma shows that we have a normal central limit theorem.
\begin{lemma}
    \label{lem:clt}
    Suppose we are in the setting of \cref{thm:multi-triangular-llt}. Then
    \begin{equation*}
        \frac{1}{\sqrt{n}} \sum_{i=1}^n (\vX_{n, i} - \E[\vX_n]) \todist N(0, \Sigma)
    \end{equation*}
    as $n \to \infty$.
\end{lemma}
\begin{proof}
    We use the Lindeberg-Feller central limit theorem. We will use the notation $\Sigma = \cov(\vX)$, $\Sigma_n = \cov(\vX_n)$, $\hat{\vX}_{n, i} = \vX_{n, i} - \E[\vX_n]$ and $\hat{\vX}_n = \vX_n - \E[\vX_n]$. We will reduce the problem to the one-dimensional case. By the Cramér--Wold device it is sufficient to show that 
    \begin{equation*}
        \frac{1}{\sqrt{n}} \sum_{i=1}^n \vu \cdot \hat{\vX}_{n, i} \todist
        N(0, \vu \cdot \Sigma \vu)
    \end{equation*}
    for all $\vu \in \R^d$. Define
    \begin{equation*}
        A_{n, i} = \frac{1}{\sqrt{n}} \vu \cdot \hat{\vX}_{n, i}.
    \end{equation*}
    Then by the version of the Lindeberg--Feller central limit theorem stated by Durrett in \cite[P.128-129, Theorem 3.4.10]{durrettProbabilityTheoryExamples2019}, to complete the proof it suffices to check that
    \begin{enumerate}
        \item $\lim_{n \to \infty} \sum_{i=1}^n \E[A_{n, i}^2] = \vu \cdot \Sigma \vu$.
        \item For all $\epsilon > 0$, $\lim_{n \to \infty} \sum_{i=1}^n \E\Big[A_{n, i}^2 \one\left\{ \abs{A_{n, i}} > \epsilon \right\}\Big] = 0$.
    \end{enumerate}
    To check condition (1),
    \begin{align*}
        \lim_{n \to \infty} \sum_{i=1}^n \E[A_{n, i}^2]
        = \lim_{n \to \infty} \E\Big[ (\vu \cdot \hat{\vX}_n)^2 \Big]
        = \lim_{n \to \infty} \vu \cdot \Sigma_n \vu
        = \vu \cdot \Sigma \vu
    \end{align*}
    by \cref{lem:l2-cvgc-corrs}. To check condition (2), for all $\epsilon > 0$
    \begin{align*}
        \lim_{n \to \infty} \sum_{i=1}^n \E\Big[A_{n, i}^2 \one\left\{ \abs{A_{n, i}} > \epsilon \right\}\Big]
        &= \lim_{n \to \infty} \E\Big[ (\vu \cdot \hat{\vX}_n)^2 \one\left\{ (\vu \cdot \hat{\vX}_n)^2 > \epsilon^2 n \right\}\Big] \\
        &\leq \norm{\vu}^2 \lim_{n \to \infty} \E\Bigg[ \norm{\hat{\vX}_n}^2 \one\left\{ \norm{\hat{\vX}_n}^2 > \frac{\epsilon^2}{\norm{\vu}^2} n \right\}\Bigg] \\
        &\leq \norm{\vu}^2 \lim_{n \to \infty} \sup_k \E\Bigg[ \norm{\hat{\vX}_k}^2 \one\left\{ \norm{\hat{\vX}_k}^2 > \frac{\epsilon^2}{\norm{\vu}^2} n \right\}\Bigg] \\
        &= 0
    \end{align*}
    by \cref{eq:ui-mean-center}.
\end{proof}

The last lemma we prove provides bounds on the absolute value of the characteristic functions of $\vX_n$. This will be used to apply the dominated convergence theorem in the main proof.
\begin{lemma}
    \label{lem:dom-cf}
    Suppose we are in the setting of \cref{thm:multi-triangular-llt}. Moreover assume that the common main lattice $\lattice$ is $\Z^d$. Let $\phi_n(\vu)$ be the characteristic function of $\hat{\vX_n} = \vX_n - \E[\vX_n]$. Then there exist $\delta, c> 0$, $\rho \in (0, 1)$ and $N$ such that for all $n \geq N$
    \begin{enumerate}
        \item $\abs{\phi_n(\vu)} \leq 1 - c \norm{\vu}^2$ for all $\vu \in S(\delta)$, and
        \item $\abs{\phi_n(\vu)} \leq \rho$ for all $\vu \in S(\pi) \setminus S(\delta)$
    \end{enumerate}
    where, for all $r > 0$, $S(r) = [-r, r]^d$.
\end{lemma}

\begin{proof}
    Firstly we use a analytical lemma stated by Durrett in \cite[P.116, Lemma 3.3.19]{durrettProbabilityTheoryExamples2019}. By that lemma, there exists a constant $A > 0$ such that 
    \begin{equation*}
        \abs*{e^{ix} - \left(1 + i x - \tfrac{1}{2}x^2\right)} \leq A \min\{\abs{x}, 1\} x^2
    \end{equation*}
    for all $x \in \R$. Then applying this with $x = \vu \cdot (\vX_n - \E[\vX_n])$
    \begin{equation*}
        \abs{\phi_n(\vu)}
        \leq \abs*{1 - \tfrac{1}{2} \vu \cdot \cov(\vX_n) \vu } + R_n(\vu)
    \end{equation*}
    where
    \begin{equation*}
        R_n(\vu) \leq A \E\left[ 
            \min\{\abs{\vu \cdot \hat{\vX}_n}, 1\} (\vu \cdot \hat{\vX}_n)^2
        \right].
    \end{equation*}
    We provide bounds on $R_n$ and $\abs{1 - \frac{1}{2} \vu \cdot \cov(\vX_n) \vu}$, starting with $\abs{1 - \frac{1}{2} \vu \cdot \cov(\vX_n) \vu}$.

    Let $\mineval_n$ and $\maxeval_n$ be the minimum and maximum eigenvalues of $\cov(\vX_n)$ respectively. Then, by standard theory for quadratic forms,
    \begin{equation*}
        \mineval_n \norm{\vu}^2
        \leq \vu \cdot \cov(\vX_n) \vu
        \leq \maxeval_n \norm{\vu}^2.
    \end{equation*}
    Moreover, let $\mineval$ and $\maxeval$ be the minimum and maximum eigenvalues of $\cov(\vX)$ respectively. The eigenvalues of a matrix are continuous in its entries and $\cov(\vX_n) \to \cov(\vX)$ by \cref{lem:l2-cvgc-corrs}. Therefore $\mineval_n \to \mineval$ and $\maxeval_n \to \maxeval$ as $n \to \infty$.

    We have assumed that $\cov(\vX)$ is non-degenerate thus $\mineval > 0$. Hence, there exists $N$ such that for all $n \geq N$,
    \begin{equation*}
        \frac{1}{2} \mineval \leq \mineval_n \leq \maxeval_n \leq 2 \maxeval.
    \end{equation*}
    There also exists $\delta_1 > 0$ sufficiently small that $\maxeval \norm{\vu}^2 < 1$ for all $\vu \in S(\delta_1)$. Then for all $n \geq N$ and $\vu \in S(\delta_1)$,
    \begin{equation}
        \abs*{1 - \tfrac{1}{2} \vu \cdot \cov(\vX_n) \vu }
        = 1 - \tfrac{1}{2} \vu \cdot \cov(\vX_n) \vu
        \leq 1 - \tfrac{1}{4} \mineval \norm{\vu}^2. \label{eq:qf-bound}
    \end{equation}
    To bound $R_n$, by the Cauchy-Schwarz inequality
    \begin{equation*}
        R_n(\vu) \leq A E_n(\vu) \norm{\vu}^2
        \quad \text{where} \quad
        E_n(\vu) = \E[\min\{\norm{\vu} \norm{\centeredX_n}, 1\} \norm{\centeredX_n}^2].
    \end{equation*}
    Then for all $L > 0$, splitting the expectation into the case where $\norm{\centeredX_n}^2 \leq L^2$ and the case when $\norm{\centeredX_n}^2 > L^2$,
    \begin{align*}
        \sup_n E_n(\vu)
        &\leq L^2 \min\{L \norm{\vu}, 1\} +
        \sup_n \E\left[ \norm{\centeredX_n}^2 \one\left\{ \norm{\centeredX_n}^2 > L^2 \right\}
        \right] \\
        &\to \sup_n \E\left[ \norm{\centeredX_n}^2 \one\left\{ \norm{\centeredX_n}^2 > L^2 \right\}
        \right] 
    \end{align*}
    as $\vu \to 0$. This holds for all $L > 0$, hence taking the limit $L \to \infty$ and using \cref{eq:ui-mean-center} we obtain that $\lim_{\vu \to 0} \sup_n E_n(\vu) = 0$. Thus, there exists $\delta_2$ such that for all $\vu \in S(\delta_2)$
    \begin{equation}
        \label{eq:rn-bound}
        R_n(\vu) \leq \frac{1}{8} \mineval \norm{\vu}^2.
    \end{equation}
    Thus setting $\delta = \min\left\{ \delta_1, \delta_2 \right\}$, for all $n \geq N$ and $\vu \in S(\delta)$
    \begin{equation*}
        \abs{\phi_n(\vu)} \leq 1 - c \norm{\vu}^2,
    \end{equation*}
    where $c = \frac{1}{8} \mineval$.

    We now address the second bound. let $\phi$ be the characteristic function of $\vX$. We assume $\vX$ has main lattice $\Z^d$, thus $\abs{\phi(\vu)} = 1$ if and only if every entry of $\vu$ is a multiple of $2 \pi$ by \cref{lem:cf-periodicity}. In particular $\abs{\phi(\vu)} < 1$ for all $\vu \in S(\pi) \setminus S(\delta)$. $\phi$ is continuous and $S(\pi) \setminus S(\delta)$ is compact. Therefore there exists $\epsilon > 0$ such that $\sup_{\vu \in S(\pi) \setminus S(\delta)} \abs{\phi(\vu)} \leq 1 - \epsilon$.

    Since $\vX_n \todist \vX$ as $n \to \infty$, $\phi_n \to \phi$ uniformly on compact sets. Therefore there exists $N$ such that for all $n \geq N$
    \begin{equation*}
        \sup_{\vu \in S(\pi) \setminus S(\delta)} \abs{\phi_n(\vu)} \leq \rho = 1 - \tfrac{1}{2} \epsilon. \qedhere
    \end{equation*}
\end{proof}

We are finally ready to prove \cref{thm:multi-triangular-llt}

\begin{proof}[Proof of \cref{thm:multi-triangular-llt}]
    We first address the case where the main lattice of $\vX$ and all $\vX_n$ is $\Z^d$. The main trick in the proof is to notice that if $n$ is integer valued then
    \begin{equation*}
        \one\{n = 0\} = \frac{1}{2\pi} \int_{-\pi}^{\pi} e^{i n u} \dif u.
    \end{equation*}
    For all $\vy \in \vc_n + \Z^d$, $\sum_{i=1}^n \vX_{n, i} - \vy \in \Z^d$, so 
    \begin{align*}
        \P\left( \sum_{i=1}^n \vX_{n, i} = \vy \right)
        &= \E\left[ 
            \frac{1}{(2 \pi)^d} \int_{S(\pi)} e^{i \vu \cdot \left(\sum_{i=1}^n \vX_{n, i} - \vy\right)} \dif \vu
         \right] \\
        &= \frac{1}{(2 \pi)^d} \int_{S(\pi)} \phi_n(\vu)^n e^{-i \vu \cdot (\vy - n \E[\vX_n])} \dif \vu,
    \end{align*}
    where $\phi_n(\vu) = \E[e^{i \vu \cdot (\vX_n - \E[\vX_n])}]$ and $S(r) = [-r, r]^d$ for all $r > 0$. Recall
    \begin{equation*}
        \vx_n = n^{-1/2}(\vy - n \E[\vX_n]).
    \end{equation*}
    Then, changing variables with $\vs = \sqrt{n} \vu$,
    \begin{equation*}
        n^{d/2} \P\left( \sum_{i=1}^n \vX_{n, i} = \vy \right)
        = \frac{1}{(2 \pi)^d} \int_{S(\pi\sqrt{n})} \phi_n(\vs/\sqrt{n})^n e^{-i \vs \cdot \vx_n} \dif \vs.
    \end{equation*}
    By the Fourier inversion theorem,
    \begin{equation*}
        f(\vx) = \frac{1}{(2 \pi)^d} \int_{\R^d} \psi(\vs) e^{-i \vs \cdot \vx} \dif \vs
    \end{equation*}
    where $\psi$ is the characteristic function of the $N(0, \cov(\vX))$ distribution. Therefore
    \begin{align*}
        &\sup_{\vy \in \vc_n + \lattice} \abs*{
            n^{d/2} \P\left( \textstyle \sum_{i=1}^n \vX_{n, i} = \vy \right) - f(\vx_n(\vy))
            } \nonumber \\
        &\hspace{6em} = \sup_{\vy \in \vc_n + \lattice} \abs*{
            \int_{\R^d} \left(\one_{S(\pi \sqrt{n})}(\vs) \phi_n(\vs/\sqrt{n})^n - \psi(\vs)\right) e^{-i\vs \cdot \vx_n(\vy)} \dif \vs
            } \\
        &\hspace{6em} \leq \int_{\R^d} \abs*{ \one_{S(\pi \sqrt{n})}(\vs) \phi_n(\vs/\sqrt{n})^n - \psi(\vs)} \dif \vs.
    \end{align*}
    We apply the dominated convergence theorem. To dominate the integrand, first note that $\psi$ is integrable. Secondly let $\delta$, $c$, $\rho$ and $N$ be as in \cref{lem:dom-cf}. For all $n \geq N$ and for all $\vs \in S(\delta \sqrt{n})$, 
    \begin{equation*}
        \abs{\phi_n(\vs/\sqrt{n})}^n
        \leq (1 - c \norm{s}^2/n)^n
        \leq e^{-c \norm{s}^2}.
    \end{equation*}
    Let $C = - \log(\rho)$. Note if $\vs \in S(\pi \sqrt{n})$ then $\norm{\vs}^2 \leq \pi^2 d n$. Thus for all $n \geq N$ and $\vs \in S(\pi \sqrt{n}) \setminus S(\delta \sqrt{n})$
    \begin{equation*}
        \abs{\phi_n(\vs/\sqrt{n})}^n
        \leq e^{-Cn} 
        \leq e^{- \frac{C}{\pi^2 d} \norm{\vs}^2}.
    \end{equation*}
    Hence for all $n \geq N$,
    \begin{equation*}
        \abs*{ \one_{S(\pi \sqrt{n})}(\vs) \phi_n(\vs/\sqrt{n})^n - \psi(\vs)}
        \leq e^{-c \norm{\vs}^2} + e^{- \frac{C}{\pi^2 d} \norm{\vs}^2} + \abs{\psi(\vs)}
    \end{equation*}
    where, in particular, the right hand side is integrable. By \cref{lem:clt},
    \begin{equation*}
        \phi_n(\vs/\sqrt{n})^n \to \psi(\vs)
    \end{equation*}
    as $n \to \infty$ for all $\vs \in \R^d$. Thus for all $\vs \in \R^d$
    \begin{equation*}
        \one_{S(\pi\sqrt{n})}(\vs) \phi(\vs/\sqrt{n})^n \to \psi(\vs)
    \end{equation*}
    as $n \to \infty$. Hence by the dominated convergence theorem
    \begin{equation*}
        \lim_{n \to \infty} \sup_{\vy \in \vc_n + \lattice} \abs*{
            n^{d/2} \P\left( \textstyle \sum_{i=1}^n \vX_{n, i} = \vy \right) - f(\vx_n)
            } = 0,
    \end{equation*}
    as required.

    Finally we generalise to any main lattice $\lattice$. Suppose that $\lattice$ is generated by the columns of the invertible matrix $A$. Then $A$, viewed as a linear transform, is a isomorphism mapping $\Z^d$ to $\lattice$. Thus $A^{-1}\vX_n$ and $A^{-1}\vX$ will have common lattice $\Z^d$ for all $n$. Moreover we can check the remaining assumptions of \cref{thm:multi-triangular-llt} still hold, thus uniformly for $\vy$ in the translation of $\lattice$ containing the support of $\sum_{i=1}^n \vX_{n, i}$,
    \begin{equation*}
        \P\Bigg(\sum_{i=1}^n A^{-1} \vX_{n, i} = A^{-1}\vy\Bigg)
        = \frac{1}{\sqrt{(2 \pi n)^{d} \det\tilde{\Sigma}}} \exp \left( 
            -\frac{1}{2} (A^{-1} \vx_n)^T \tilde{\Sigma}^{-1} (A^{-1} \vx_n)
        \right) + \littleo(n^{-d/2}).
    \end{equation*}
    where $\tilde{\Sigma} = \cov(A^{-1} \vX)$. This simplifies to
    \begin{equation*}
        \P\Bigg(\sum_{i=1}^n \vX_{n, i} = \vy\Bigg)
        = \frac{1}{\sqrt{(2 \pi n)^{d} \det\tilde{\Sigma}}} \exp \left( 
            -\frac{1}{2} \vx_n^T (A \tilde{\Sigma} A^T)^{-1} \vx_n 
        \right) + \littleo(n^{-d/2}).
    \end{equation*}
    We have that
    \begin{equation*}
        \tilde{\Sigma}
        = \cov(A^{-1} \vX)
        = A^{-1} \cov(\vX) (A^{-1})^T.
    \end{equation*}
    Therefore
    \begin{equation*}
        \det(\tilde{\Sigma}) = \det(A)^{-2} \det(\cov(\vX)) = \det(\lattice)^{-2} \det(\cov(\vX))
    \end{equation*}
    and so
    \begin{align*}
        \P\Bigg(\sum_{i=1}^n \vX_{n, i} = \vy\Bigg)
        &= \frac{\det(\lattice)}{\sqrt{(2 \pi n)^{d} \det(\cov \vX)}} \exp \left( 
            -\frac{1}{2} \vx_n^T \cov(\vX)^{-1} \vx_n 
         \right) + \littleo(n^{-d/2}),
    \end{align*}
    as required.
\end{proof}
\section{Proof of technical lemmas}\label{app.technical}
\begin{proof}[Proof of Lemma \ref{lem.technicalhittingtimes}]
Denote $g_n(s)=\inf\{t:f_n(t)>s\}$ and $g(s)=\inf\{t:f(t)>s\}$. By Proposition 3.6.5 in the book by Ethier and Kurtz \cite{ethierMarkovProcessesCharacterization1986}, it is sufficient to show that for any $s>0$, for any $s_n\to s$, 
\begin{enumerate}
    \item $\max\{|g_n(s_n)-g(s)|,|g_n(s_n)-g(s-)|\}\to 0$;
    \item If $u_n\leq s_n$ for all $n$, $s_n\to s$, $u_n\to s$ and $g_n(s_n)\to g(s-)$, then $g_n(u_n)\to g(s-)$;
    \item If $u_n\geq s_n$ for all $n$, $s_n\to s$, $u_n\to s$ and $g_n(s_n)\to g(s)$, then $g_n(u_n)\to g(s)$. 
\end{enumerate}
Fix $s>0$. If $g(s-)=g(s)$, the result is straightforward, so we focus on $g(s-)<g(s)$.  

We start by proving the first property. Fix $\epsilon>0$ and suppose $s_n\to s$. We observe that $g(s-)<g(s)$ implies that $f$ has a local maximum at $g(s-)$ and that $f(g(s-))=f(g(s))=s$. By the uniqueness of local maxima of $f$ and the definition of $g$, there exists a $\delta_1>0$ such that for all $t<g(s-)-\epsilon$, we have that $f(t)<s-\delta_1$. Similarly, there exists a $\delta_2>0$ such that for all $g(s-)+\epsilon<t<g(s)-\epsilon$, we have that $f(t)<s-\delta_2$. Moreover, define $$\delta_3=\sup\left\{f(t):g(s)<t<g(s)+\epsilon\right\}-s,$$ so that, by definition of $g$, we have that $\delta_3>0$. Define $\delta=\min\{\delta_1,\delta_2,\delta_3\}$. 
Now, let $n$ be large enough such that $\sup_{t\in [0,g(s)+\epsilon]}|f_n(s)-f(s)|<\delta/2$ and $|s_n-s|<\delta/2$. 
Then, it holds that 
\begin{enumerate}
    \item $f_n(t)<s-\delta/2<s_n$ for all $t<g(s-)-\epsilon$;
    \item $f_n(t)<s-\delta<s_n$ for all $g(s-)+\epsilon<t<g(s)-\epsilon$;
    \item There is a $g(s)<t<g(s)+\epsilon$ such that $f_n(t)>s+\delta/2>s_n$.
\end{enumerate}
These tree facts imply that $g_n(s_n)\subset [g(s-)-\epsilon,g(s-)+\epsilon]\cup[g(s)-\epsilon,g(s)+\epsilon]$, which proves the first of the three conditions.

Then, the second and third property follow immediately from the first property and the monotonicity of $g_n$ and $g$. \end{proof}

\begin{proof}[Proof of Lemma \ref{lemma.extractexcursions}]
First, note that $g_i^n$, $\sigma_i^n$, $g_i$, and $\sigma_i$ are well-defined for all $i\in [m]$, $n\geq 1$ by $\inf\{f(t):t\leq T\}<\inf\{f(t):t\leq x_m\}$ and $\inf\{f_n(t):t\leq T\}<\inf\{f_n(t):t\leq x^n_m\}$.

Fix $i$. We will first show that $g^n_i\to g_i$ and $\sigma_i^n\to \sigma_i$ as $n\to \infty$. Firstly, note that by the assumption that $f(x_i)-\inf\{f(s):s\leq x_i\}>0$ and the continuity of $f$, $g_i<x_i<g_i+\sigma_i$. Fix $0<\epsilon<\min\{x_i-g_i,g_i+\sigma_i-x_i\}/2$. We claim that the following conditions are sufficient for $g^n_i\to g_i$ and $\sigma_i^n\to \sigma_i$ as $n\to \infty$. For all $n$ large enough,
\begin{enumerate}
    \item \label{cond.excursions1} $g_i+\epsilon<x^n_i<g_i+\sigma_i-\epsilon$
    \item \label{cond.excursions2}$\inf\left\{f_n(s):s\in (g_i-\epsilon, g_i+\epsilon)\right\}<\inf\left\{f_n(s):s\in [g_i+\epsilon,g_i+\sigma_i-\epsilon] \right\}$, 
    \item \label{cond.excursions3}$\inf\left\{f_n(s):s\in (g_i-\epsilon, g_i+\epsilon)\right\}<\inf\left\{f_n(s):s\in [0,g_i-\epsilon]\right\}$,
    \item \label{cond.excursions4} $\inf\left\{ f_n(s):s\in (g_i+\sigma_i-\epsilon,g_i+\sigma_i+\epsilon)\right\}<\inf\left\{f_n(s):s\in [0,g_i+\sigma_i-\epsilon]\right\}$
\end{enumerate}
Indeed, conditions \ref{cond.excursions1}, \ref{cond.excursions2} and \ref{cond.excursions3} imply $|g^n_i-g_i|<\epsilon$, while conditions \ref{cond.excursions1}, \ref{cond.excursions2} and \ref{cond.excursions4} imply $|(g^n_i+\sigma^n_i)-(g_i+\sigma_i)|<\epsilon$. Note that condition \ref{cond.excursions1} holds for $n$ large enough by definition of $\epsilon$ and the convergence of $x_i^n$ to $x_i$. To show the other conditions, define
\begin{align*}\delta_1&=\inf\left\{f(s):s\in [g_i+\epsilon,g_i+\sigma_i-\epsilon]\right\}-\inf\left\{f(s):s\in (g_i-\epsilon,g_i+\epsilon)\right\}\\
\delta_2&=\inf\left\{f(s):s\in [0,g_i-\epsilon]\right\}-\inf\left\{f(s):s\in (g_i-\epsilon,g_i+\epsilon)\right\}\\
\delta_3&=\inf\left\{f(s):s\in [0,g_i+\sigma_i-\epsilon]\right\}-\inf\left\{f(s):s\in (g_i+\sigma_i-\epsilon,g_i+\sigma_i+\epsilon)\right\}.
\end{align*}
By uniqueness of local minima and the definition of $g_i$ and $\sigma_i$, we have $\delta:=\min\{\delta_1,\delta_2,\delta_3\}/3>0$. Then, note that for $n$ large enough, $\sup\{|f_n(s)-f(s)|:s\leq g_i+\epsilon\}<\delta$, which implies conditions \ref{cond.excursions2}, \ref{cond.excursions3}, and \ref{cond.excursions4} for such $n$. \\
Since $i$ was arbitrary, and $m$ is finite, we find that $$(g_i^n,\sigma_i^n)_{1\leq i\leq m}\to (g_i,\sigma_i)_{1\leq i\leq m}$$
in $\R^{2m}$ as $n\to \infty$. \\
We now claim that $g_i^n\to g_i$ and $g_j^n\to g_i$ implies that $g_i^n=g_j^n$ for $n$ large enough. Indeed, by definition of $g_i^n$, $g_j^n$ and $\sigma_i^n$, we have that $g_i^n<g_j^n$ implies that $g_j^n-g_i^n\geq \sigma_i^n$, and by the argument above, $\sigma_i^n\to \sigma_i>0$, so $g_i^n-g_j^n\to 0$ can only hold if $g_i^n=g_j^n$ for $n$ large enough. This implies that 
$$\#\left\{(g_i^n,\sigma_i^n):1\leq i \leq m\right\}\to \#\left\{(g_i,\sigma_i):1\leq i \leq m\right\}.$$
Then, the result follows.
\end{proof}
\end{appendices}

\bibliography{Content/Support/Bibliography.bib}

\end{document}